\documentclass{ut-thesis}

\usepackage{amsmath}
\usepackage{amsthm}
\usepackage{amssymb}
\usepackage[all]{xy}

\usepackage[ps2pdf,bookmarks=true]{hyperref}

\degree{Doctor of Philosophy}
\department{Mathematics}
\gradyear{2003}
\author{Peter G. Bubenik}
\title{{Cell Attachments and the Homology of Loop Spaces and
Differential Graded Algebras}}

\CompileMatrices

\newtheorem{thm}{Theorem}
\newtheorem{lemma}[thm]{Lemma}
\newtheorem{prop}[thm]{Proposition}

\newtheorem{cor}[thm]{Corollary}
\newtheorem{conj}[thm]{Conjecture}
\newtheorem{thma}{Theorem}

\theoremstyle{definition}
\newtheorem{eg}[thm]{Example}

\newtheorem{rem}[thm]{Remark}

\numberwithin{thm}{chapter}
\numberwithin{equation}{chapter}

\newcommand {\Fp} {\ensuremath {\mathbb{F}_p} }
\newcommand {\F} {\ensuremath {\mathbb{F}} }
\newcommand {\freeT} {\ensuremath {\mathbb{T}} }

\newcommand {\Z} {\ensuremath {\mathbb{Z}} }
\newcommand {\Q} {\ensuremath {\mathbb{Q}} }
\newcommand {\R} {\ensuremath {\mathbb{R}} }
\newcommand {\AH}[1] {\ensuremath {\mathbf{A} (#1)} }

\newcommand {\lX} {\ensuremath {\Omega X} }
\newcommand {\lY} {\ensuremath {\Omega Y} }
\newcommand {\lZ} {\ensuremath {\Omega Z} }
\newcommand {\LL}[1] {\ensuremath {\mathbf{L} (#1)} }

\newcommand {\freeL} {\ensuremath {\mathbb{L}} }
\newcommand {\p} {\ensuremath {^{\prime}} }
\newcommand {\pp} {\ensuremath {^{\prime\prime}} }
\newcommand {\isom} {\ensuremath {\cong} }
\newcommand {\tensor} {\ensuremath {\otimes} }
\newcommand {\incl} {\ensuremath {\hookrightarrow} }

\newcommand {\sdp} {\ensuremath {\ltimes} }
\newcommand {\PiS} {\ensuremath {\prod \mathcal{S}}}
\newcommand {\bL} {\ensuremath {\mathbf{L}} }
\newcommand {\bA} {\ensuremath {\mathbf{A}} }
\newcommand {\bB} {\ensuremath {\mathbf{B}} }
\newcommand {\nP} {\ensuremath {\tilde{P}} }
\newcommand {\fpnP} {\ensuremath {\forall p \in \nP } }

\newcommand {\fanypnP} {for any $p \in \nP$}
\newcommand {\feachpnP} {for each $p \in \nP$}
\newcommand {\onto} {\ensuremath {\twoheadrightarrow} }
\newcommand {\eL} {\ensuremath {\mathbf{\underline{L}}} }

\newcommand {\isomto} {\ensuremath {\xrightarrow{\isom}} }
\newcommand {\xto}[1] {\ensuremath {\xrightarrow{#1}} }
\newcommand {\eolBox} {\vspace{2.5mm} \ensuremath {\hfill \Box} }

\DeclareMathOperator{\im}{im}
\DeclareMathOperator{\Rank}{Rank}
\DeclareMathOperator{\rank}{rank}
\DeclareMathOperator{\gr}{gr}

\DeclareMathOperator{\coker}{coker}
\DeclareMathOperator{\Map}{Map}
\DeclareMathOperator{\gldim}{gldim}

\begin{document}

\begin{preliminary}

\maketitle

\begin{abstract}
Let $R$ be a subring of $\Q$ containing $\frac{1}{6}$ or $R=\Fp$ with
$p>3$.
Let $({A},{d})$ be a differential graded algebra (dga) such that its
homology, $H({A},{d})$, is $R$-free and is 
the universal enveloping algebra of some Lie algebra $L_0$.  
Given a free $R$-module $V$ and a map $d:V \to {A}$ such that
there is an induced map $d\p:V \to L_0$, let $\bB$ denote the
canonical dga extension $({A} \amalg \freeT V,d)$ where $\freeT V$ is
the tensor algebra generated by $V$.
$H\bB$ is studied in the case that the Lie ideal $[d\p V] \subset L_0$ is
a free Lie algebra.
This is a broad generalization of the previously studied \emph{inert}
condition.
An intermediate \emph{semi-inert} condition is introduced under which
the algebraic structure of $H\bB$ is determined.

This problem is suggested by the following \emph{cell attachment
problem}, first studied by J.H.C. Whitehead around 1940.
For a simply-connected finite-type CW-complex $X$, denote by $Y$ the space
obtained by attaching a finite-type wedge of cells to $X$.
Then how is the loop space homology $H_*(\lY;R)$ related to $H_*(\lX;R)$?
In some cases this topological problem is described by the previous
algebraic situation. 
Under a further hypothesis it is shown that if $H_*(\lX;R)$ is generated by
Hurewicz images then so is $H_*(\lY;R)$ and if $R \subset \Q$ then the
localization of $\lY$ at $R$ is homotopy equivalent to a product of
spheres and loop spaces on spheres. 

It is well-known that if 
the coefficient ring is a field,
Lie subalgebras of free Lie algebras are also free. 
To help implement the above results this fact is generalized, giving a
simple condition guaranteeing that Lie subalgebras of more general
Lie algebras are free. 
\end{abstract}


\begin{acknowledgements}
First, I would like to thank my advisor Paul Selick who has been an
invaluable mentor. 
His patience, encouragement and guidance made this thesis
possible.

I would like to thank Steve Halperin for introducing me to homotopy
theory and for many helpful discussions.
Many thanks to Jean-Michel Lemaire who kindly agreed to be my external
examiner and whose insight illuminated several aspects of my thesis.
Thanks to Lisa Jeffrey for her careful proofreading.
I would also like to thank John Holbrook who first showed me the joys
of mathematics.

Many thanks to the helpful staff at the Department of Mathematics and
especially to Ida Bulat whose cheerful presence I greatly appreciated.

I am indebted to NSERC, the government of Ontario, the University of
Toronto, Paul Selick, Steve Halperin and the AMS whose financial
support made my studies possible.

I am grateful to my parents who always encouraged my curiosity and are
a constant source of support.
Many thanks to my sisters Helen and Susan, whose friendship I value, and
thanks to all those who were friends while I was at the
University of Toronto. 

Finally I would like to thank my wife Jodi, whose love and support 
made my studies such a memorable and enjoyable experience.
\end{acknowledgements}

\tableofcontents

\end{preliminary}

\part {Introduction and Mathematical Review}

\chapter {Introduction and Summary of Results} \label{chapter-intro}

In the chapter we introduce and summarize our results and give an
outline of the thesis.

\section{Introduction}

In algebraic topology one seeks to understand topological spaces by
studying associated algebraic objects. 
For example each space comes with a set of homology groups and
homotopy groups. 
These groups are invariants of the collection of spaces homotopy
equivalent to the given space.

It is often the case that the homology groups are easy to calculate.
Unfortunately they do not carry as much information as one might like
about the associated space.
The homotopy groups on the other hand typically carry much more
information but are usually impossible to calculate.

An intermediate approach is to study the loop space homology groups of
a given space. 
Under the Pontrjagin product these invariants have an algebra
structure, giving us much more information about the space 
under consideration~\cite{felix:H*LoopX}.
For example, knowing the loop space homology is generated by Hurewicz
images may allow one to write the loop space as a product of spheres
and loops on spheres (see Section~\ref{section-hsb}). 
If this happens it gives the homotopy groups of the space in terms of
the homotopy groups of spheres, which have been extensively studied. 

This thesis presents a method for calculating the loop space homology
in some instances of the \emph{cell attachment problem} which we now
describe. 

Since topological spaces are weak homotopy equivalent to spaces which
are built out of cells (see Section~\ref{section-basic-top}), there is an
obvious need to understand how attaching one or more cells affects the
algebraic invariants associated to the space. 

The effect on homology is straight-forward.
In the torsion-free case attaching a cell will either add or remove
one homology generator. 
However the effect on homotopy and loop space homology is much more
mysterious.

\noindent
\textbf{The Cell Attachment Problem:} 
\emph{Given a simply-connected topological space $X$ what is the effect on
the loop space homology and the homotopy type if one attaches one or
more cells to $X$? } 

In some cases this problem is equivalent to the following purely
algebraic problem.

\noindent
\textbf{The Differential Graded Algebra Extension Problem:}
\emph{Given a differential graded algebra $(A,d)$ what is the effect on
the homology if one adds one or more generators to $(A,d)$? }

This thesis analyzes these problems under some assumptions.

The cell attachment problem was perhaps first considered by
J.H.C. Whitehead \cite{whitehead:addingRelations},
\cite[Section~6]{whitehead:simplicialSpaces}, around 1940 and more 
recently these two problems have been studied by Anick, F\'{e}lix,
Halperin, Hess, Lemaire and Thomas among others.
 
There are parallel theories analyzing the two problems.
For the topological problem, let $X$ be a finite-type simply-connected
topological space and let $Y$ be some space obtained by attaching
cells to $X$. 
For the algebraic problem let $(A,d)$ be a finite-type connected
differential graded algebra (dga) and let $\bB$ be some dga
obtained by adding generators to $(A,d)$.
In the first case let $I$ denote the  two-sided ideal of $H_*(\lX;R)$
generated by the image of the attaching maps (see
Section~\ref{section-cap}). 
In the second case let $I$ denote the two-sided ideal of $H(A,d)$
generated by the boundaries of the added generators.
The conditions under which $H_*(\lY;R) \isom H_*(\lX;R)/I$ and $H\bB
\isom H(A,d)/I$ have been extensively studied \cite{lemaire:autopsie,
anick:stronglyFree, halperinLemaire:inert, felixThomas:attach,
halperinLemaire:attach}.   
This is called the \emph{inert} condition.
The slightly weaker \emph{lazy} condition and the related \emph{nice}
condition have also been studied~\cite{felixLemaire:pontryagin,
halperinLemaire:attach, hessLemaire:nice}. 
(We will not use or define lazy, but we will define nice in
Remark~\ref{rem-nice} and in Section~\ref{section-cap}.) 
For more details see Sections \ref{section-cap} and \ref{section-inert-extns}.
Most of the results in the above papers were proved for the case where the
coefficient ring is a field. 
We will prove our results for both fields and subrings of the rational
numbers.

In this thesis we study a broad generalization of the inert
condition which we call the \emph{free} condition, and which we now define.
We assume that $H(A,d)$ is $R$-free and as algebras $H(A,d) \isom
UL_0$ for some Lie algebra $L_0$.
As explained in the next section this assumption is often satisfied
when the algebra extension is derived from a topological cell attachment.
Let $J \subset L_0$ be the Lie ideal generated by the boundaries of the
added generators.
The dga extension $\bB$ satisfies the free condition if $J$ is a free
Lie algebra.
If $\bB$ is a dga extension then it has an obvious filtration, which
we define in the next section.
If $\bB$ is a dga extension which satisfies the free condition
then we calculate $H\bB$ as an $R$-module and give the multiplicative
structure of the \emph{graded object} associated to the filtration (see
Section~\ref{section-basic-alg}). 

For dga extensions which satisfy the free condition we define a
\emph{semi-inert} condition which is a weaker generalization of the
inert condition.
Under this condition we are able to determine $H\bB$ including its
multiplicative structure. 

In case the topological problem for a space $Y$ is equivalent
to the algebraic problem for $\bB$ we obtain analogous results for
$H_*(\lY;R)$.
Conditions under which this happens are described in the next
section.
Previously studied special cases include the case where $R = \Q$ and
the case where $X$ is a wedge of spheres.

Under further topological conditions we show that if
$\lX$ is homotopy-equivalent to a product of spheres and loop spaces
of spheres then so is $\lY$.

It follows from our results that in certain cases the Lie algebra of
Hurewicz images automatically contains a free Lie algebra on two
generators.

Most of our results depend on showing that certain Lie subalgebras $J$
of a Lie algebra $L$ are free Lie algebras.
When $L$ is a free Lie algebra it is a well know fact that $J$ is a
free Lie algebra (when the coefficient ring is a field). 
However in our cases $L$ is generally not free.
Nevertheless we were able to prove a simple condition under which Lie
subalgebras of the given Lie algebras are free.

Our approach to the cell attachment problem generalizes work that has
been done in the special case that $X$ is a wedge of spheres,
in which case $Y$ is called a \emph{(spherical) two cone}.
This situation has been closely studied \cite{lemaire:monograph,
anick:thesis, fht:cat2, bogvad:envelopingAlg, felixThomas:cat2,
anick:cat2, felixThomas:attach, fht:lshCat12, bft:gammaFunctor,  
felixLemaire:2level, felixThomas:lshSmallCat, popescu:2cones}. 
Our approach is particularly indebted to~\cite{anick:cat2}.


We apply our results to various topological examples, construct an
infinite family of finite CW-complexes using only semi-inert cell
attachments, and give an algebraic version of Ganea's fiber-cofiber
construction.

\section{Summary of results}

In this summary we briefly introduce notation and definitions that will
be used throughout this thesis. 
A more detailed exposition will be given in Chapters \ref{chapter-dga},
\ref{chapter-ah} and \ref{chapter-cap}.

Let $R$ be a 
principal ideal domain containing $\frac{1}{6}$.
If $R$ is a subring of $\Q$ let $P$ be the set of invertible primes in
$R$ and let $\nP = \{p \in \Z \ \mid \ \text{p is prime and }p \notin
P\} \cup \{0\}$. 
All of our $R$-modules $M$ are graded, connected and have finite type.
Even when we do not explicitly mention it, our notation and
definitions depend on the choice of $R$.

Let $(\check{A},\check{d})$ be a differential graded algebra (dga)
over $R$. 
Let $Z \check{A}$ denote the subalgebra of cycles.

Let $V_1$ be a free $R$-module together with a map $d:V_1 \to Z
\check{A}$.
Let $\freeT V_1$ denote the tensor algebra on $V_1$ and let
$\check{A} \amalg \freeT V_1$ denote the coproduct or free product of
the algebras $\check{A}$ and $\freeT V_1$.
$\check{d}$ and $d$ can be canonically extended to a differential on
the coproduct which we also denote $d$.
Let $\bA = (\check{A} \amalg \freeT V_1, d)$ be the resulting dga.

There is an increasing filtration $\{F_k \bA\}$ on $\bA$ defined by
taking $F_{-1} \bA = 0$, $F_0 \bA = \check{A}$ and for $k\geq 0$,
$F_{k+1} = \sum_{i=0}^k F_i \bA \cdot V_1 \cdot F_{k-i}\bA$.
This induces a filtration on $H\bA$, the homology of $\bA$.
Each of these filtrations will be denoted by $\{F_k (\_\!\_) \}$.
If $M$ has an increasing filtration $\{F_k M\}$ then there exists an
associated graded object
\[ \gr(M) = \bigoplus_k ( F_k M / F_{k-1} M ).
\]
Let $\gr_k(M) = F_k M / F_{k-1} M$.

Assume that $H(\check{A},\check{d})$ is $R$-free and that as algebras
$H(\check{A},\check{d}) \isom UL_0$ for some Lie algebra $L_0$.
By~\cite[Theorem 8.3]{halperin:uea}, if $(\check{A},\check{d}) \isom
U(\check{L},\check{d})$ and $H(\check{A},\check{d})$ is $R$-free then
$H(\check{A},\check{d}) \isom UL_0$ as algebras 
(indeed as Hopf algebras) for some free $R$-module $L_0$.
There is an induced map 
\[ d': V_1 \to Z\check{A} \onto H(\check{A},\check{d}) \isomto UL_0.
\]
Assume one can choose $L_0$ so that $d'(V_1) \subset L_0$.  
In this case we say that $\bA$ is a \emph{dga extension of
$((\check{A},\check{d}),L_0)$}, which we will sometimes abbreviate to
a \emph{dga extension of $\check{A}$}. 

Letting $d' L_0 = 0$ one can canonically extend $d'$ to a
differential on $L_0 \amalg \freeL V_1$.
Let $\eL = (L_0 \amalg \freeL V_1, d')$.
$\eL$ is bigraded as a dgL, where the usual grading is called
dimension and the second grading, called degree, is given by letting
$L_0$ and $V_1$ be in degrees $0$ and $1$, respectively.
The sign conventions in a bigraded dgL use dimension.
$d'$ is a differential of bidegree $(-1,-1)$.
There is an induced bigrading on $U\eL$, $H\eL$ and $HU\eL$.
We record the following observation.

\begin{lemma} \label{lemma-dA} 
If $\check{d}=0$ and $\check{A} \isom U\check{L}$ then one can choose
$L_0 = \check{L}$. 
Then $\bA \isom U\eL$ and hence $\bA$ and $H\bA$ are bigraded.
As a result $H\bA \isom \gr (H\bA)$ as algebras in this case.
\end{lemma}

We will use the following notation throughout the paper:
If $M \subset X$ where $X$ is a Lie algebra or algebra then $[M]$ and
$(M)$ are respectively, the Lie ideal and two-sided ideal in $X$
generated by $M$. 
If $M \subset L$ it is a standard lemma that $U(L/[M]) \isom UL/(M)$
(see \cite[p.153; Theorem V.1(4)]{jacobson:lieAlgebras} for example).
If $M$ is bigraded then $M_i$ will denote the component of $M$ in
degree $i$.
Let $\Fp$ denote the finite field with $p$ elements and we use the
convention that $\F_0 = \Q$. 
For the $R$-module $M$, given any $p \in \nP$ denote $M \tensor \Fp$
by $\bar{M}$ and denote $d' \tensor \Fp$ by $\bar{d}$, with $p$
omitted from the notation.

Given Lie algebras $A$ and $B$ over $R$ such that $B$ is an
$A$-module, one can define the \emph{semi-direct product} $L = A \sdp
B$ as follows. 
As an $R$-module $L \isom A \times B$. 
The bracket between $A$ and $B$ is given by the action of $A$ on $B$. 
Equivalently (see Lemma~\ref{lemma-sdp}) there is a short exact
sequence of Lie algebras 
\[ 0 \to B \to L \xrightarrow{g} A \to 0
\]
with a Lie algebra section for $g$.

If $R = \F$ where $\F = \Q$ or $\Fp$ with $p>3$, define a dga extension of
$((\check{A},\check{d}),L_0)$, $\bA = (\check{A} \amalg \freeT V_1,
d)$, to be \emph{free} if $[d' V_1] \subset L_0$ is a free Lie algebra. 
If $R$ is a subring of $\Q$ define a dga extension of
$((\check{A},\check{d}),L_0)$, $\bA = (\check{A} \amalg \freeT V_1,
d)$, to be \emph{free} if \feachpnP, $[\bar{d} 
\bar{V}_1] \subset \bar{L}_0$ is a free Lie algebra.
In either case we introduce the following new condition.
Call $\bA$ \emph{semi-inert} if furthermore $\gr_1(H\bA)$ is a free
$\gr_0(H\bA)$-bimodule.
We will show that this is equivalent to both the condition that
$(H\eL)_1$ is a free $(H\eL)_0$-module and the condition that  
\[ (H\eL)_0 \sdp \freeL(H\eL)_1 \isom (H\eL)_0 \amalg \freeL K
\]
as Lie algebras for some free $R$-module $K \subset (H\eL)_1$.

This is a generalization of the \emph{inert} condition (defined in
Section~\ref{section-inert-extns}) of Anick~\cite{anick:stronglyFree} and
Halperin and Lemaire~\cite{halperinLemaire:inert} (see
Lemma~\ref{lemma-inertIsSemiInert}). 
It is a strict generalization.
For example we will show (see Example~\ref{eg-semi-inert-extn}) that when
$L = (\freeL \langle x,y,a,b \rangle, d)$ where $|x|=|y|=2, \ dx=dy=0,
\ da = [[x,y],x]$ and $db=[[x,y],y]$,
then $UL$ is not an inert extension of $\freeT \langle x,y \rangle$
but is a semi-inert extension of $\freeT \langle x,y \rangle$.

Our first main theorem is the following.

\begin{thma} \label{thm-a}
Let $R = \F$ where $\F= \Q$ or $\Fp$ with $p>3$.
Let $(\check{A}, \check{d})$ be a connected finite-type dga and let
$V_1$ be a connected finite-type $\F$-module with a map $d:V_1 \to
\check{A}$. 
Let $\bA = (\check{A} \amalg \freeT V_1, d)$.
Assume that there exists a Lie algebra $L_0$ such that
$H(\check{A},\check{d}) \isom UL_0$ as algebras and $d' V_1 \subset
L_0$ where $d'$ is the induced map.
Also assume that $[d' V_1] \subset L_0$ is a free Lie algebra.
That is, $\bA$ is a \emph{dga extension} of
$((\check{A},\check{d}),L_0)$ which is \emph{free}. 
Let $\eL = (L_0 \amalg \freeL V_1, d')$. \\ 
(i) Then as algebras
\[ \gr(H\bA) \isom U ((H\eL)_0 \sdp \freeL (H\eL)_1)
\]
with $(H\eL)_0 \isom L_0/[d' V_1]$ as Lie algebras.  \\
(ii) Furthermore if $\bA$ is \emph{semi-inert} then as algebras
\[ H\bA \isom U ( (H\eL)_0 \amalg  \freeL K' )
\]
for some $K' \subset F_1 H\bA$.
\end{thma}

\begin{rem}
Baues, F\'{e}lix and Thomas~\cite{bft:gammaFunctor,
felixThomas:lshSmallCat} showed that if $\bB$ is a \emph{two-level dga}
then $H\bB \isom \freeT_{H_0\bB} H_1\bB$ as algebras. 
Using Theorem~\ref{thm-a}(i) one can show that
$\gr(H\bA) \isom \freeT_{\gr_0(H\bA)} \gr_1{H\bA}$ as algebras.
\end{rem}

If $R$ is a subring of $\Q$, then we prove Theorem~\ref{thm-b}, a
slightly stronger analogue of Theorem~\ref{thm-a}. 

Let $\underline{\psi}: H\eL \to HU\eL$ be the map induced by
the inclusion $\eL \incl U\eL$.
This is a map of Lie algebras where $HU\eL$ is a Lie algebra under the
commutator bracket. 
Theorem~\ref{thm-b} is stronger than Theorem~\ref{thm-a} in that in
addition to the conclusions of Theorem~\ref{thm-a}, we will 
identify $\underline{\psi} H\eL$.

\begin{thma} \label{thm-b}
Let $R \subset \Q$ be a subring containing $\frac{1}{6}$.
Let $(\check{A}, \check{d})$ be a connected finite-type dga and let
$V_1$ be a connected finite-type free $R$-module with a map $d:V_1 \to
\check{A}$. 
Let $\bA = (\check{A} \amalg \freeT V_1, d)$.
Assume that there exists a Lie algebra $L_0$ such that
$H(\check{A},\check{d}) \isom UL_0$ as algebras and $d' V_1 \subset
L_0$ where $d'$ is the induced map.
Also assume that $L_0/[d' V_1]$ is a free $R$-module and that
\fanypnP, $[d' V_1] \subset L_0$ is a free Lie algebra. 
That is, $\bA$ is a \emph{dga extension} of
$((\check{A},\check{d}),L_0)$ which is \emph{free}. 
Let $\eL = (L_0 \amalg \freeL V_1, d')$. \\ 
(i) Then $H\bA$ and $\gr(H\bA)$ are $R$-free and as algebras
\[ \gr(H\bA) \isom U ( (H\eL)_0 \sdp \freeL (H\eL)_1 )
\]
with $(H\eL)_0 \isom L_0/[d' V_1]$ as Lie algebras.
Additionally $(H\eL)_0 \sdp \freeL (H\eL)_1 \isom \underline{\psi}
H\eL$ as Lie algebras.\\
(ii) Furthermore if $\bA$ is \emph{semi-inert} then as algebras 
\[ H\bA \isom U ( (H\eL)_0 \amalg  \freeL K' )
\]
for some $K' \subset F_1 H\bA$.
\end{thma}

Let $M(z)$ denote the \emph{Hilbert series} (see
Section~\ref{section-basic-alg}) for the $R$-module $M$.
If $M(0) = 1$ then let $M(z)^{-1}$ denote the series $1/M(z)$.
Using Anick's formula~\cite[Theorem 3.7]{anick:thesis}, the following
corollary follows from both Theorems \ref{thm-a} and \ref{thm-b}.

\begin{cor} \label{cor-semiInertEqnAlg}
Let $\bA$ be a semi-inert dga extension satisfying the hypotheses of
either Theorems \ref{thm-a} or \ref{thm-b}.
Let $K'$ be the $R$-module in $(ii)$ of the Theorem.
Then
\begin{equation} \label{eqn-semiInertEqnAlg}
K' (z) = V_1(z) + z[UL_0(z)^{-1} - U(H\eL)_0(z)^{-1}].
\end{equation}
\end{cor}

The above equation provides a necessary condition for semi-inertness.
If a dga extension is semi-inert then the right hand side of
\eqref{eqn-semiInertEqnAlg} should give a Hilbert series all of whose
terms have non-negative coefficients.
While this condition is not sufficient, 
for a `generic' dga extension which is not semi-inert there is no
reason to expect that the terms in the right hand side of
\eqref{eqn-semiInertEqnAlg} should have non-negative coefficients.

We now introduce a topological situation for which we prove Theorems
\ref{thm-c} and \ref{thm-d} analogous to Theorems \ref{thm-a} and
\ref{thm-b} and, with an additional hypothesis, the stronger
Theorems \ref{thm-hurewiczF} and \ref{thm-hurewiczR}.

Let $X$ and $W$ be simply-connected topological spaces with the
homotopy type of a finite-type CW-complex.
Using the Samelson product $\pi_*(\lX) \tensor R$ is a Lie algebra,
and under the commutator bracket $H_*(\lX;R)$ is a Lie algebra.
Using these products the Hurewicz map $h_X: \pi_*(\lX) \tensor R \to
H_*(\lX;R)$ is a Lie algebra map~\cite{samelson:products}.
Let $L_X$ denote its image.
A map $f: W \to X$ induces a map $L_W \to L_X$. 
Denote the image of this map by $L^W_X$, and let $[L^W_X]$ be the Lie
ideal in $L_X$ generated by $L^W_X$.
Note that the map $f$ is omitted from the notation.

Consider a continuous map $W \xrightarrow{f} X$, where $W = \bigvee_{j
\in J} S^{n_j}$ is a finite-type wedge of spheres and $f = \bigvee_{j \in
J} \alpha_j$.
The \emph{attaching map construction} (see
Section~\ref{section-basic-top}) gives the \emph{adjunction space} 
\[ Y = X \cup_f \left( \bigvee_{j \in J} e^{n_j+1} \right).
\]

Assume that $H_*(\lX;R)$ is torsion-free and that $H_*(\lX;R) \isom
UL_X$ as algebras. 
If $R$ is a field then there is no torsion and if $R=\Q$ then the
latter condition is trivial by the Milnor-Moore Theorem
(Theorem~\ref{thm-milnorMoore})~\cite{milnorMoore:hopfAlgebras}. 
If $R\subset \Q$ then in many examples one can reduce to the
torsion-free case by \emph{localizing}  (see
Section~\ref{section-lsd}) away from a finite set of primes. 
McGibbon and Wilkerson~\cite{mcgibbonWilkerson} have shown that any
finite CW-complex $X$ such that $\pi_*(X) \tensor \Q$ is finite
satisfies both of these conditions after localizing away from finitely
many primes (see Section~\ref{section-lsd}). 
Even if the loop space homology has torsion for infinitely many primes
\cite{anick:torsion, avramov:torsion}, one might be able to study the
space by including it in a larger torsion-free space~\cite{anick:cat2}.
If $H_*(\lX;R)$ is torsion-free then it is has been conjectured that
$H_*(\lX;R') \isom UL_X$~\cite{anick:slsd} where $R'$ is obtained
from $R$ by inverting finitely many primes.
This thesis will provide methods for verifying this conjecture
in certain cases (see Theorem~\ref{thm-hurewiczR} and
Corollary~\ref{cor-finite-complex}).

Let $\hat{\alpha}_j: S^{n_j-1} \to \lX$ denote the adjoint of
$\alpha_j$.
Let $\eL = (L_X \amalg \freeL \langle y_j \rangle_{j \in J}, d')$,
where $d' L_X = 0$ and $d' y_j = h_X(\hat{\alpha}_j)$.
We will justify reusing the notation $\eL$ in
Section~\ref{section-correspondence}. 
$\eL$ is a bigraded differential graded Lie algebra.

If $R = \F$ where $\F = \Q$ or $\Fp$ with $p>3$, define the attaching map $f$
to be \emph{free} if the Lie ideal $[L^W_X] \subset L_X$ is a free Lie
algebra.  
If $R$ is a subring of $\Q$ define the attaching map $f$ to be
\emph{free} if \feachpnP, $[\bar{L}^W_X]$ is a free Lie algebra.
In either case we introduce the following new condition.
Call $f$ \emph{semi-inert} if in addition
$(H\eL)_1$ is a free $(H\eL)_0$-module.
(We will see in Section~\ref{section-correspondence} why this
corresponds to same terminology defined earlier.)
We will show that this is equivalent to the condition that
\[ (H\eL)_0 \sdp \freeL (H\eL)_1 \isom (H\eL)_0 \amalg \freeL K
\]
as Lie algebras for some free $R$-module $K \subset (H\eL)_1$.
We will also show that a certain Adams-Hilton model (see
Section~\ref{section-ah}) provides a filtration for $H_*(\lY;R)$ under
which a free attaching map is semi-inert if and only if $\gr_1(H_*(\lY;R))$
is a free $\gr_0(H_*(\lY;R))$-bimodule.

Again this is a generalization of the inert condition \cite{anick:stronglyFree,
halperinLemaire:inert} (see Section~\ref{section-cap}).
For example we will show that the attaching map of the top cells in
the $6$-skeleton of $S^3 \times S^3 \times S^3$ 
($= (S^3_a \vee S^3_b \vee S^3_c) \cup_f \left(\bigvee_{i=1}^3 e^6\right)$, 
where $f = [\iota_b, \iota_c] \vee [\iota_c,\iota_a] \vee
[\iota_a,\iota_b]$) is not inert but is semi-inert (see
Example~\ref{eg-fat-wedge}).

In the statements of Theorems \ref{thm-c} and \ref{thm-d} below
$\gr(H_*(\lY;R)$ refers to the graded object associated to the
filtration mentioned above.

The following are two of our four main topological theorems. 
Recall that $L^X_Y$ is the image of the induced map $L_X \to L_Y$
between Hurewicz images.
Corresponding to Theorem~\ref{thm-a} we have

\begin{thma} \label{thm-c}
Let $R = \F$ where $\F = \Q$ or $\Fp$ with $p > 3$.
Let $X$ be a finite-type simply-connected CW-complex such that
$H_*(\lX;R)$ is torsion-free and as algebras $H_*(\lX;R) \isom UL_X$
where $L_X$ is the Lie algebra of Hurewicz images. 
Let $W = \bigvee_{j \in J} S^{n_j}$ be a finite-type wedge of spheres and
let $f:W \to X$.
Let $Y = X \cup_f \left(\bigvee_{j \in J} e^{n_j+1}\right)$.
Assume that $[L^W_X] \subset L_X$ is a free Lie algebra.
That is, $f$ is \emph{free}. \\
(i) Then as algebras 
\[ \gr(H_*(\lY; \F)) \isom U( {L}^X_Y \sdp \freeL (H\eL)_1)
\]
with ${L}^X_Y \isom {L}_X/[{L}^W_X]$ as Lie algebras. \\ 
(ii) Furthermore if $f$ is \emph{semi-inert} then as algebras 
\[ H_*(\lY;\F) \isom U ( L^X_Y \amalg \freeL K' )
\]
for some $K' \subset F_1 H_*(\lY; \F)$. 
\end{thma}

The following theorem corresponds to Theorem~\ref{thm-b}.

\begin{thma} \label{thm-d}
Let $R \subset \Q$ be a subring containing $\frac{1}{6}$.
Let $X$ be a finite-type simply-connected CW-complex such that
$H_*(\lX;R)$ is torsion-free and as algebras $H_*(\lX;R) \isom UL_X$
where $L_X$ is the Lie algebra of Hurewicz images. 
Let $W = \bigvee_{j \in J} S^{n_j}$ be a finite-type wedge of spheres and
let $f:W \to X$.
Let $Y = X \cup_f \left(\bigvee_{j \in J} e^{n_j+1}\right)$.
Assume that $L_X/[L^W_X]$ is $R$-free and that \feachpnP,  
$[\bar{L}^W_X] \subset \bar{L}_X$ is a free Lie algebra. 
That is, $f$ is \emph{free}. \\
(i) Then $H_*(\lY;R)$ and $\gr(H_*(\lY;R))$ are torsion-free and as algebras 
\[ \gr(H_*(\lY; R)) \isom U ( L^X_Y \sdp \freeL (H\eL)_1 )
\] 
with $L^X_Y \isom L_X/[L^W_X]$ as Lie algebras. \\ 
(ii) If in addition $f$ is \emph{semi-inert} then 
as algebras 
\[ H_*(\lY; R) \isom U ( L^X_Y \amalg \freeL K' )
\] 
for some $K' \subset F_1 H_*(\lY; R)$. 
\end{thma}

\begin{rem} \label{rem-nice}
Since $UL_X/(L^W_X) \isom U(L_X/[L^W_X])$ it follows from (i) in both
Theorems \ref{thm-c} and \ref{thm-d} that
$UL_X/(L^W_X)$ injects in $H_*(\lY;R)$. 
This is the definition of Hess and Lemaire's~\cite{hessLemaire:nice}
\emph{nice} condition for $f$.
\end{rem}

\begin{rem}
In Theorems \ref{thm-a} and \ref{thm-b} (respectively Theorems
\ref{thm-c} and \ref{thm-d}) if one assumes that
$H(\check{A},\check{d}) \isom UL_0$ (respectively $H_*(\lX;R) \isom
UL_X$) as Hopf algebras then 
it is straightforward to check that the
isomorphisms in parts (i) the Theorems are indeed Hopf algebra
isomorphisms.   
We will not need to use this fact.
\end{rem}

The following corollary follows from Theorems \ref{thm-c} and
\ref{thm-d}.

\begin{cor}
If $f$ is a non-inert free cell attachment satisfying the conditions
of either Theorem~\ref{thm-c} or Theorem~\ref{thm-d} and $\dim
(H\eL)_1 \neq 1$ then $H_*(\lY;R)$ contains a tensor algebra on two
generators. 
\end{cor}

\begin{rem}
F\'{e}lix and Lemaire~\cite[Corollary 1.4]{felixLemaire:pontryagin}
showed%
\footnote{The trivial case, which is illustrated by
Example~\ref{eg-semi-inert-extn}, is omitted
in~\cite{felixLemaire:pontryagin}.}
that for a non-inert free attachment either there exists an $N$ such
that $\dim H_i(\lY;R) \leq N$ for all $i$, or $H_*(\lY;R)$ contains a
tensor algebra on two generators. 
\end{rem}

Let $\tilde{H}_*(W)(z)$ be the Hilbert series (see
Section~\ref{section-basic-alg}) for the reduced homology of $W$. 
The following corollary follows from both Theorems \ref{thm-c} and
\ref{thm-d}.

\begin{cor} \label{cor-semiInertEqnTop}
Let $f$ be a semi-inert attaching map satisfying the hypotheses of
either Theorems \ref{thm-c} or \ref{thm-d}.
Let $K'$ be the $R$-module in $(ii)$ of the Theorem.
Then
\begin{equation} \label{eqn-semiInertEqnTop}
K' (z) = \tilde{H}_*(W)(z) + z[UL_X(z)^{-1} - U(L^X_Y)(z)^{-1}].
\end{equation}
\end{cor}

Like Corollary~\ref{cor-semiInertEqnAlg}, this corollary provides a
necessary condition for an attaching map to be semi-inert.
Note that $W$ and $L_X$ are already known and $L^X_Y
\isom L_X / [L^W_X]$, which in some cases is easy to calculate.

We now proceed to our next two main topological theorems.

The following is an important collection of spaces.
Let $\mathcal{S} = \{ S^{2m-1},\,\Omega S^{2m+1} \:|\: m \geq 1\}$.
A space $Y$ is called \emph{atomic} if for some $r$ it is $r$-connected,
$\pi_{r+1}(Y)$ is a cyclic abelian group, and any self-map $f:Y\rightarrow Y$
inducing an isomorphism on $\pi_{r+1}(Y)$ is a homotopy equivalence.
The spaces in $\mathcal{S}$ are atomic.
Let $\PiS$ be the collection of spaces homotopy equivalent
to a \emph{weak product} (see Section~\ref{section-basic-top} for the
definition) of spaces in $\mathcal{S}$. 

Let $W$, $X$, and $Y$ be the spaces in Theorems \ref{thm-c} and
\ref{thm-d}.
If in addition to the hypotheses of Theorems \ref{thm-c} and \ref{thm-d}
there exists a Lie algebra map $\sigma_X: L_X \to \pi_*(\lX) \tensor R$ right
inverse to $h_X$ then we also have the stronger topological results below.
In the case where $X$ is a wedge of spheres, studied by
Anick~\cite{anick:cat2}, such a map always exists (see
Example~\ref{eg-2cones}). 
In Theorem~\ref{thm-hurewiczR} we will give a sufficient conditions for
extending such a map from a subspace, but for other cases its
existence is an open problem.  
When $R \subset \Q$, in Section~\ref{section-ip} we will use such a
map to define the set of \emph{implicit primes} labeled $P_Y$ which
contains the invertible primes in $R$.
Intuitively the implicit primes are those primes $p$ for which
$p$-torsion is used in the cell-attachments.
If $W$ is finite and the set of invertible primes in $R$ is finite
then $P_Y$ is finite as well (see Lemma~\ref{lemma-ip}).
The next two theorems are analogues of the fundamental Milnor-Moore Theorem
(Theorem~\ref{thm-milnorMoore}) which holds for $R = \Q$.

\begin{thma} \label{thm-hurewiczF}
Let $R = \F$ where $\F = \Q$ or $\Fp$ with $p >3$.
Let $\bigvee_{j \in J} S^{n_j} \xto{\bigvee \alpha_j} X$ be a cell
attachment satisfying the hypotheses of Theorem~\ref{thm-c}.
Let $Y = X \cup_{\bigvee \alpha_j} \left(\bigvee e^{n_j+1}\right)$
and let $\hat{\alpha}_j$ denote the adjoint of $\alpha_j$. 
In addition assume that there exists a map $\sigma_X$ right inverse to
the Hurewicz map $h_X$ and that $\forall j \in J$, $\sigma_X h_X
\hat{\alpha}_j = \hat{\alpha}_j$. 
Then the canonical algebra map
\[ UL_Y \to H_*(\lY;\F)
\] 
is a surjection.
That is, $H_*(\lY;\F)$ is generated as an algebra by Hurewicz images.
\end{thma}

Let $R \subset \Q$ be a subring with invertible primes $P$ and let $X$
be a simply-connected topological space.
Then there is a topological analogue to the localization of
$\Z$-modules at $R$ (see Section~\ref{section-hsb}).
We will call this the \emph{localization at $R$} or the
\emph{localization away from $P$}.
We denote the localization of $X$ at $R$ by $X_{(R)}$.

\begin{rem} \label{rem-hurewiczF}
If Theorem~\ref{thm-hurewiczF} holds for $\Q$ and for $\Fp$ for
all non-invertible primes in some $R \subset \Q$ containing $\frac{1}{6}$
and $H_*(\lY;R)$ is torsion-free then by the Hilton-Serre-Baues
Theorem (Theorem~\ref{thm-hsb}) $H_*(\lY;R) \isom UL_Y$ as algebras
and localized at $R$, $\lY \in \PiS$. 
\end{rem}

\begin{thma} \label{thm-hurewiczR}
Let $R \subset \Q$ be a subring containing $\frac{1}{6}$.
Let $\bigvee_{j \in J} S^{n_j} \xto{\bigvee \alpha_j} X$ be a cell
attachment satisfying the hypotheses of Theorem~\ref{thm-d}.
Let $Y = X \cup_{\bigvee \alpha_j} \left(\vee e^{n_j+1}\right)$.
Furthermore assume that there exists a map $\sigma_X$ right inverse to
the Hurewicz map $h_X$.
Let $P_Y$ be the set of \emph{implicit primes} and let $R' =
\Z[{P_Y}^{-1}]$. 
Then \\
(i)  $H_*(\lY;R')$ is torsion-free and as algebras 
\[ H_*(\lY;R') \isom UL_Y
\] 
where $\gr(L_Y) \isom L^X_Y \sdp \freeL (H\eL)_1$ as Lie algebras, and \\ 
(ii) localized away from $P_Y$, $\lY \in \PiS$. \\
(iii) If in addition $f$ is semi-inert then localized away from $P_Y$,
$L_Y \isom L^X_Y \amalg \freeL \hat{K}$ as Lie algebras for some
$\hat{K} \subset F_1 L_Y$, and there exists a map $\sigma_Y$ right
inverse to $h_Y$.
\end{thma}

\begin{rem}
By the Hilton-Serre-Baues Theorem (Theorem~\ref{thm-hsb}), the
assumption that $H_*(\lX;R)$ is torsion-free and that $H_*(\lX;R) \isom
UL_X$ as algebras, is equivalent to the assumption that localized at
$R$, $\lX \in \PiS$.
Therefore under the conditions of Theorem~\ref{thm-hurewiczR}, if $\lX_{(R)}
\in \PiS$ and there are no new invertible primes then $\lY_{(R)} \in
\PiS$.
\end{rem}

This theorem is important for the following reason. If we know that
$\lY \in \PiS$ then we can apply results about the homotopy groups of
spheres, which have been extensively studied, to make some strong
conclusions about $Y$.
For example, localized away from $P_Y$, $Y$ satisfies the \emph{Moore
conjecture} (Conjecture~\ref{conj-wmc})~\cite{selick:mooreConj} which
relates the torsion of the homotopy groups of $Y$ to the torsion-free
behaviour of the $\pi_*(Y)$.
Furthermore, it creates the possibility that one might be able to
demonstrate the Moore conjecture for other spaces by following Anick's
plan~\cite{anick:cat2} of embedding them in spaces for which $\lY \in
\PiS$.



The following two corollaries follow from Theorem \ref{thm-hurewiczR}.

\begin{cor} \label{cor-finite-complex}
If $Y$ is a finite complex constructed out of semi-inert attachments
then there is a finite set of primes $P_Y$ such that localized away
from $P_Y$, $\lY$ is homotopy equivalent to a product of spheres and
loops on spheres. 
Equivalently $H_*(\lY;\Z[{P_Y}^{-1}]) \isom UL_Y$.
\end{cor}

\begin{cor}
If $f$ is a non-inert free cell attachment satisfying the conditions
of Theorem~\ref{thm-hurewiczR} and $\dim (H\eL)_1 \neq 1$ then $L_Y$
contains a free Lie algebra on two generators. 
\end{cor}

To apply our results one needs to know that certain Lie ideals are
free Lie algebras.
This is difficult to check in general.
However, it is a well-known result that over a field $\F$, any Lie subalgebra
of a free (graded) Lie algebra is a free Lie algebra
\cite{shirshov:subalgebras, mikhalev:subalgebras, mikhalevZolotykh}.
This is referred to as the \emph{Schreier property}~\cite{mikhalevZolotykh}.

Unfortunately our Lie ideals are in Lie algebras that are generally not
free so we cannot use this property.
By the results of this thesis, our ideals are often in Lie algebras
whose associated graded object is a certain semi-direct product. 
For this more general situation we have succeeding in giving a simple
condition from which one can conclude that a given Lie
ideal is in fact a free Lie algebra.

\begin{thma}  \label{thm-schreier2}
Over a field $\F$, let $L$ be a finite-type graded Lie algebra with filtration
$\{F_k L\}$ such that $\gr(L) \isom L_0 \sdp \freeL V_1$ as Lie
algebras, where $L_0 = F_0 L$ and $V_1  = F_1 L / F_0 L$.
Let $J \subset L$ be a Lie subalgebra such that $J \cap F_0 L = 0$.
Then $J$ is a free Lie algebra.
\end{thma}

As an application of our results we obtain an
algebraic analogue of Ganea's fiber-cofiber construction, which we
discuss in Chapter~\ref{chapter-ganea}.
Like Ganea's construction, our construction can be iterated.

We apply our results to various topological examples in
Chapter~\ref{chapter-eg}. 
For example we illustrate our results by calculating the loop space
homology of a particular \emph{spherical three-cone} and show that its
loop space is homotopy equivalent to a product of spheres and loop spaces on
spheres. 
By induction we produce an infinite family of finite CW-complexes
constructed out of semi-inert cell attachments. 
We also give a more abstract analysis of (spherical) two-cones,
three-cones and $N$-cones, as well as some other examples.

\ignore
{

By Lemma~\ref{lemma-ip}, $P_Y = \{2,3\}$.

Now 
\begin{equation} \label{eqn-LXinSIeg} 
L_X \isom L_{Z_1} \amalg L_{Z_2} \isom L^{A_1}_{Z_1} \amalg
L^{A_2}_{Z_2} \amalg \freeL \langle w_1, w_2 \rangle.
\end{equation}
Therefore by Theorem~\ref{thm-schreier2} $f$ is a free attaching map.
Thus $Y$ satisfies the hypotheses of Theorem~\ref{thm-d}.

By Corollary~\ref{cor-semiInertEqnTop}, if $f$ is semi-inert then 
\[
K' (z) = \tilde{H}_*(W)(z) + z[UL_X(z)^{-1} - U(L^X_Y)(z)^{-1}].
\]
$L_X$ is given by~\eqref{eqn-LXinSIeg} and $L^X_Y \isom L_X /
[h_X(\hat{\beta_1}), h_X(\hat{\beta_2})]$. 
Since the Hurewicz images are $[[w_1,w_2],w_1]$ and $[[w_1,w_2],w_1]$
which are contained in $\freeL \langle w_1, w_2 \rangle$ we have
\[ 
L^X_Y \isom L^{A_1}_{Z_1} \amalg L^{A_2}_{Z_2} \amalg \tilde{L}
\text{ where } \tilde{L} = \freeL \langle w_1, w_2 \rangle / [
R \{ [[w_1,w_2],w_1],[[w_1,w_2],w_2] \} ].
\]

Since $(A \amalg B)(z)^{-1} = A(z)^{-1} + B(z)^{-1} -
1$~\cite[Lemma 5.1.10]{lemaire:monograph} it follows that
\[
K'(z) = \tilde{H}_*(W)(z) + z[U\freeL \langle w_1,w_2 \rangle -
U\tilde{L}(z)^{-1}].
\]
Now $U\freeL \langle w_1,w_2 \rangle \isom \freeT \langle w_1,w_2
\rangle$ and as $R$-modules $\tilde{L} \isom R\{
w_1,w_2,[w_1,w_1],[w_1,w_2],[w_2,w_2]\}$ and $U\tilde{L} \isom \mathbb{S}
\tilde{L}.$
Therefore
\begin{eqnarray*}
K'(z) & = & 2z^{28} + z \left[ 1-2z^9 -
\frac{(1-z^{18})^3}{(1+z^9)^2} \right] \\ 
& = & z^{37}.
\end{eqnarray*}

Recall that $\eL = (L_X \amalg \freeL \langle e,g \rangle, d')$ where
$d' e = [[w_1,w_2],w_1]$ and $d' g = [[w_1,w_2],w_2]$.
Also recall (from Theorem~\ref{thm-c}) that $(H\eL)_0 \isom L^X_Y$.
Let $\bar{u} = [e,w_2] + [g,w_1]$ (with $|\bar{u}| = 37$).
Then it is easy to check that
\[ (H\eL)_0 \sdp \freeL (H\eL)_1 \isom (H\eL)_0 \amalg \freeL \langle
\bar{u} \rangle
\]
and $f$ is indeed semi-inert.
Since $\sigma_{Z_1}$ and $\sigma_{Z_2}$ are right inverses to
$h_{Z_1}$ and $h_{Z_2}$ it follows that $\sigma_X = \sigma_{Z_1} \vee
\sigma_{Z_2}$ is right inverse to $h_X$.
Let $P_Y$ be the finite set of implicit primes for $Y$.
Depending on the choice of $\omega$ this may be contain a finite set of primes
in addition to $2$ and $3$.
Let $R' = \Z[{P_Y}^{-1}]$.
As a result by Theorem~\ref{thm-hurewiczR}
\[ H_*(\lY;R') \isom UL_Y \text{ where } L_Y \isom L^X_Y \amalg \freeL
\langle u \rangle
\]
with $u = h_Y(\hat{\nu})$ for some map $\nu: S^{38} \to Y$.
Therefore localized away from $P_Y$, $\lY \in \PiS$.
Furthermore localized away from $P_Y$, there exists a map $\sigma_Y$
right inverse to $h_Y$. 
\eolBox
\end{eg}

}

\ignore
{

\section{Future directions}

The results of this thesis lead to three areas for future research as
well as a number of specific open questions.

\subsection{Torsion in loop space homology}

One of the motivations for this thesis was to give loop space
decompositions for spaces with torsion (possibly at infinitely many
primes) in their loop space homology.
The methods of this thesis can applied to this study.
For example I can identify the loop space for Anick's space with
torsion of all orders~\cite{anick:torsion} as a product of simpler
factors.

\subsection{The Moore conjecture}

Another motivation for this thesis was to generalize Anick's
verification of the weak Moore conjecture for finite
two-cones~\cite{anick:cat2} to arbitrary CW-complexes.
It is hoped that this line of research we lead to a generalization of
Anick's results.

\subsection{Rational Homotopy Theory}

Localized over the rationals, I can write any CW-complex with finite
cone length as a finite sequence for free cell attachments. 
Using this framework one should be able to derive many of the main
results of rational homotopy theory.
In addition one may be able to obtain new results using this approach.

}

\section{Outline of the thesis}

In Chapter \ref{chapter-dga} we review some algebraic constructions
and introduce the differential graded algebra (dga) extension problem.

In Chapter \ref{chapter-ah} we review some topology, introduce Adams-Hilton
models~\cite{adamsHilton} and prove some of their properties which will be
needed in Chapter~\ref{chapter-hurewicz}. 

In Chapter \ref{chapter-cap} we introduce the cell attachment problem and
compare it to the previously discussed dga extension problem. 
We also review some topological theory.

In Chapter \ref{chapter-HUL} we prove Theorems \ref{thm-a} and
\ref{thm-b}, the main algebraic theorems. 

In Chapter \ref{chapter-attach} we prove Theorems \ref{thm-c} and
\ref{thm-d}, the topological analogues of Theorems \ref{thm-a} and
\ref{thm-b}. 

In Chapter \ref{chapter-hurewicz} we prove Theorems
\ref{thm-hurewiczF} and \ref{thm-hurewiczR} which identify Hurewicz images.

In Chapter \ref{chapter-freeLa} we prove Theorem~\ref{thm-schreier2}
which will be useful in applying our results to examples.

In Chapter \ref{chapter-ganea} we show that our results can be applied
to give an algebraic construction analogous to Ganea's Fiber-Cofiber
construction~\cite{ganea:construction}. 

In Chapter \ref{chapter-eg} we apply our results to some topological
examples and conclude with some open questions.

\chapter{Differential Graded Algebra} \label{chapter-dga}

In this chapter we review some standard algebraic constructions and
results. 
We also motivate and define our main algebraic objects of study:
\emph{dga extensions} which are \emph{free} and \emph{semi-inert}.
Furthermore we introduce the spectral sequence with which we study
these objects.

Let $R$ be a principal ideal domain containing $\frac{1}{2}$ and
 $\frac{1}{3}$.
This is equivalent to saying that $R$ is a principal ideal domain
 containing $\frac{1}{6}$.
This assumption is made so that the \emph{Hurewicz images} (defined in
 Section~\ref{section-basic-top}) have the
 structure of a \emph{(graded) Lie algebra} (defined in
 Section~\ref{section-basic-alg}). 

\section{Filtered, differential and graded objects} \label{section-basic-alg}

A \emph{graded $R$-module} is an $R$-module $M$ such that $M =
\bigoplus_{i \in \Z} M_i$ where $M_i$ is an $R$-module.
All of our $R$-modules will be assumed to be graded.
$R$ itself is a graded $R$-module concentrated in degree $0$.
If $a \in M_i$ let $|a| = i$.
$M$ is said to have \emph{finite type} if for all $i$, $M_i$ is a finite
dimensional $R$-module.
$M$ is said to be connected if $M_i = 0$ for $i \leq 0$.
Since $R$ is a principal ideal domain, $M$ is a \emph{free} $R$-module
if it is torsion-free. 
That is, for $m \in M$, $r \in R$, $rm=0 \implies r=0$.
Equivalently we will also say that $M$ is \emph{$R$-free}.
We will always assume that our $R$-modules are free and have finite type.
Given a set $S$, let $RS$ be the free $R$-module with basis $S$.
A \emph{graded $R$-module map} preserves the gradation and the
$R$-module structure.
That is, $f(r\cdot m) = r \cdot f(m)$ and $|f(a)|=|a|$.
All of our maps will be assumed to be graded $R$-module maps.

Though all of our definitions depend on $R$, we will usually leave
this dependence implicit.

For a torsion-free graded $R$-module $M$, the \emph{Hilbert series}
for $M$ is the power series $M(z) = \Sigma_{i \in \Z} (\Rank_R M_i) z^i$.

A \emph{(graded) Lie algebra} is a graded $R$-module $L$ together with a
linear map $[ \ , \ ]: L_i \tensor L_j \to L_{i+j}$ (called the
bracket) satisfying the
following relations: 
\begin{enumerate}
\item $[x,y] = -(-1)^{|x||y|}[y,x]$ (anti-commutativity) and 
\item $[x,[y,z]] = [[x,y],z] + (-1)^{|x||y|} [y,[x,z]]$ (the Jacobi
identity). 
\end{enumerate}
Note that since $\frac{1}{6} \in R$, $[x,x]=0$ when $|x|$ is even and
$[[x,x],x] = 0$ for all $x$.

A \emph{Lie algebra map} is a map between Lie algebras which preserves
the Lie algebra structure.
That is, $f([a,b]) = [f(a), f(b)]$.
A Lie algebra is said to be \emph{connected} if it is connected as an
$R$-module. 
A \emph{Lie subalgebra} $J \subset L$ is an $R$-submodule which is
closed under the bracket.
That is, for $a,b \in J$, $[a,b] \in J$.
$J$ is a \emph{Lie ideal} of $L$ if furthermore for all $a \in J$ and
for all $b \in L$, $[a,b] \in J$.
If $J$ is a Lie ideal then there is an induced Lie algebra structure
on the quotient $L/J$.

We will use the following notation for iterated brackets 
\[
[[a_1, a_2, \ldots a_{n-1}, a_n] = [[\cdots [a_1, a_2], \ldots
a_{n-1}], a_n].
\]
If $J$ is the Lie ideal generated by $a \in L$ we will sometime denote
$L/J$ by $L/a$.

A \emph{(graded) algebra} is a graded $R$-module $A$ together with $R$-module
maps $\mu: A_i \tensor A_j \to A_{i+j}$ (multiplication) and $\nu: R \to A_0$
(unit) such that: 
\begin{enumerate}
\item the composite $A \isom A \tensor R \xto{A \tensor \nu} A
\tensor A \xto{\mu} A$ equals $1_{A}$;
\item the composite $A \isom R \tensor A \xto{\nu \tensor A} A
\tensor A \xto{\mu} A$ equals $1_{A}$; and 
\item $\mu \circ (\mu \tensor 1_A) = \mu \circ (1_A \tensor \mu): A
\tensor A \tensor A \to A$. 
\end{enumerate}
We remark that $\mu(a \tensor b)$ is usually denoted by $a \cdot b$ or
just $ab$. 
An algebra $A$ is said to be \emph{connected} if $A_i=0$ for $i<0$ and
$A_0 \isom R$.
An \emph{algebra map} is a map between algebras which preserves
the algebra structure.
That is, $f(a \cdot b) = f(a) \cdot f(b)$.

We will assume that all of our $R$-modules, Lie algebras and algebras
are connected, $R$-free and have finite type.

A \emph{subalgebra} $B \subset A$ is an $R$-submodule which is closed
under multiplication (that is, for $a,b \in B$, $\mu(a \tensor b) \in
B$) and which contains the unit (that is, $\nu(R) \subset B$).
A \emph{two-sided ideal} of $A$ is an $R$-submodule $I \subset A$ such
that for all $a \in A$ and for all $b \in I$, $\mu(a \tensor b) \in I$
and $\mu(b \tensor a) \in I$.

Let $T$ denote the (twist) map $V \tensor W \to W \tensor V$ given by $v
\tensor w \mapsto (-1)^{|v||w|} w \tensor v$.

A \emph{(graded, cocommutative) Hopf algebra} is an algebra $H$
together with $R$-module maps $\Delta: H \to H \tensor H$
(comultiplication) and $\epsilon: H \to R$ (counit) such that 
\begin{enumerate}
\item the composite $H \xto{\Delta} H \tensor H \xto{H \tensor \epsilon}
H \tensor R \isom H$ equals $1_H$;
\item the composite $H \xto{\Delta} H \tensor H \xto{\epsilon \tensor H}
R \tensor H \isom H$ equals $1_H$; 
\item $(\Delta \tensor 1_H) \circ \Delta = (1_H \tensor \Delta) \circ
\Delta: H \to H \tensor H \tensor H$;
\item $T \circ \Delta = \Delta: H \to H \tensor H$; and 
\item $\Delta$ is a homomorphism of algebras, where the multiplication
on $H \tensor H$ is given by $H \tensor H \tensor H \tensor H \xto{1
\tensor T \tensor 1} H \tensor H \tensor H \tensor H \xto{\mu \tensor
\mu} H \tensor H$.
\end{enumerate}

A \emph{differential (graded) $R$-module (dgm)} is a $R$-module $M$
together with a linear map $d: M_i \to M_{i-1}$ such that $d^2 = 0$,
called the \emph{differential}.
Equivalently we say $M$ is an $R$-module with a differential.
A \emph{differential graded Lie algebra (dgL)} is a Lie algebra $L$ together
with a differential $d$ such that $d[x,y] = [dx,y] +
(-1)^{|x|}[x,dy]$ (the Leibniz rule).
A \emph{differential graded algebra (dga)} is an algebra $A$ together
with a differential $d$ such that $d(xy) = (dx)y + (-1)^{|x|}x(dy)$
(also the Leibniz rule).


A \emph{filtered dgm} is a dgm $M$ together with a filtration $\{F_i M\}_{i
\in \Z}$ such that $ \ldots \subset F_i M \subset F_{i+1} M \subset
\ldots \subset M$ and $d (F_i M) \subset F_i M$.
A \emph{filtered dgL} is a dgL which has a dgm filtration such that
$[F_i L, F_j L] \subset F_{i+j} L$.
A \emph{filtered dga} is a dga which has a dgm filtration such that
$F_i A \cdot F_j A \subset F_{i+j} A$.
A map of filtered objects is a map $f: M \to N$ such that $f(F_i M)
\subset F_i N$. That is, $f$ preserves the filtration.
Given a filtered dgm there exists an \emph{associated graded object }
\[ \gr(M) = \bigoplus_{i \in \Z} F_i M / F_{i-1} M,
\]
which becomes a dgm with the differential induced from the
differential on $M$.
If $L$ is a filtered dgL then the Lie bracket on $L$ induces a Lie
algebra structure on $\gr(L)$.
Similarly if $A$ is a filtered dga then the product on $A$ induces an
algebra structure on $\gr(A)$.
If $f: M \to N$ is a filtered map then there exists an induced map
$\gr(f): \gr(M) \to \gr(N)$.
We will sometimes refer to $\gr(M)$ by $\gr_*(M)$ where $\gr_i(M) =
F_i M / F_{i-1} M$. 

An $R$-module $M$ 
is \emph{bigraded} if it has a second grading. 
That is, $M = \bigoplus_{i,j \in \Z} M_{i,j}$ where $M_{i,j}$ is an
$R$-module.
A (graded) algebra is \emph{bigraded} if it is bigraded
as an $R$-module and $\mu: A_{i,j} \tensor A_{k,l} \to A_{i+j,k+l}$.
In this thesis we will say that a (graded) Lie algebra is
\emph{bigraded} if it is bigraded as an $R$-module and $[ \cdot,
\cdot]: L_{i,j} \tensor L_{k,l} \to L_{i+j,k+l}$.
Note that the anti-commutativity and Jacobi relations need not hold with
respect to the second grading. 
A \emph{differential bigraded $R$-module} is a bigraded $R$-module
with a map $d: M_{i,j} \to M_{i-1,j-1}$ such that $d^2=0$.
We will say that a dga/dgL is \emph{bigraded} if it is a bigraded
algebra/Lie algebra and a differential bigraded $R$-module.
Note that the Leibniz rule need not hold for the second grading.

We will usually denote the first grading by dimension and the
second grading by degree.
The dimension and degree of $a$ will be denoted by $|a|$ and $\deg(a)$
respectively. 
The anti-commutativity and Jacobi relations and the Leibniz rule all
use $|\cdot|$ and not $\deg(\cdot)$.

A \emph{short exact sequence} of $R$-modules is a sequence of
$R$-modules $A,B,C$ together with $R$-module maps 
\[
0 \to A \xto{f} B \xto{g} C \to 0
\] 
such that $\ker(f) = 0$, $\ker(g) = \im(f)$, and $\im(g) = 0$.
It is a short exact sequence of Lie algebras if each of the
modules is a Lie algebra and each of the maps is a Lie algebra map.
The short exact sequence of Lie algebras 
$0 \to A \xto{f} B \xto{g} C \to 0$
is said to split if there exists a Lie algebra map 
$h: C \to B$ such that $g \circ h = id_{C}$.

\section{Free algebras and Lie algebras and products} \label{section-sdp}

Let $V$ be a free $R$-module.
Let $\freeT V$ be the free algebra generated by $V$.
As an $R$-module $\freeT V \isom \bigoplus_{k=0}^{\infty} \freeT^k V$
where $\freeT^0 V = R$ and for $k\geq 1$, $\freeT^k V = \underbrace{V
\tensor \cdots \tensor V}_{k times}$.
The unit is the inclusion $R = \freeT^0 V \incl \freeT V$.
Multiplication $\freeT^i V \tensor \freeT^j V \to \freeT^{i+j} V$ is
given by $a \cdot b = a \tensor b$.
Any algebra of this form is called a \emph{tensor algebra}.

Given a free $R$-module $V$ we define the \emph{free commutative
algebra} or \emph{symmetric algebra} by $\mathbb{S} V = \freeT V / I$,
where $I$ is the two-sided ideal generated by the commutator brackets
$[a,b] = a \cdot b - (-1)^{|a||b|} b \cdot a$.

Recall that the commutator bracket gives any algebra the structure of
a Lie algebra.
In particular $\freeT V$ is a Lie algebra with the commutator bracket.\
Let $\freeL V$ be the Lie subalgebra of $\freeT V$ generated by
$V$.
We call any Lie algebra of this form a \emph{free Lie algebra}.

The coproduct between two algebras and the coproduct between two Lie
algebras is defined categorically~\cite{hiltonStammbach:book}.
One can check that the following constructions satisfy the
definitions. 
Given two tensor algebras $\freeT V$ and $\freeT W$ we define the
\emph{coproduct} or \emph{free product} $\freeT V \amalg \freeT W =
\freeT (V \oplus W)$. 
Similarly we define $\freeL V \amalg \freeL W = \freeL (V \oplus W)$.
More generally given two algebras $A$ and $B$ we can define the
coproduct $A \amalg B$ to be the free $R$-module on words of the
form $a_1 b_1 a_2 b_2 \cdots a_n b_n$, $a_1 b_1 a_2 b_2 \cdots b_{n-1}
a_n$, $b_1 a_1 \cdots a_n b_n$ and $b_1 a_1 \cdots b_{n-1} a_n$ (where
$a_i \in A$ and $b_i \in B$) modulo the relations $w a_k 1 a_{k+1} w\p
= w (a_k a_{k+1}) w\p$ and  $w b_k 1 b_{k+1} w\p = w (b_k b_{k+1})
w\p$ where $w, w\p$ are words in the coproduct.
Multiplication is given by concatenation.
Similarly we define $L_1 \amalg L_2$ for two Lie algebras $L_1$ and
$L_2$.

Given a Lie algebra $L_0$, an $R$-module $M$ is an \emph{$L_0$-module} if
there exists a map $[\cdot,\cdot]: L_0 \tensor M \to M$ such that
$[x_1, [x_2, y]] = [[x_1,x_2],y] + (-1)^{|x_1||x_2|} [x_2,[x_1,y]]$.

\begin{rem} \label{rem-L0-module}
Given a Lie algebra $L_0$ and a $L_0$-module $M$ then $\freeL M$ is
also a $L_0$-module where the action of $L_0$ on $\freeL M$ is given
inductively by the Jacobi identity and the action of $L_0$ on $M$.
If $x \in L_0$ and $[[y_1,\ldots,y_{n-1},y_n] \in M$ then
\[
[x,([[y_1,\ldots,y_{n-1},y_n])] = [[x,([[y_1,\ldots,y_{n-1]}],y_n] 
+ (-1)^{|x||y|}[([[y_1,\ldots,y_{n-1}]),[x,y_n]].
\]
\end{rem}

Given a Lie algebra $L_0$ and a $L_0$-module $N$ we define the
\emph{semi-direct product} $L = L_0 \sdp N$ to be the following
Lie algebra.
As an $R$-module $L \isom L_0 \times N$.
For $a,b \in L_0$ (respectively $N$), $[a,b]$ is given by the
Lie bracket in $L_0$ (respectively $N$).
For $x \in L_0$ and $y \in N$, $[x,y]$ is given by the action
of $L_0$ on $N$.

The \emph{direct product} $L_0 \times N$ is the special case of
the semi-direct product where $[x,a]=0$ if $x \in L_0$ and $a \in
N$.

\begin{lemma} \label{lemma-sdp}
Let $L$ be a Lie algebra. Then $L$ is a semi-direct product $L_0 \sdp
N$ if and only if there exists a split short exact sequence of
Lie algebras 
\[ 0 \to N \xto{f} L \xto{g} L_0 \to 0.
\]
\end{lemma}

\begin{proof}
($\Rightarrow$) Since $L_0$ acts on $N$. That is, $[L_0, N]
\subset N$, $N$ is a Lie ideal.
Thus there exists a quotient map $\rho: L \to L_0$.
Hence the inclusion and projection give a short exact sequence of Lie
algebras 
\begin{equation} \label{eqn-sdp-rho} 
0 \to N \to L \xto{\rho} L_0 \to 0.
\end{equation}
Let $h: L_0 \to L$ be the inclusion map.
Then $\rho \circ h = id_{L_0}$. So~\eqref{eqn-sdp-rho} is a split
short exact sequence of Lie algebras. \\
($\Leftarrow$) Since $0 \to \freeL \xto{f} L \xto{g} L_0 \to 0$ is a
short exact sequence of $R$-modules, $L \isom L_0 \times N$ as
$R$-modules. 
Let $h: L_0 \to L$ be the splitting. 
Since $g \circ h = id_{L_0}$, $h$ is an injection.
Since $f$ and $h$ are injections of Lie algebras, for $a,b \in L_0$
(respectively $N$), $[a,b]$ within $L$ is given by the bracket in $L_0$
(respectively $N$).
Since $g$ is a Lie algebra map, $\ker(g)$ is a Lie ideal.
Since $f$ is an injection $N = \im(f) = \ker(g)$ and $N$
is a Lie ideal. 
Therefore $[L_0, N] \subset N$. That is, $N$ is a
$L_0$-module.
So $L \isom L_0 \sdp N$.
\end{proof}

\begin{lemma} \label{lemma-sdp-extension}
Assume $L_1$ is a $L_0$-module.
Given a Lie algebra map $u_0: L_0 \to L$ and an $L_0$-module map $u_1:
L_1 \to L$ (eg. for $x \in L_0$, $y \in L_1$, $u_1([y,x]) =
[u_1(y),u_0(x)]$) then there exists an induced Lie algebra map $u: L_0
\sdp \freeL L_1 \to L$.
\end{lemma}

\begin{proof}
Since $L_1$ is a $L_0$-module, by Remark~\ref{rem-L0-module} so is
$\freeL L_1$.
Since there exists a map $u_1: L_1 \to L$ there is a canonical induced Lie
algebra map $\tilde{u}: \freeL L_1 \to L$ (given by $\tilde{u}|_{L_1}
= u_1$ and otherwise defined inductively by $\tilde{u}([a,b]) =
[\tilde{u}(a), \tilde{u}(b)]$).
It is easy to prove by induction that this is an $L_0$-module map.
That is, $\tilde{u}([y,x]) = [\tilde{u}(y),u_0(x)]$.
Since $L_0 \sdp \freeL L_1 \isom L_0 \times \freeL L_1$ as $R$-modules
this defines a map $u: L_0 \sdp \freeL L_1 \to L$.
Finally for $x \in L_0$ and $y \in \freeL L_1$, $[x,y] \in \freeL L_1$
so $u([x,y]) = \tilde{u}([x,y]) = [u_0(x), \tilde{u}(x)] = [u(x),u(y)]$.
Thus $u$ is a Lie algebra map.
\end{proof}

\section{Universal enveloping algebras}

Let $L$ be a Lie algebra.
Define $UL$ the \emph{universal enveloping algebra} on $L$ as follows. 
$UL = \freeT L / I$ where $I$ is the two-sided ideal generated by 
$\{x \tensor y - (-1)^{|x||y|} y \tensor x - [x,y]\}$. 
If $L$ is a dgL then $\freeT L$ is a dga under the induced differential:
$d(ab) = (da)b + (-1)^{|a|}a(db)$.
This dga structure induces a dga structure on $UL$.
Similarly, if $L$ is filtered there is an induced filtration on
$\freeT L$ which induces a filtration on $UL$.

It is a basic result~\cite[p.168; Theorem V.7]{jacobson:lieAlgebras} that
$U\freeL V \isom \freeT V$. 

The \emph{commutator bracket} $[x,y] = xy - (-1)^{|x||y|} yx$ gives
any algebra the structure of a Lie algebra. 

The universal enveloping algebra is universal in the following sense.
Given any algebra (dga) $A$ and a Lie algebra (dgL) map $f$
\[
\xymatrix{
L \ar[r]^f \ar@{^{(}->}[d] & A \\
UL \ar@{-->}[ur]_g
}
\]
there exists a unique map $g$ making the diagram commute.
In fact $UL$ can be defined as the unique object satisfying this
universal property.

The Poincar\'{e}-Birkhoff-Witt Theorem is the central theorem in the
study of universal enveloping algebras (see for example
\cite{jacobson:lieAlgebras, selick:book}). 
We give a weak form of this Theorem.

\begin{thm}[Poincar\'{e}-Birkhoff-Witt Theorem] \label{thm-pbw}
Let $L$ be a Lie algebra. 
Let $\mathbb{S}L$ denote the free commutative algebra (symmetric
algebra) (see Section~\ref{section-sdp}) on the $R$-module $L$. 
Then there is an $R$-module isomorphism $\mathbb{S} L \isomto UL$
which restricts to the identity on $L$.
\end{thm}

\begin{rem}
Although we will only need the weak form given above, the full theorem
specifies that is isomorphism is an isomorphism of \emph{coalgebras}.
\end{rem}

\begin{lemma} \label{lemma-Ugr=grU}
If $L$ is a filtered Lie algebra such that $\gr(L) \isom L$ as
$R$-modules then the canonical map $U\gr(L) \to \gr(UL)$ is an isomorphism. 
\end{lemma} 

\begin{proof}
If $L$ is filtered then $UL$ has an induced filtration and the
inclusion $L \incl UL$ preserves the filtration.
This induces a map $\gr(L) \to \gr(UL)$ which induces a map $U\gr(L)
\to \gr(UL)$.
Similarly there is an induced map $\mathbb{S}\gr(L) \to \gr(\mathbb{S}L)$.\

Using the Poincar\'{e}-Birkhoff-Witt Theorem, these fit into the
following commutative diagram of $R$-linear morphisms.
\[
\xymatrix{
{\mathbb{S}} \gr(L) \ar[rr]^{\phi} \ar[dd]_{\isom} & &
\gr({\mathbb{S}}L) \ar[dd]^{\isom} \\ 
& \gr(L) \ar@{_{(}->}[ul] \ar@{^{(}->}[dl] \ar[ur] \ar[dr] \\
U \gr(L) \ar[rr] & & \gr(UL)
}
\]

Since $\mathbb{S}L$ depends only on the $R$-module structure on $L$, the
$R$-module isomorphism $L \isomto \gr(L)$ induces an algebra isomorphism
\begin{equation} \label{eqn-sl}
\mathbb{S}L \isomto \mathbb{S}\gr(L).
\end{equation}
From this we get the following canonical commutative diagram.
\[
\xymatrix{
\gr(L) \ar[r] \ar@{^(->}[d] & \gr({\mathbb{S}}L) \ar[r]^{\isom} &
\gr({\mathbb{S}}\gr(L)) \\
{\mathbb{S}}\gr(L) \ar[ur]^{\phi} \ar[urr]_{\isom}
}
\]


Thus $\phi$ is an isomorphism, and therefore the canonical map
$U\gr(L) \to \gr(UL)$ is an isomorphism.
\end{proof}

Let $\bL$ be a differential graded Lie algebra.
The canonical inclusion $\iota: \bL \to U\bL$ induces a map on homology
$H(\iota): H\bL \to HU\bL$ which induces the natural map
\[ \psi: UH\bL \to HU\bL.
\]
Since $H(\iota)$ is the restriction of $\psi$ to $H\bL$ we will usually
also refer to this map as $\psi$.

If $R = \Q$ then Quillen~\cite{quillen:rht} showed that $\psi: UH\bL \to
HU\bL$ is an isomorphism. 
We can also state this result as follows: 
\[
HU\bL \isom U \psi (H\bL) \text{, when } R = \Q.
\]
In Theorems \ref{thm-a} and \ref{thm-b} we show that under certain
hypotheses $HU\bL$ can be calculated from $H\bL$ for more general
coefficient rings.
See~\cite{popescu:UHLandHUL} for other results of this type.

\section{Inert, free and semi-inert dga extensions} \label{section-inert-extns}

Let $R = \Q$ or $\Fp$ where $p>3$ or $R \subset \Q$ is a subring
containing $\frac{1}{6}$.
If $R \subset \Q$ let $P$ be the set of invertible primes in $R$. 
Let $\nP = \{ p \in \Z | \text{ p is prime and } \ p \notin P\} \cup \{0\}$.
Furthermore for an $R$-module $V$ and a $R$-module map $f$, for each
$p \in \nP$ denote $V \tensor \Fp$ by $\bar{V}$ and $f \tensor \Fp$ by
$\bar{f}$. 
Note that the prime is omitted from the notation.

Let $(\check{A},\check{d})$ be a dga such that
$H(\check{A},\check{d})$ is $R$-free and that as algebras
$H(\check{A},\check{d}) \isom UL_0$ for some Lie algebra $L_0$.
If $H(\check{A},\check{d})$ is $R$-free and $(\check{A},\check{d})
\isom U(\check{L},\check{d})$ then by~\cite{halperin:uea},
$H(\check{A},\check{d}) \isom UL_0$ as algebras for some Lie algebra
$L_0$.
Let $Z\check{A}$ denote the cycles in $\check{A}$.
Let $\bA = (\check{A} \amalg \freeT V, d)$ where $d|_{\check{A}} =
\check{d}$ and $dV \subset Z \check{A}$.
Then there is an induced map 
\[ d\p: V \xto{d} Z \check{A} \to H(\check{A},\check{d}) \isomto
UL_0.
\]
If for some choice of $L_0$, $d\p:V \subset L_0$ then call $\bA$ a
\emph{dga extension of $((\check{A},\check{d}),L_0)$}. 
We will sometimes abbreviate this to a dga extension of $\check{A}$.

\begin{rem} 
The technical condition of the existence of a choice of $L_0$ such
that $d\p:V \subset L_0$ is satisfied in all the examples in this thesis.
We do not know if there are examples arising from cell attachments or
such that $(\check{A},\check{d}) \isom U(\check{L},\check{d})$ for
which there does not exist such a choice of $L_0$.
\end{rem}

If $R = \Q$ or $\Fp$ where $p > 3$, call $\bA$ an \emph{inert}
extension if $H\bA \isom U(L_0/[dV])$ 
\cite{halperinLemaire:inert, anick:stronglyFree} (Anick uses the
terminology \emph{strongly-free} and he gives a nice
combinatorial characterization in the case where $L_0$ is a free Lie
algebra).
Recall that all of our $R$-modules have finite type.
If $R \subset \Q$ we will say that $\bA$ is \emph{inert} if
$L_0/[dV]$ is torsion-free and $H\bA \isom U(L_0/[dV])$.

\begin{rem} \label{rem-inertR}
Note that if $R \subset \Q$ and $\bA$ is inert then $H\bA$ is
$R$-free. 
Thus by the Universal Coefficient Theorem,  for each $p \in \nP$,
$H(\bA \tensor \Fp) \isom U(\bar{L}_0/ [\bar{d} \bar{V}])$. That is
$\bA \tensor \Fp$ is inert over $\Fp$ ($p$-inert). 
\end{rem}

\begin{eg} \label{eg-inert-dga}
$(\check{A},\check{d}) = U (\freeL \langle x,y \rangle, 0)$.
$\bA = (\check{A} \amalg \freeT \langle a \rangle, d)$, where $da =
[[x,y],y]$
\end{eg}

Using the results of~\cite{anick:stronglyFree}, $\bA$ is an inert dga
extension (since in Anick's terminology $\{xyy\}$ is a
\emph{combinatorially free} set).  
It follows that $H\bA \isom U(\freeL \langle x,y \rangle /
[[x,y],y])$.
\eolBox

The following theorem gives an equivalent characterization of the
inert condition.
We will use it as motivation to generalize the inert condition.
Let $J$ be the Lie ideal $[dV] \subset L_0$.
The following is proved by Halperin and
Lemaire~\cite[Theorem 3.3]{halperinLemaire:inert} in the case when $R
= \Q$, and by Felix and Thomas~\cite[Theorem 1]{felixThomas:attach} in
the case when $R$ is a field of characteristic $\neq 2$. 

\begin{thm}[\cite{halperinLemaire:inert
}] \label{thm-inert}
If $R=\Q$ or $\Fp$ where $p>3$ then $\bA$ is an inert extension if and
only if \\ 
(i) $J$ is a free Lie algebra, and \\
(ii) $J/[J,J]$ is a free $U(L/J)$-module.
\end{thm}

We will study extensions which in general do not satisfy the second
condition. 
If $R = \Q$ or $\Fp$ where $p>3$ we define $\bA$ to be a \emph{free}
dga extension 
if $[d\p V] \subset L_0$ is a free Lie algebra.
If $R \subset \Q$ we define the dga extension $\bA$ to be \emph{free}
if for each $p \in \nP$, $[\bar{d} \bar{V}] \subset \bar{L}_0$ is a
free Lie algebra. 

This condition is a broad generalization of the inert condition and we
will prove results about $H\bA$ if $\bA$ is a free dga extension
(see part (i) of Theorems \ref{thm-a} and \ref{thm-b}).
However many free dga extensions satisfy a second condition.
We define this condition below and call it the \emph{semi-inert}
condition.
Under this condition we are able to determine $H\bA$ as an algebra
(see part (ii) of Theorems \ref{thm-a} and \ref{thm-b}).
We will show that it too is a generalization of the inert condition.

We will use the following filtration.
$\bA$ is filtered by taking $F_{-1}\bA = 0$, $F_0 \bA = \check{A}$ and for
$n\geq 0$, $F_{n+1} \bA = \sum_{i=0}^n F_i \bA \cdot V_1 \cdot F_{n-i} \bA$.
Since $dF_n \bA \subset F_n \bA$ and $F_n \bA \cdot F_m \bA \subset
F_{n+m} \bA$ this makes $\bA$ into a filtered dga.
There is an induced dga filtration on $H\bA$.
Let $\gr_*(H\bA)$ be the associated graded object.
For degree reasons $\gr_1(H\bA)$ is a $\gr_0(H\bA)$-bimodule.

Let $\eL = (L_0 \amalg \freeL V, d\p)$.
Then $\eL$ and $H\eL$ are bigraded.
So for degree reasons $(H\eL)_1$ is a $(H\eL)_0$-module.
Thus by Remark~\ref{rem-L0-module} we can form the semi-direct
product $(H\eL)_0 \sdp \freeL (H\eL)_1$. 

We will use the following lemma to define the semi-inert condition.
Since this condition is not used in the statements or the proofs of
Theorem~\ref{thm-a}(i) and Theorem~\ref{thm-b}(i) we will take the
liberty of using these results. 

\begin{lemma} \label{lemma-semi-inert-dga}
Let $\bA$ be a free dga extension. Then the following conditions are
equivalent: \\
(a) $(H\eL)_0 \sdp \freeL (H\eL)_1 \isom (H\eL)_0 \amalg \freeL K$ 
for some free $R$-module $K \subset (H\eL)_1$, \\
(b) $(H\eL)_1$ is a free $(H\eL)_0$-module, and \\
(c) $\gr_1(H\bA)$ is a free $\gr_0(H\bA)$-bimodule.
\end{lemma}

\begin{proof}
(b) $\implies$ (a)
Let $K$ be a basis for $(H\eL)_1$ as a free $(H\eL)_0$-module.
Then $(H\eL)_0 \sdp \freeL (H\eL)_1 \isom (H\eL)_0 \amalg \freeL
K$. \\
(a) $\implies$ (c)
Since $\bA$ is a free dga extension, by Theorem~\ref{thm-a}(i) or Theorem~\ref{thm-b}(i), 
$\gr_*(H\bA) \isom U \left( (H\eL)_0 \sdp \freeL (H\eL)_1 \right)$.
So by (a), 
\[
\gr_*(H\bA) \isom U \left( (H\eL)_0 \amalg \freeL K \right) \isom
\gr_0(H\bA) \amalg \freeT K,
\]
for some free $R$-module $K \subset (H\eL)_1$. 
Therefore 
\[
\gr_1(H\bA) \isom \left[ \gr_0(H\bA) \amalg \freeT K \right]_1 \isom
\gr_0(H\bA) \tensor K \tensor \gr_0(H\bA).
\]
(c) $\implies$ (b)
Let $L' = (H\eL)_0 \sdp \freeL (H\eL)_1$.
Then by Theorem~\ref{thm-a}(i) or Theorem~\ref{thm-b}(i), $\gr_*(H\bA)
\isom UL'$ and $\gr_1(H\bA) \isom (UL')_1$.
By (c), $(UL')_1$ is a free $(UL')_0$-bimodule.
Then we claim that it follows that $L'_1$ is a free $L'_0$-module.
Indeed, if there is a nontrivial degree one relation in $L'$ then
there is a corresponding nontrivial degree one relation in $UL'$.
\end{proof}

We say that a dga extension $\bA$ is \emph{semi-inert} if $\bA$ is
free and it satisfies the conditions of the previous lemma.
We justify this terminology with the following.

\begin{lemma} \label{lemma-inertIsSemiInert}
An inert dga extension is semi-inert.
\end{lemma}

\begin{proof}
Let $\bA$ be an inert dga extension.
If $R = \Q$ or $\Fp$ then by Theorem~\ref{thm-inert} $[dV]$ is a
free Lie algebra.
If $R \subset \Q$ then by Remark~\ref{rem-inertR} \feachpnP, $\bA
\tensor \Fp$ is $p$-inert. 
So by Theorem~\ref{thm-inert} \feachpnP, $[\bar{d}\bar{V}]$ is a free
Lie algebra.
So in either case $\bA$ is a free dga extension.

Since $(H\eL)_1 = 0$, the semi-inert condition is trivially satisfied.
\end{proof}

\begin{eg} \label{eg-semi-inert-extn}
A semi-inert dga extension which is not inert.
\end{eg}

Let $R = \Fp$ where $p>3$ or $R \subset \Q$ containing $\frac{1}{6}$.
Let $L = (\freeL \langle x,y,a,b \rangle, d)$ where
$|x|=|y|=2, \ |a|=|b|=7, \ dx=dy=0, \ da = [[x,y],x]$ and $db=[[x,y],y]$.
We will show that $UL$ is semi-inert dga extension of $\freeT\langle
x,y \rangle$ but not an inert dga extension of $\freeT \langle x,y \rangle$. 

Let $(\check{A}, \check{d}) = U(\freeL \langle x,y \rangle, 0)$.
Then $H(\check{A}, \check{d}) \isom U\freeL \langle x,y \rangle$, $d\p
= d$ and $UL$ is a dga extension of $U \freeL \langle x,y \rangle$.
Since $\freeL \langle x,y \rangle$ is a free Lie algebra, if $R$ is a
field then by the Schreier property the Lie ideal $[R\{da,db\}]
\subset \freeL \langle x,y \rangle$ is automatically a free Lie algebra.
If $R \subset \Q$ then \feachpnP, the Lie ideal
$[\Fp\{\bar{d}a,\bar{d}b\}] \subset \freeL_{\Fp} \langle x,y \rangle$ is
automatically a free Lie algebra.
Thus in either case $UL$ is a free dga extension.

Let $w = [a,y] - [b,x]$.
Then by anti-commutativity,
\begin{eqnarray*}
dw & = & [[[x,y],x],y] - [[[x,y],y],x] \\
& = & [[x,y],[x,y]] \\
& = & 0.
\end{eqnarray*}
Since $w$ is not a boundary $0 \neq [w] \in (HL)_1$ and $0 \neq [w]
\in (HUL)_1$.
Thus $UL$ is not an inert dga extension.
By the definition of homology 
\[
(HL)_0 \isom \freeL \langle x,y
\rangle \left/ \, \bigl[ R\{[[x,y],x],[[x,y],y]\} \bigr]. \right.
\]

One can check that $(HL)_1$ is freely generated by the $(HL)_0$-action
on $[w$].
Thus $(HL)_0 \sdp \freeL (HL)_1 \isom (HL)_0 \amalg \freeL \langle [w]
\rangle$.
That is, $UL$ is a semi-inert extension.
Therefore by part (ii) of either Theorem~\ref{thm-a} or
Theorem~\ref{thm-b},  
\[ 
HUL \isom U( (HL)_0 \amalg \freeL \langle [w] \rangle )
\]
as algebras.
%
\eolBox

\begin{eg} \label{eg-alg-fat-wedge}
Another semi-inert dga extension which is not inert.
\end{eg}

Let $R = \Fp$ where $p>3$ or $R \subset \Q$ containing $\frac{1}{6}$.
Let $L = (\freeL \langle x,y,z,a,b,c \rangle, d)$ where
$|x|=|y|=|z|=2, \ |a|=|b|=|c|=5 \ dx=dy=dz=0, \ da = [y,z], \
db=[z,x]$, and $dc=[x,y]$.   
$UL$ is a free dga extension of $U\freeL \langle x,y,z \rangle$.
As in the previous example, we will show that $UL$ is a semi-inert
extension but not an inert extension.


Let $w = [x,a] + [y,b] + [z,c]$.
By the Jacobi identity $dw = 0$.
Since $w$ is not a boundary $0 \neq [w] \in (HL)_1$ and $0 \neq [w]
\in (HUL)_1$.
Thus $UL$ is not an inert dga extension.
By the definition of homology $(HL)_0 \isom \freeL_{ab} \langle x,y,z
\rangle$, where $\freeL_{ab}$ denotes the free abelian Lie algebra
(that is all brackets are zero).

One can check that $(HL)_1$ is freely generated by the $(HL)_0$-action
on $[w$].
Therefore $UL$ is a semi-inert extension and by part (ii) of either Theorem
\ref{thm-a} or Theorem \ref{thm-b}, 
\[ HUL \isom U( (HL)_0 \amalg \freeL \langle [w] \rangle )
\]
as algebras.
\eolBox

\begin{eg}
A free extension which is not semi-inert.
\end{eg}

Let $R = \Fp$ where $p>3$ or $R \subset \Q$ containing $\frac{1}{6}$.
Let $L = (\freeL \langle x,a \rangle, d)$, where $|x|$ is odd, $dx=0$
and $da=[x,x]$.
Then $UL$ is a free dga extension of $U \freeL \langle x \rangle$.
Let $u = [x,a]$
Then $du = [x,[x,x]] = 0$ by the Jacobi identity.
$u$ is not a boundary so $0 \neq [u] \in (HL)_1$ and $0 \neq [u] \in
(HUL)_1$.
Thus $UL$ is not an inert extension.
By the definition of homology $(HL)_0 \isom \freeL_{ab} \langle [x]
\rangle$, where $\freeL_{ab}$ denotes the free abelian Lie algebra
(where all brackets are zero).
However $d(\frac{1}{4}[a,a]) = [u,x]$, so $\left[ [u],[x] \right] = 0$
and $HUL \neq U( \freeL_{ab} \langle [x] \rangle \amalg
\freeL \langle [u] \rangle)$.

Since the cycles of $L$ in degree one are just $R\{u,[u,x]\}$ and
 $[u,x]$ is a boundary, by the definition of homology $(HL)_1 =
 R\{[u]\}$.
So by part (i) of either Theorem \ref{thm-a} or Theorem \ref{thm-b}, 
$HUL \isom U( \freeL_{ab}  \langle [x] \rangle \sdp \freeL \langle [u]
 \rangle)$ as algebras, where the  semi-direct product is given by
 $[[x],[u]] = 0$. 
Thus 
\[
HUL \isom U(\freeL_{ab} \langle [x] \rangle \times \freeL
 \langle [u] \rangle) \isom U(\freeL_{ab} \langle [x],[u] \rangle)
\]
as algebras.
\eolBox

\ignore
{

\section{Inert, free and semi-inert dga extensions} \label{section-dga-extns}

When applying our algebraic results to topological examples (via
\emph{Adams-Hilton models}) we will actually need a slightly more
general situation than the dga extensions discussed above. 

Let $(\check{A},\check{d})$ be a differential graded algebra such that
$H(\check{A},\check{d})$ is torsion-free and $H(\check{A},\check{d})
\isom UL_0$ as algebras.
Let $\bA = (\check{A} \amalg \freeT V,d)$ where $d|_{\check{A}} =
\check{d}$ and $dV \subset \check{A}$.
Then there is an induced map 
\[
d': V \xto{d} Z\check{A} \to H(\check{A},\check{d}) \isomto UL_0.
\]
If for some choice of $L_0$, $d'V \subset L_0$ then call $\bA$ a
\emph{dga extension of $((\check{A},\check{d}),L_0)$}.

We use the same definitions as for dga extensions, to define inert,
free and semi-inert dga extensions.
In particular,
if $R = \Q$ or $\Fp$ where $p > 3$, call $\bA$ an \emph{inert}
extension if $H\bA \isom U(L_0/[d'V])$ 
If $R \subset \Q$ we will say that $\bA$ is \emph{inert} if
$L_0/[d'V]$ is torsion-free and $H\bA \isom U(L_0/[dV])$.

If $R = \Q$ or $\Fp$ we define $\bA$ to be a \emph{free} dga extension
if $[d\p V] \subset L_0$ is a free Lie algebra.
If $R \subset \Q$ we define the dga extension $\bA$ to be \emph{free}
if for each $p \in \nP$, $[\bar{d} \bar{V}] \subset \bar{L}_0$ is a
free Lie algebra. 

Let $\eL = (L_0 \amalg \freeL V, d\p)$.
Then $\eL$ and $H\eL$ are bigraded.
So for degree reasons $(H\eL)_1$ is a $(H\eL)_0$-module.
Thus by Remark~\ref{rem-L0-module} we can form the semi-direct
product $(H\eL)_0 \sdp \freeL (H\eL)_1$. 
We say that $\bA$ is \emph{semi-inert} if $\bA$ is free and 
\[ 
(H\eL)_0 \sdp \freeL (H\eL)_1 \isom (H\eL)_0 \amalg \freeL K
\] 
for some free $R$-module $K \subset (H\eL)_1$.

}

\section{The spectral sequence of a dga extension} \label{section-ss}

Let $\bA$ be a dga extension of $((\check{A},\check{d}),L_0)$.
Then by definition $\bA = (\check{A} \amalg \freeT V_1,d)$
where $d|_{\check{A}} = \check{d}$ and $dV_1 \subset Z \check{A}$. 
Also by definition $H(\check{A},\check{d}) \isom UL_0$ and $d\p V_1
\subset L_0$ where $d\p$ is the induced differential $d\p: V_1 \xto{d}
Z \check{A} \to H(\check{A},\check{d}) \isomto UL_0$. 

$\bA$ is filtered by taking $F_{-1}\bA = 0$, $F_0 \bA = \check{A}$ and for
$n\geq 0$, $F_{n+1} \bA = \sum_{i=0}^n F_i \bA \cdot V_1 \cdot F_{n-i} \bA$.
Since $dF_n \bA \subset F_n \bA$ and $F_n \bA \cdot F_m \bA \subset
F_{n+m} \bA$ this makes $\bA$ into a filtered dga.
There is an induced dga filtration on $H\bA$.

From the filtration of $\bA$ there is an associated first quadrant
spectral sequence 
\[ \gr(\bA) \implies \gr(H\bA)
\]
which thus converges \cite{mccleary:usersGuide, selick:book}.
Anick studied this spectral sequence~\cite{anick:thesis} and showed
that it collapses under certain conditions. 
In Chapter~\ref{chapter-HUL} we will show that Anick's conditions are
satisfied if the dga extension is free.

The $E^0$ term is given by $E^0_{p,q} = F_p \bA_{p+q} / F_{p-1}
\bA_{p+q}$, where $\bA_k$ denotes the component of $\bA$ in dimension
$k$.
The differential $d^0$ is the induced differential (from $d$) on
$\gr(\bA)$.
Since $\check{d}$ does not lower filtration but $d|_{V_1}$ does,
$d^0|_{\check{A}} = \check{d}$ and $d^0|_{V_1} = 0$.
In fact $(E^0 \bA, d^0) = (\check{A} \amalg \freeT V_1, \check{d})$, where
$\check{d}|_{\freeT V_1} = 0$.
Therefore 
\[
E^1 \bA = H(E^0 \bA, d^0) \isom H(\check{A},\check{d}) \amalg \freeT
V_1 \isom U(L_0 \amalg \freeL V_1).
\]
One can check that the induced differential $d^1$ is just the induced
differential $d\p$.
Therefore 
\[ E^2 \bA \isom HU\eL \text{ where } \eL = (L_0 \amalg \freeL
V_1,d\p).
\]

Since the spectral sequence converges, $E^{\infty} \isom \gr H\bA$.
Our main algebraic result (see Chapter~\ref{chapter-HUL}) will involve
showing that when the dga extension is free, the associated spectral
sequence collapses at the $E^2$-term. That is, $E^{\infty} = E^2$.

Unfortunately since $\bA$ is not bigraded as a dga, it is not
necessarily the case that $\gr(H\bA) \isom H\bA$ as algebras.
The following dga extension (which is not free) illustrates this.

\begin{eg}
$(\check{L}, \check{d}) = (\freeL \langle
x_1,x_2,y_1,y_2,u_1,u_2,v_1,v_2 \rangle, \check{d})$ where for
$i=1,2$, $|x_i| = |y_i| = 2$, $\check{d} x_i = \check{d} y_i = 0$, 
$\check{d} u_i = [[x_1,x_2],y_i]$ and $\check{d} v_i =
[[y_1,y_2],x_i]$.
Let $(\check{A},\check{d}) = U(\check{L},\check{d})$ and let
$\bA = (\check{A} \amalg \freeT \langle a,b \rangle, d)$ where $da = x_1$
and $db = [y_1,y_2]$. 
\end{eg}

Let $\{ F_i \bA \}$ be the usual filtration.
Let $L_0 = H(\check{L}, \check{d})$. 
Abusing notation we will refer to the homology classes represented by
$x_i$ and $y_i$ by $x_i$ and $y_i$.

For $i=1,2$ let $\alpha_i = [[a,x_2],y_i] - u_i$ and $\beta_i =
[b,x_i] - v_i$. 
One can check that these are cycles in $F_1 \bA$ which are not
boundaries in $\bA$.

Let $\gamma = -[\beta_1,x_2] - [\alpha_1,y_2] + [\alpha_2,y_1] \in
Z F_1 \bA$.
Let $\epsilon = [v_1,x_2] - [v_2,x_1] + [u_1,y_2] - [u_2,y_1] \in Z
F_0 \bA$.
Then one can check $\gamma$ is not a boundary in $\bA$.
Therefore $[\gamma] \neq 0 \in F_1 H\bA$.
Thus 
\begin{equation} \label{eqn-ungraded} 
-[[\beta_1],[x_2]] - [[\alpha_1],[y_2]] + [[\alpha_2],[y_1]] \neq 0 \in
H\bA.
\end{equation}

However
\[ d([[a,b],x_2] - [a,v_2]) = \gamma - \epsilon.
\]
As a result $[\gamma] = [\epsilon] \in F_0 H\bA$.
For a cycle $\bar{z}$ in $\bA$ let $[\bar{z}]$ denote the corresponding
homology class in $\gr(H\bA)$.
Therefore
\begin{equation} \label{eqn-graded} 
-[[\bar{\beta}_1],[\bar{x}_2]] - [[\bar{\alpha}_1],[\bar{y}_2]] +
[[\bar{\alpha}_2],[\bar{y}_1]] = 0 \in \gr(H\bA).
\end{equation}

Comparing \eqref{eqn-ungraded} and \eqref{eqn-graded} we see that the
multiplicative structure in $H\bA$ and $\gr(H\bA)$ is not the same.

\chapter{Basic Topology and Adams-Hilton Models} \label{chapter-ah}

In this chapter we review some basic topological results and Adams-Hilton
models.
We also prove some results on Adams-Hilton models which we will need
in Chapter~\ref{chapter-hurewicz}.

\section {Basic topology} \label{section-basic-top}

We will work in the usual category studied in algebraic topology, that of
\emph{compactly generated} topological
spaces~\cite{may:bookCCAlgTop, selick:book, spanier:book, whitehead:book}.
We will assume that all of our spaces $X$ are \emph{simply-connected}.
That is, $X$ is path-connected and has trivial fundamental group ($\pi_1(X)
= 0$) and that our spaces are \emph{pointed} (also called
\emph{based}).
That is, they come with a chosen point called the
\emph{basepoint} and usually denoted $*$. 
The pointed space $(X,*)$ will usually be referred to as just $X$.
All our maps are assumed to be continuous and pointed.
That is, $f:(X,*) \to (Y,*)$ satisfies $f(*)=*$.
Furthermore we will assume that all of our spaces have `the weak
homotopy type of a finite-type CW complex.' 
We explain this statement below.

Whenever we take products of infinitely many spaces we will always
mean the \emph{weak infinite product}.
That is, $x \in \prod_i X_i$ implies that $x_i = *$ for all but finitely
many $i$.
A \emph{weak product} is a finite product or a weak infinite product.

Given spaces $X$ and $Y$ define the \emph{wedge} $X \vee Y$ to be the
space obtained by attaching $X$ and $Y$ at the basepoints.
That is, $X \vee Y = \{ (x,y) \in X \times Y \ | \ x=* \text{ or } y=*
\}$ with basepoint $(*,*)$.  
Let an \emph{$(n+1)$-cell} be the unit disc $e^{n+1} = \{ x \in \R^{n+1} \ | \
\|x\| \leq 1 \}$ and let the \emph{$n$-sphere} be its boundary $S^n = \{
x \in \R^{n+1} \ | \ \|x\| = 1 \}$.

Let $W$ be a subspace of $Z$ and $f:W \to X$.
Then the \emph{attaching map construction} builds a new space $Y = (Z
\amalg X) / \! \sim$ where $f(w) \sim w, \ \forall w \in W$.
This space is called an \emph{adjunction space} and we
denote\footnote{
In the literature this space is sometimes also denoted by $Z \cup_f
X$.
}  
it by $X \cup_f Z$ and call $f$ the \emph{attaching map}. 
For example if $Z = \bigvee_{j \in J} e^{n_j+1}$ and $W = \bigvee_{j \in J}
S^{n_j}$ then we say that the space $Y$ is obtained from $X$ by
\emph{attaching cells along $W$}. 
In this case $W \to X \to Y$ is a \emph{homotopy cofibration}
\cite{may:bookCCAlgTop, selick:book, spanier:book, whitehead:book}.
A cofibration $W \to X$ is an inclusion and given any homotopy $H: W \times I
\to A$ and an extension $g:X \to A$ of $H|_{W \times \{0\}}$ there
exists an extension $H\p: X \times I \to A$ of $H$ such that $H\p|_{X
\times \{0\}} = g$.

A \emph{finite-type CW complex} is any space $X$ that can be constructed by
the following inductive procedure.
Let $X^{(0)}$ be a finite set of points. 
Obtain $X^{(n+1)}$ from $X^{(n)}$ by attaching a finite number of
$(n+1)$-cells along their boundary.
Let $X = \cup_n X^{(n)}$.
Call $X^{(n)}$ the \emph{$n$-skeleton} of $X$. 
Call a map $f:X \to Y$ between $CW$-complexes \emph{cellular} if
$f:X^{(n)} \to Y^{(n)}, \ \forall n$.

Let $I$ denote the unit interval $[0,1]$.
Two (pointed) maps $f,g: X \to Y$ are said to be \emph{homotopy equivalent},
written $f \simeq g$ if there is a map (called a \emph{homotopy}) $H:
X \times I \to Y$ such that $H(x,0) = f(x)$, $H(x,1) = g(x)$ and
$H(*,t) = *, \ \forall x \in X$ and $\forall t \in I$.
This is an equivalence relation.
Let $X$ and $Y$ be two topological spaces. 
If there exist maps $f:X \to Y$ and $g:Y \to X$ such that $g
\circ f \simeq id_X$ and $f \circ g \simeq id_Y$ then we say $X$ is
\emph{homotopy equivalent} to $Y$ or $X$ has the same \emph{homotopy
type} as $Y$, denoted $X \approx Y$.  
This is also an equivalence relation.

Let $\pi_n(X)$ be the set of homotopy classes of pointed maps from
$S^n \to X$. 
For $n \geq 1$ this is a group \cite{selick:book, may:bookCCAlgTop,
spanier:book, whitehead:book}. 

A map $f:X \to Y$ is called a \emph{weak homotopy equivalence} if
$f_{\#}:\pi_n{X} \to \pi_n(Y)$ is an isomorphism for all $n$.
Topological spaces $X$ and $Y$ have the same \emph{weak homotopy type}
if there exists a chain of weak homotopy equivalences  $X \leftarrow Z_1
\to \cdots \leftarrow Z_n \to Y$.
The following theorem shows that any topological space has the weak
homotopy type of a CW-complex.

\begin{thm}[The CW Approximation Theorem~{\cite[Theorem
V.2.2]{whitehead:book}}] 
Given a topological space $Y$ there is a CW-complex $X$ and a map $f:X
\to Y$ such that $f_{\#}: \pi_n(X) \to \pi_n(Y)$ is an isomorphism for
all $n$. 
\end{thm}

If $Y$ is simply-connected then one can choose $X$ such that $X^{(1)} = *$. 

In addition to the CW Approximation Theorem, the following two
theorems show why we can work in the category of CW-complexes and why
it is convenient to do so.

\begin{thm}[The Cellular Approx. Theorem {\cite[Theorem
II.4.5]{whitehead:book}}] 
Let $f:X \to Y$ be a map between CW-complexes. Then there is a
cellular map $g: X \to Y$ such that $f \simeq g$.
\end{thm}

\begin{thm}[The Whitehead Theorem for simply-connected CW-complexes
{\cite[Theorems IV.7.13 and V.3.5]{whitehead:book}}]  
\label{thm-whitehead} 
Let $X$ and $Y$ be simply-connected CW-complexes.
Then $X$ and $Y$ are homotopy equivalent if and only $X$ and $Y$ are
weak homotopy equivalent if and only if $\exists \ f: X \to Y$ such
that $f_*: H_n(X) \to H_n(Y)$ is an isomorphism $\forall n$.
\end{thm}

Given a (based) topological space $(X,*)$ let $\lX$ denote the set of
(based) loops in $X$.
That is, $\lX = \Map_*( (S^1,*) \to (X,*) )$, where
$\Map_*(X,Y)$ denotes the pointed maps from $X$ to $Y$ with the
\emph{compact-open} topology (see~\cite{selick:book, may:bookCCAlgTop,
spanier:book, whitehead:book}).
The constant map $(S^1,*) \xto{*} (*,*)$ is the basepoint of $\lX$.
$\lX$ is called the \emph{loop space} on $X$ and $H_*(\lX;R)$ is
called the \emph{loop space homology} of $X$.
A key property of $\lX$ is the following.

\begin{lemma}[{\cite[Corollary IV.8.6]{whitehead:book}}]
$\forall n \geq 1, \ \pi_n(\lX) \isom \pi_{n+1}(X)$.
\end{lemma}


Thus we can study the homotopy groups of $X$ by studying $\lX$.
From our point of view the main advantage of studying $\lX$ is that
the homology groups of $\lX$ have the structure of a Hopf algebra
whereas the homology groups of $X$ do not have the structure of an algebra.
The algebra structure of $H_*(\lX)$ can reveal important information
about the homotopy type of $X$.

There is a natural map 
\begin{equation} \label{eqn-hurewiczmapR} 
h_X: \pi_*(\lX) \tensor R \to H_*(\lX;R)
\end{equation} 
called the \emph{Hurewicz map} and defined as follows.
For $\alpha \in \pi_n(\lX) \tensor R$ choose a representative $a:S^n
\to \lX$.
Then there is an induced map $a_*: H_*(S^n;R) \to H_*(\lX;R)$.
Let $h_X(\alpha) = a_*(\iota_n)$, where $\iota_n$ is a generator for
$H_n(S^n)$.
One can check that this map is well defined and is
in fact a homomorphism.

Concatenation of loops induces an algebra structure on
$H_*(\lX;R)$~\cite{selick:book} called the \emph{Pontrjagin product}.
If $f:X \to Y$ then using this algebra structure on $H_*(\lX;R)$ and
$H_*(\lY;R)$ the induced map 
\begin{equation} \label{eqn-loopf*}
(\Omega f)_*: H_*(\lX;R) \to H_*(\lY;R)
\end{equation}
is an algebra map.
Using this algebra structure, the commutator bracket $[a,b] = ab -
(-1)^{|a||b|}ba$ gives $H_*(\lX;R)$ the structure of a Lie algebra.

Given $\alpha \in \pi_m(\lX) \tensor R$ and $\beta \in \pi_n(\lX)
\tensor R$ define $[\alpha,\beta] \in \pi_{m+n}(\lX) \tensor R$,
called the \emph{Samelson product}, as follows.
Let $f:S^m \to \lX$ and $g:S^n \to \lX$ be representatives of $\alpha$
and $\beta$.
Define $h: S^m \times S^n \to \lX$ by letting $h(x,y)$ be the loop
obtained by concatenating the loops $f(x), \ g(y), \ -f(x)$ and
$-g(y)$, where the negative of a loop is the loop traced in reverse.
Since $h|_{S^m \vee S^n}$ is contractible one can show that there is
an induced map $S^{m+n} \to \lX$.
Let $[\alpha,\beta]$ be the homotopy class of this map.
One can check that the Samelson product gives $\pi_*(\lX) \tensor R$ the
structure of a Lie algebra~\cite{selick:book}.

Using these products the Hurewicz map $h_X: \pi_*(\lX) \tensor R \to
H_*(\lX;R)$ is a Lie algebra map~\cite{samelson:products}.
Let $L_X$ denote the image of $h_X$ and call this Lie subalgebra of
the loop space homology the \emph{Hurewicz images}.
A map $f: W \to X$ induces a map $L_W \to L_X$. 
Denote the image of this map by $L^W_X$, and let $[L^W_X]$ be the Lie
ideal in $L_X$ generated by $L^W_X$.

Given $\alpha \in \pi_m(X)$ and $\beta \in \pi_n(X)$ there is a
\emph{Whitehead product} $[\alpha, \beta] \in \pi_{m+n-1}(X)$ which can
be defined as follows. 
Let $\omega_{m,n}: S^{m+n-1} \to S^m \vee
S^n$ be the attaching map of the top cell in $S^m \times S^n$.
An explicit description is given in~\cite{fht:rht}.
$[\alpha, \beta]$ is then the homotopy class represented by the
composite $S^{m+n-1} \xto{\omega_{m,n}} S^m \vee S^n \xto{a,b} X$
where $a$ and $b$ are representatives of $\alpha$ and $\beta$.
Up to sign, the Whitehead product can also be defined using the
\emph{adjoint} (defined below) of the Samelson product~\cite{selick:book}. 

Assume $R \subset \Q$ is a subring containing $\frac{1}{6}$ and assume
$H_*(\lX;R)$ is torsion-free.
Let $P$ be the set of invertible primes in $R$, and let $\nP = \{p \in
\Z \ | \ p \text{ is prime and } p \notin P\} \cup \{0\}$.
Let $F_0 = \Q$.
Then $\fpnP$ we have the commutative diagram
\[
\xymatrix{ 
\pi_*(\lX) \tensor R \ar[r]^{h_X} \ar[d]^{\_ \tensor \Fp} & H_*(\lX;R)
\ar[d]^{\_ \tensor \Fp} \\
\pi_*(\lX) \tensor \Fp \ar[r]^-{h_X \tensor \Fp} & H_*(\lX;R) \tensor
\Fp \ar[r]^-{\isom} & H_*(\lX;\Fp)
}
\] 
where the bottom right map is the isomorphism given by the Universal
Coefficient Theorem. 
Abusing notation we will refer to the composition of the bottom two
maps as $h_X \tensor \Fp$ and refer to its image as $L_X \tensor \Fp$.
It is easy to check that this is the same as the map $h_X: \pi_*(\lX)
\tensor \Fp \to H_*(\lX;\Fp)$ defined in~\eqref{eqn-hurewiczmapR} when
$R = \Fp$.
We will sometimes denote $h_X \tensor \Fp$ and $L_X \tensor \Fp$ by
$\bar{h}_X$ and $\bar{L}_X$ omitting $p$ from the notation.

If $f: W \to X$ and $H_*(\Omega W;R)$ and $H_*(\lX;R)$ are $R$-free
then $\fpnP$ there is an induced map $\bar{L}_W \to \bar{L}_X$.
Denote the image of this map by $\bar{L}^W_X$, and let $[\bar{L}^W_X]$
be the Lie ideal in $\bar{L}_X$ generated by $\bar{L}^W_X$.

When $R = \Q$ Milnor and Moore~\cite{milnorMoore:hopfAlgebras} proved
the following major result about the rational Hurewicz map.

\begin{thm}[The Milnor-Moore Theorem \cite{milnorMoore:hopfAlgebras}]
\label{thm-milnorMoore} 
The rational Hurewicz homomorphism $h_X: \pi_*(\lX) \tensor
\Q \to H_*(\lX;\Q)$ is an injection, and furthermore, as algebras
$H_*(\lX;\Q) \isom UL_X$.
\end{thm}

In Theorems \ref{thm-hurewiczF} and \ref{thm-hurewiczR} we prove
versions of this theorem for more general coefficient rings under
certain hypotheses.  
See~\cite{scott:tfmmPreprint} for another extension of the
Milnor-Moore Theorem.

Given a space $X$ define the \emph{(reduced) suspension} $\Sigma X$,
as follows. 
Let $\Sigma X = (X \times I) / \! \sim$ where $(a,0)\sim(b,0), \
(a,1)\sim(b,1)$ and $(*,s)\sim(*,t), \ \forall a,b \in X$ and $\forall
s,t \in I$. 
The basepoint is the equivalence class of $(*,0)$.
For example $\Sigma S^n \isom S^{n+1}$ and $\Sigma\left(\bigvee_j
S^{n_j}\right) \approx \bigvee_j S^{n_j+1}$. 

Given a map $g: \Sigma X \to Y$ there is an \emph{adjoint} map
$\hat{g}: X \to \lY$ defined by $\hat{g}(x)(t) = g(x,t)$ for $x \in X$
and $t \in I$. 
Let $\alpha: X \to \Omega \Sigma X$ be the adjoint of $id_{\Sigma X}$.

\begin{thm}[The Bott-Samelson Theorem~\cite{bottSamelson:thm}]
\label{thm-bott-sam} 
Let $R$ be a principal ideal domain and let $X$ be a connected space
such that $H_*(X;R)$ is torsion-free. 
Then 
\[
H_*(\Omega \Sigma X;R) \isom \freeT( \tilde{H}_*(X;R)).
\]
Furthermore $\alpha: X \to \Omega \Sigma X$ induces the canonical
inclusion $\tilde{H}_*(X;R) \incl \freeT (\tilde{H}_*(X;R))$.
\end{thm}

\begin{eg}
$H_*(\Omega S^{n+1};R) \isom \freeT \langle x \rangle \isom U\freeL
\langle x \rangle$ where $|x| = n$.
Furthermore $x = h_{S^{n+1}}([\alpha])$ where $\alpha: S^n \to
\Omega S^{n+1}$ is defined above. 
Therefore $H_*(\Omega S^{n+1};R) \isom UL_{S^{n+1}}$.
\end{eg}

\begin{eg} \label{eg-wedge-of-spheres}
$H_*(\Omega (\bigvee_{j \in J} S^{n_j+1}); R)
\isom \freeT V \isom U 
\freeL V$ where $V$ is a free $R$-module with basis $\{x_j\}_{j \in
J}$ and $|x_j|=n_j$. 
Let $\iota_j: S^{n_j} \to \bigvee_{j \in J} S^{n_j}$ denote the inclusion
of one of the spheres.
Then $x_j = h_{\bigvee S^{n_j+1}} ([\alpha \circ \iota_j])$ and hence
$H_*(\Omega (\bigvee_{j \in J} S^{n_j+1});R ) \isom
UL_{\bigvee S^{n_j+1}}$.   
\end{eg}

In Section~\ref{section-hsb} we will need the following \emph{infinite
mapping telescope} construction for the \emph{localization} of
CW-complexes.
Let 
\[ X_1 \xto{j_1} X_2 \xto{j_2} \ldots X_n \xto{j_n} \ldots
\]
be a sequence of maps between topological spaces.
The infinite mapping telescope of this sequence is the space
\[ T = \left. \left(\coprod_{n=1}^{\infty} (X_n \times [n-1,n])
\right) \right/ \sim \text{ where } (x_n \times \{n\}) \sim (j_n (x_n)
\times \{n\}). 
\]

\section{Adams-Hilton models} \label{section-ah}

Let $R$ be a principal ideal domain containing $\frac{1}{6}$.
Given a differential graded algebra (dga) $A$, a dga morphism $A\p \to
A$ is called a \emph{model} for $A$ if it induces an isomorphism on
homology.
If $A\p \isom (TV,d)$ then it is called a \emph{free model}.
A \emph{model} for a simply-connected topological space $Z$ means a model
for $CU_*^1(\Omega Z)$ where $CU_*^1(\,)$ is the first Eilenberg subcomplex
of the cubical singular chain complex.
If such a model is free it is called an \emph{Adams-Hilton model
(AH-model)}~\cite{adamsHilton}.

Every simply-connected topological space $Z$ has a (non-unique)
AH-model 
\[
\theta_Z: \AH{Z} \stackrel{\simeq}{\rightarrow} CU_*^1(\Omega Z).
\]
This induces an isomorphism of algebras $H_*(\AH{Z})
\stackrel{\cong}{\rightarrow} H_*(\Omega Z; R)$. 
We will usually denote an Adams-Hilton model by just $\AH{Z}$.
We state some basic properties of these models.

For any map $f:X\rightarrow Y$ between simply-connected CW-complexes
and any choice of Adams-Hilton models $(\AH{X},\theta_X)$ and
$(\AH{Y},\theta_Y)$ there is a dga homomorphism 
$\AH{f}: \AH{X} \rightarrow \AH{Y}$ which comes with a dga homotopy $\psi_f$
from $CU_*^1(\Omega f) \circ \theta_X$ to $\theta_Y \circ \AH{f}$.

If $X$ is a finite-type CW-complex then a CW-structure on $X$
determines an Adams-Hilton model on $X$ as follows.
The CW-structure gives a sequence of cofibrations 
\[ \bigvee_{j \in J_n} S^n_j \xto{\bigvee \alpha_{n,j}} X^{(n)} \to X^{(n+1)}
\]
where $X^{(n)}$ is the $n$-skeleton of $X$ and each index set $J_n$ is
finite.
Let $\{x_{n,j}\}_{j \in J_n}$ be a set of graded elements with
$|x_{n,j}| = n$.
Let $V_n = R\{x_{n,j}\}_{j \in J_n}$ and let $V = \bigoplus_n V_n$.
One can choose $(\freeT V,d)$ as an AH-model for
$X$, where the differential is defined inductively and satisfies
$dx_{n,j} \subset Z (\freeT \left(\bigoplus_{i=1}^{n-1} V_1\right), d)$.

Let $f: W \to X$ be a map between finite-type simply-connected
CW-complexes where $W = \bigvee_{j \in J} S^{n_j}$, $f = \bigvee_{j
\in J} \alpha_j$, $H_*(\lX;R)$ is torsion-free and $H_*(\lX;R) \isom
UL_X$ (recall $L_X$ from Section~\ref{section-basic-top}). 
Let $Y$ be the adjunction space $X \cup_f \left(\bigvee_{j \in J}
e^{n_j+1}\right)$.
Then one can take 
\begin{equation} \label{eqn-L(Y)} 
\AH{Y} = \AH{X} \amalg \freeT \langle y_j \rangle_{j \in J}, 
\end{equation}
where $\freeT \langle y_j \rangle_{j \in J}$ is the tensor algebra on the 
free $R$-module $R \{ y_j \}_{j \in J}$
and $A \amalg B$ is the coproduct or free product of $A$ and $B$.
The differential on $\AH{Y}$ is an extension of the differential on
$\AH{X}$ satisfying $dy_j \in \AH{X}$.
Furthermore if $d'$ is the induced map 
\[ d': R\{y_j\} \to Z\AH{X} \onto H\AH{X} \isomto H_*(\lX;R) \isomto
UL_X
\]
then $d'y_j = h_X(\alpha_j) \in L_X$.
That is, $\AH{Y}$ is a dga extension (see
Section~\ref{section-inert-extns}) of $(\AH{X},L_X)$.

We can filter $\AH{Y}$ as a dga by taking $F_{-1}\AH{Y} = 0$,
$F_0\AH{Y} = \AH{X}$ and 
for $i \geq 0$, $F_{i+1}\AH{Y} = \sum_{j=0}^i F_{j}\AH{Y}
\cdot R\{y_j\}_{j \in J} \cdot F_{i-j}\AH{Y}$.
This filtration induces a filtration on $H\AH{Y}$.
From the filtration on $\AH{Y}$ we get a first quadrant spectral
sequence $\gr(\AH{Y}) \implies \gr(H\AH{Y})$ (see
Section~\ref{section-ss}).

\ignore
{

\section{Anick models}

Let $R = \Fp$ with $p>3$ or let $R \subset \Q$ be a subring containing
$\frac{1}{6}$. 
Let $\rho$ be the least non-invertible prime in $R$, taking
$\rho=\infty$ if $R$ is a field.

Let $Z$ be a $r$-connected (with $r\geq 1$) topological space of
dimension $\leq n$. 
Then if $\rho \geq \frac{n}{r}$, then $Z$ has an Adams-Hilton
model $\theta_Z: U\LL{Z} \to CU^r_*(\lZ)$ where $\LL{Z}$ is
differential graded Lie algebra called an \emph{Anick
model}~\cite{anick:hah} for $Z$. 
If $R$ is a field it has the additional property that $H_*(\theta_Z;R)$
is a Hopf algebra isomorphism.

Let $f: W \to X$ be a map between finite-type simply-connected
CW-complexes where $W = \bigvee_{j \in J} S^{n_j}$, $f = \bigvee_{j
\in J} \alpha_j$ and $H_*(\lX;R) \isom UL_X$ (recall $L_X$ from
Section~\ref{section-basic-top}). 
Let $Y$ be the adjunction space $X \cup_f \left(\bigvee_{j \in J}
e^{n_j+1}\right)$.

Assume that $X$ and $Y$ are $r$-connected (with $r\geq 1$) and of
dimension $\leq n$. 
Then if $\rho \geq \frac{n}{r}$, $X$ and $Y$ have Anick models
$\LL{X}$ and $\LL{Y}$ satisfying 
\begin{equation}
\LL{Y} = \LL{X} \amalg \freeL \langle y_j \rangle_{j \in J}, 
\end{equation}
where $\freeL$ denotes the free Lie algebra.
The differential on $\LL{Y}$ is an extension of the differential on
$\LL{X}$ satisfying $dy_j \in \LL{X}$.
Furthermore if $d'$ is the induced map 
\[ d': R\{y_j\} \to Z\LL{X} \incl ZU\LL{X} \onto HU\LL{X} \isomto
H_*(\lX;R) \isomto UL_X
\]
then $d'y_j = h_X(\alpha_j) \in L_X$.
That is, $\LL{Y}$ is a dgL extension of $(\LL{X},L_X)$. 

So if $\rho \geq \frac{n}{r}$ then one can use Anick models instead of
Adams-Hilton models.
As a result one has dgL extensions instead of just dga extensions and
if $R$ is a field then one gets Hopf algebra isomorphisms.

}

\section{Some lemmas}

We will prove some results using Adams-Hilton models which we will
need in Chapter~\ref{chapter-hurewicz}.
Let $R$ be a principal ideal domain containing $\frac{1}{6}$.

First we need the following Adams-Hilton models of some standard spaces.
Recall from Section~\ref{section-ah} that there is an Adams-Hilton
model corresponding to a CW structure.
Consider $S^m$ with the CW structure $S^m_0 \cong * \cup e^m$. 
Its corresponding AH model is $(\freeT \langle a \rangle, 0)$ where
$|a|=m-1$.
Attach another cell to get the disk $D^{m+1}_0 \cong S^m_0 \cup e^{m+1}$.
Its corresponding AH model is $(\freeT \langle a,b \rangle, d)$
where $|a|=m-1, \ |b|=m, \ da=0$ and $db=a$.
We fix the AH model for the  inclusion 
$\iota_0: S^m_0 \hookrightarrow D^{m+1}_0$ to be the homomorphism
defined by $\AH{\iota_0}(a) = a$.

For any $s$, the sphere also has the following more complicated CW structure
\[ S^{m+1}_{s} \isom \left(\bigvee_{k=1}^s D^{m+1}_0\right)
\cup_{\sum_{k=1}^s \iota_k} e^{m+1} \]
where the attaching map for the last cell is given by the inclusions 
of the spheres into $\bigvee_{k=1}^{s} S^m_0$.  
Its corresponding AH model is $(\freeT \langle a_1, \ldots a_s, b_0, \ldots b_s
\rangle, d)$ where $|a_i|=m-1 \ |b_i| = m, \ da_i=0, \ db_i = a_i$
for $1\leq i\leq s$, and $db_0 = a_1 + \ldots + a_s$.
Furthermore there is a homeomorphism 
\begin{equation} \label{eqn-Psi}
\Psi: S_0^{m+1} \stackrel{\cong}{\rightarrow} S_s^{m+1}
\end{equation}
which has an Adams-Hilton model that sends $a$ to $b_1 + \ldots + b_s - b_0$.

Recall from Section~\ref{section-basic-top} that the Whitehead
product is defined using a map  
$\omega_{m,n}: S^{m+n-1}_0 \to S^m_0 \vee S^n_0$.
If we take $(\freeT \langle a_{m+n-1} \rangle,0)$ and $(\freeT \langle
a_m, a_n \rangle, 0)$ to be the Adams-Hilton models of those spaces then
by~\cite[Corollary 2.4]{adamsHilton} 
we may take $\AH{\omega_{m,n}}(a_{m+n-1}) = \pm [a_m, a_n]$. 
We will make use of the following fact.

\begin{lemma} \label{lemma-whiteheadextension}
\cite[Lemma 5.2]{anick:cat2}
Let $(\freeT \langle a,b \rangle, d)$ with $db=a$ be an AH model for
$D_0^{m+n}$. 
Let $(\freeT \langle a\p,b\p,c\p \rangle, d)$ be an AH model for
$D_0^{m+1} \vee S_0^n$ where $c\p$ corresponds to the $S^n_0$ and
$d(b\p) = a\p$. 
Then there is an extension of $\pm \omega_{m,n}: S_0^{m+n-1} \to S_0^m \vee
S_0^n$ to a map $f: D_0^{m+n} \to D_0^{m+1} \vee S_0^{n+1}$ whose AH
model may be chosen so as to satisfy
\[ 
\AH{f}(a)=[a\p,c\p], \;\; \AH{f}(b)=[b\p,c\p], \;\; \text{and} \;\;
\psi_f(a) \in CU_*(\Omega( S_0^m \vee S_0^n)).
\]
\end{lemma}

From this we prove the following lemma which we will need in
Chapter~\ref{chapter-hurewicz} to construct a map $S^{m+1} \to Y$. 
It is a slight generalization of \cite[Lemma 5.3]{anick:cat2} and we copy
Anick's proof.

\begin{lemma} \label{lemma-ah}
Suppose we are given a simply-connected space $X$ for which there
exists a Lie algebra map $\sigma_X$ right inverse to $h_X$. 
In addition we are given $c \in R$, a map $\alpha: S^{n+1}_0
\rightarrow X$, and $x_i \in Z\AH{X}$, for $i=1,\ldots t$ such that
$[x_i] \in L_X$ and $\beta_i: S^{n_i+1} \rightarrow X$ is the adjoint
of $\sigma_X ([x_i])$. 
Choose an Adams-Hilton model $\AH{\alpha}$.
Let $z = \AH{\alpha} (a) \in \AH{X}$ and 
let $Y = X \cup_\alpha e^{n+2}$.
Then we can choose
$\AH{Y} = \AH{X} \amalg \freeT \langle y \rangle$ with $dy = z$.
In addition there exists a map $g: (D_0^{m+1},S_0^m) \rightarrow (Y, X)$
such that $g|_{S^m_0} = \epsilon c [[\alpha, \beta_1, \ldots \beta_t]$, where
$\epsilon=\pm 1$.
Furthermore \AH{g}\!\! may be chosen so that 
$\AH{g}(a) = c[[z, x_1, \ldots x_t]$ and 
$\AH{g}(b) = c[[y, x_1, \ldots x_t]$ and the dga homotopy 
$\psi_g(a)$ lies in the submodule $CU_*(\lX) \subset CU_*(\lY)$. 
\end{lemma} 

\begin{proof}
The proof is by induction on $t$, the length of the list of indices.
When $t=0$, \cite[Theorem 3.2]{adamsHilton} tells us that the AH model
\AH{\alpha} may be extended over \AH{D_0^{n+2}} such that the
generators $a$ and $b$ are sent to $z$ and $y$.
Composing this with the degree $c$ map from $D_0^{n+2}$ to itself
gives the map $g$ for which $\AH{g}(a)=cz$ and $\AH{g}(b)=cy$.

For the inductive step, suppose the result to be true when the list
has $t-1$ elements.
Let $m\p = n +1 + n_1 + \ldots + n_{t-1}$.
Then the inductive hypothesis gives us a map
\[ g\p: (D_0^{m\p+1}, S_0^{m\p}) \to (Y,X) \]
satisfying  $\AH{g\p}(a\p) = c[[z, x_1, \ldots x_{t-1}]$,
$\AH{g\p}(b\p) = c[[y, x_1, \ldots x_{t-1}]$ and furthermore
$\psi_g(a\p) \in CU_*(\lX)$.
Define a map
\[ g\pp = g\p \vee \beta_t: (D_0^{m\p+1} \vee S_0^{n_t+1}, S_0^{m\p}
\vee S_0^{n_t+1}) \to (Y,X) \]
with $\AH{g\pp}$ an extension of $\AH{g\p}$ such that 
$\AH{g\pp}(c\p) = x_t$ where $c\p$ is the generator corresponding
to $S_0^{n_t+1}$.
Let $f$ denote the map of Lemma~\ref{lemma-whiteheadextension}, where
$m = m\p$ and $n=n_t+1$, and set $g = g\pp \circ f$.
Take $\AH{g} = \AH{g\pp} \circ \AH{f}$.
Then $g$ has the desired properties.
\end{proof}

\chapter{The Cell Attachment Problem} \label{chapter-cap}

In this chapter we review various topological constructions and
results. 
We also motivate and define our main topological objects of study:
\emph{cell attachments} which are \emph{free} and \emph{semi-inert}.

\section{Whitehead's cell-attachment problem} \label{section-cap}

One of the oldest questions in homotopy theory asks what effect 
attaching one or more cells has on the homotopy groups and loop space
homology groups of a space. 
This questions was perhaps first considered by
J.H.C. Whitehead \cite{whitehead:addingRelations},
\cite[Section~6]{whitehead:simplicialSpaces}, around 1940.

\noindent
\textbf{The cell attachment problem:}
Given a simply-connected topological space $X$ and a cofibration
\[ \bigvee_{j \in J} S^{n_j} \xto{f = \bigvee \alpha_j} X \xto{i} Y
\]
how is $H_*(\lY;R)$ related to $H_*(\lX;R)$ and how is $\pi_*(Y)$
related to $\pi_*(X)$?

We assume that $H_*(\lX;R)$ is torsion-free.
This condition is trivial if $R$ is a field.
For $R \subset \Q$ one can often reduce to this case by
\emph{localizing} (see Section~\ref{section-lsd}) away from a finite
set of primes. 
Even if the loop space homology of a given space has torsion at
infinitely many primes~\cite{anick:torsion, avramov:torsion} one might
be able to study the given space by including into a space $X$ such
that $H_*(\lX;R)$ is torsion-free~\cite{anick:cat2}.

The inclusion $i$ induces an algebra map $(\Omega i)_*: H_*(\lX;R) \to
H_*(\lY;R)$ (see \eqref{eqn-loopf*}). 
Recall (from Section~\ref{section-basic-top}) that $h_X$ denotes the
Hurewicz map and that $\hat{\alpha}_j$ is the adjoint of $\alpha_j$.
Let $W = \bigvee_{j \in J} S^{n_j}$, let $a_j = h_X(\hat{\alpha}_j)$ and let
$V_1 = R\{y_j\}_{j \in J}$ where $|y_j| = |a_j| + 1$.
Using the notation of Section~\ref{section-basic-top}, $L^W_X$ is the
Lie subalgebra of $L_X \subset H_*(\lX;R)$ generated by $\{a_j\}_{j
\in J}$.
Let $(L^W_X) \subset H_*(\lX;R)$ be the two-sided ideal generated by $L^W_X$.
Since $(\Omega i)_*$ is an algebra map $(\Omega i)_*((L^W_X)) = 0$ and
hence $(\Omega i)_*$ factors through the quotient map.
\begin{equation} \label{eqn-inert-am} 
\xymatrix{
H_*(\lX;R) \ar[dr] \ar[rr]^{(\Omega i)_*} & & H_*(\lY;R) \\
& H_*(\lX;R)/(L^W_X) \ar[ur]_g
}
\end{equation}
We will say that the attaching map $f$ is 
\emph{inert}~\cite{halperinLemaire:inert} if $H_*(\lX;R)/(L^W_X)$ is
torsion-free and $g$ is an isomorphism. 

\begin{rem} \label{rem-inertR-am}
Note that if $R \subset \Q$ and $f$ is inert then $H_*(\lY;R)$ is
$R$-free. 
Thus by the Universal Coefficient Theorem  $\fpnP$, $H_*(\lY;\Fp)
\isom H_*(\lX;\Fp) / (\bar{L}^W_X)$.
That is, $f$ is inert over $\Fp$ ($p$-inert in~\cite{hessLemaire:nice}).  
\end{rem}

\begin{eg}
$X = S^3_a \vee S^3_b$ and $Y = X \cup_f e^8$, where $f$ is the
iterated Whitehead product $[[\iota_a, \iota_b], \iota_b]$ with
$\iota_a$ and $\iota_b$ the inclusions of $S^3_a$ and $S^3_b$ in $X$.
\end{eg}

By the Bott-Samelson Theorem (Theorem~\ref{thm-bott-sam}) $H_*(\lX;R) \isom
\freeT \langle x,y \rangle$. 
Let $I$ be the two-sided ideal in $H_*(\lX;R)$ generated by the image
 of $f$.
That is, $I$ is the two-sided ideal generated by $[[x,y],y]$.

For an Adams-Hilton model of $Y$ we can take (see
Section~\ref{section-ah}) $U\LL{Y}$ where $\LL{Y} = (\freeL \langle
x,y,a \rangle, d), \ dx=dy=0$ and  $da=[[x,y],y]$.
So $H_*(\lY;R) \isom HU\LL{Y}$ as algebras.
By Example~\ref{eg-inert-dga} $HU\LL{Y} \isom U(\freeL \langle x,y
\rangle / [[x,y],y]) \isom U\freeL \langle x,y \rangle / ([[x,y],y])$
 as algebras. 
Thus $H_*(\lY;R) \isom H_*(\lX;R)/I$ as algebras and $f$ is an inert
attaching map.
\eolBox

If $R$ is a field then Halperin and
Lemaire~\cite{halperinLemaire:inert} showed that $g$ is surjective if
and only if it is an isomorphism.
In general the map $g$ in \eqref{eqn-inert-am} need not be injective
or surjective.
The attaching map $f$ is said to be
\emph{nice}~\cite{hessLemaire:nice} if the map $g$~\eqref{eqn-inert-am}
is injective. 
As we will see in Example~\ref{eg-fat-wedge} $f$ can be nice but not
inert.

The following theorem gives an equivalent characterization of the
inert condition.
We will use it as motivation to generalize the inert condition.
It is proved by Halperin and Lemaire~\cite[Theorem
3.3]{halperinLemaire:inert} if $R = \Q$, and F\'{e}lix and
Thomas~\cite[Theorem 1]{felixThomas:attach} if $R$ is a field of
characteristic $\neq 2$.  

\begin{thm}[\cite{halperinLemaire:inert, felixThomas:attach}]
\label{thm-inert-am} 
Let $J = [L^W_X] \subset L_X$ be the Lie ideal of $L_X$ generated by
$L^W_X$.
If $R=\Q$ or $\Fp$ where $p>3$ then $f$ is an inert attaching map if and
only if \\ 
(i) $J$ is a free Lie algebra, and \\
(ii) $J/[J,J]$ is a free $U(L_X/J)$-module.
\end{thm}

We will study attaching maps which in general do not satisfy the second
condition. 
If $R = \Q$ or $\Fp$ we define $f$ to be a \emph{free} attaching map
if $[L^W_X] \subset L_X$ is a free Lie algebra.
If $R \subset \Q$ we define the attaching map $f$ to be \emph{free}
if \feachpnP, $[\bar{L}^W_X] \subset \bar{L}_X$ is a free Lie algebra.

This condition is a broad generalization of the inert condition and we
will prove results about $H_*(\lY;R)$ if $f$ is a free cell attachment
(see part (i) of Theorems \ref{thm-c} and \ref{thm-d}).
However many free cell attachments satisfy a second condition.
We define this condition below and call it the \emph{semi-inert}
condition.
Under this condition we are able to determine $H_*(\lY;R)$ as an algebra
(see part (ii) of Theorems \ref{thm-c} and \ref{thm-d}).
We will show that it too is a generalization of the inert condition.

We will use the following filtration on $H_*(\lY;R)$.
Recall from Section~\ref{section-ah} that $X$ and $Y$ have
Adams-Hilton models $\AH{X}$ and $\AH{Y}$ satisfying 
$\AH{Y} = \AH{X} \amalg \freeT V_1$ where $V_1 = R\{ y_j \}_{j \in
J}$.
$\AH{Y}$ is filtered by taking $F_{-1}\AH{Y} = 0$, $F_0 \AH{Y} =
\AH{X}$ and for $n\geq 0$, $F_{n+1} \AH{Y} = \sum_{i=0}^n F_i \AH{Y}
\cdot V_1 \cdot F_{n-i} \AH{Y}$. 
Since $dF_n \AH{Y} \subset F_n \AH{Y}$ and $F_n \AH{Y} \cdot F_m
\AH{Y} \subset F_{n+m} \AH{Y}$ this makes $\AH{Y}$ into a filtered dga.
There is an induced dga filtration on $H\AH{Y}$.
Since $H_*(\lY;R) \isom H\AH{Y}$ there is a corresponding filtration
on $H_*(\lY;R)$.
Let $\gr_*(H_*(\lY;R))$ be the associated graded object.
For degree reasons $\gr_1(H_*(\lY;R))$ is a $\gr_0(H_*(\lY;R))$-bimodule.

Let $\eL = (L_X \amalg \freeL V_1, d\p)$ where $d'y_j = h_X(\hat{\alpha}_j)$.
Note that $[dV_1] = [L^W_X]$.
If we let $L_X$ be in degree $0$ and $V_1$ be in degree $1$ then $\eL$
is a bigraded dgL and $H\eL$ is a bigraded Lie algebra.
Let $(H\eL)_i$ denote the component of $H\eL$ in degree $i$. 
Then for degree reasons $(H\eL)_1$ is a $(H\eL)_0$-module.
Thus by Remark~\ref{rem-L0-module} there exists a semi-direct
product $(H\eL)_0 \sdp \freeL (H\eL)_1$.

We will use the following lemma to define the semi-inert condition.
Since this condition is not used in the statements or the proofs of
Theorem~\ref{thm-c}(i) and Theorem~\ref{thm-d}(i) we will take the
liberty of using these results. 

\begin{lemma} \label{lemma-semi-inert-attach}
Let $f:W \to X$ be a free cell attachment and let $Y$ be the
corresponding adjunction space. 
Then the following conditions are equivalent: \\
(a) $(H\eL)_0 \sdp \freeL (H\eL)_1 \isom (H\eL)_0 \amalg \freeL K$ 
for some free $R$-module $K \subset (H\eL)_1$, \\
(b) $(H\eL)_1$ is a free $(H\eL)_0$-module, and \\
(c) $\gr_1(H_*(\lY;R))$ is a free $\gr_0(H_*(\lY;R))$-bimodule.
\end{lemma}

\begin{proof}
The proof of this lemma is the same as the proof of
Lemma~\ref{lemma-semi-inert-dga}. \\
(b) $\implies$ (a)
Let $K$ be a basis for $(H\eL)_1$ as a free $(H\eL)_0$-module.
Then $(H\eL)_0 \sdp \freeL (H\eL)_1 \isom (H\eL)_0 \amalg \freeL
K$. \\
(a) $\implies$ (c)
Since $f$ is a free cell attachment, by Theorem~\ref{thm-c}(i) or
Theorem~\ref{thm-d}(i),  
$\gr_*(H_*(\lY;R)) \isom U \left( (H\eL)_0 \sdp \freeL (H\eL)_1 \right)$.
So by (a), 
\[
\gr_*(H_*(\lY;R)) \isom U \left( (H\eL)_0 \amalg \freeL K \right) \isom
\gr_0(H_*(\lY;R)) \amalg \freeT K,
\]
for some free $R$-module $K \subset (H\eL)_1$. 
Therefore 
\[
\gr_1(H_*(\lY;R)) \isom \left[ \gr_0(H_*(\lY;R)) \amalg \freeT K
\right]_1 \isom \gr_0(H_*(\lY;R)) \tensor K \tensor \gr_0(H_*(\lY;R)).
\]
(c) $\implies$ (b)
Let $L' = (H\eL)_0 \sdp \freeL (H\eL)_1$.
Then by Theorem~\ref{thm-c}(i) or Theorem~\ref{thm-d}(i), $\gr_*(H_*(\lY;R))
\isom UL'$ and $\gr_1(H_*(\lY;R)) \isom (UL')_1$.
By (c), $(UL')_1$ is a free $(UL')_0$-bimodule.
We claim that it follows that $L'_1$ is a free $L'_0$-module.
Indeed, if there is a nontrivial degree one relation in $L'$ then
there is a corresponding nontrivial degree one relation in $UL'$.
\end{proof}

We say that a cell attachment $f$ is \emph{semi-inert} if $f$ is
free and it satisfies the conditions of the previous lemma.
%
We justify this terminology with the following lemma.

\begin{lemma}
An inert attaching map is semi-inert.
\end{lemma}

\begin{proof}
Let $f$ be an inert attaching map.
If $R = \Q$ or $\Fp$ then by Theorem~\ref{thm-inert-am} $[dV_1]$ is a
free Lie algebra.
If $R \subset \Q$ then by Remark~\ref{rem-inertR-am} $\fpnP$, $f \tensor
\Fp$ is $p$-inert.  
So by Theorem~\ref{thm-inert-am} \feachpnP, $[\bar{d}\bar{V}_1]$ is a free
Lie algebra.
So in either case $f$ is free.

Since $(H\eL)_1 = 0$, the semi-inert condition is trivially satisfied.
\end{proof}

\begin{rem}
For two-cones there is a nice equivalent condition to the semi-inert
condition.
It is known that the attaching map of a two cone $Y$ is inert iff
$\gldim H_*(\lY) \leq 2$.
Using the main result of \cite{felixLemaire:2level} it follows
directly that the attaching map of a two cone $Y$ is semi-inert iff
$\gldim H_*(\lY) \leq 3$.
\end{rem}

\begin{eg} \label{eg-semi-inert-top} 
Let $X = S^3_a \vee S^3_b$ and let $\iota_a, \ \iota_b$ denote the
inclusions of the spheres into $X$.
Let $Y = X \cup_{\alpha_1 \vee \alpha_2} (e^8 \vee e^8)$ where the
attaching maps are given by the iterated Whitehead products $\alpha_1
= [[\iota_a, \iota_b], \iota_a]$ and $\alpha_2 = [[\iota_a,
\iota_b],\iota_b]$.
\end{eg}

Let $W = S^7 \vee S^7$ and for $i=1,2$ let $\hat{\alpha}_i: S^6 \to
\lX$ denote the adjoint of $\alpha_i$.
Then $W \xto{\alpha_1 \vee \alpha_2} X \to Y$ is a homotopy cofibration and
$[L^W_X] = [h_X(\hat{\alpha}_1),h_X(\hat{\alpha}_2)]$.

$Y$ has an Adams-Hilton model (see Section~\ref{section-ah}) $U(L,d)$
where $L = \freeL \langle x,y,a,b \rangle, \ |x|=|y|=2, \ dx=dy=0, \
da=[[x,y],x]$ and $db=[[x,y],y]$.
Furthermore $h_X(\hat{\alpha}_1) = [da]$ and
$h_X(\hat{\alpha}_2)=[db]$.

By Example~\ref{eg-semi-inert-extn}, as algebras 
\[ HU(L,d) \isom \left( H_*(\lX;R)/ (L^W_X) \right) \amalg \freeT \langle [w]
\rangle,
\]
where $w = [a,y] - [b,x]$.
Thus $\alpha_1 \vee \alpha_2$ is nice but not inert.
From Example~\ref{eg-wedge-of-spheres} $H_*(\lX;R) \isom UL_X$.
Thus as algebras 
\[ H_*(\lY;R) \isom U \left( L_X/[L^W_X] \amalg \freeL \langle [w]
\rangle \right).
\]
Therefore $\alpha_1 \vee \alpha_2$ is a semi-inert attaching map.
\eolBox

\begin{eg} \label{eg-fat-wedge}
The $6$-skeleton of $S^3 \times S^3 \times S^3$.
\end{eg}

This space $Y$ is also known as the \emph{fat
wedge}~\cite{selick:book} $FW(S^3, S^3, S^3)$. 
Let $X$ = $S^3_a \vee S^3_b \vee S^3_c$.
Let $\iota_a, \iota_b$ and $\iota_c$ be the inclusions of the
respective spheres.
Let $W = \bigvee_{j=1}^3 S^5$.
Then $Y = X \cup_f (\bigvee_{j=1}^3 e^6_j)$ where $f: W \to X$ is given by
$\bigvee_{j=1}^3 \alpha_j$ with $\alpha_1 = [\iota_b,\iota_c]$, $\alpha_2 =
[\iota_c,\iota_a]$ and $\alpha_3 = [\iota_a,\iota_b]$.
Let $\hat{\alpha}_j: S^2 \to \lX$ denote the adjoint of $\alpha_j$.
Let $I$ denote the two-sided ideal of $H_*(\lX;R)$ generated by
$\{h_X(\hat{\alpha}_j)\}_{j \in J}$.

Then $Y$ has Adams-Hilton model (see Section~\ref{section-ah}) $U(L,d)$
where $L = \freeL \langle x,y,z,a,b,c \rangle$,  $|x| = |y| = |z| = 2,  \
dx=dy=dz=0, \ da=[y,z], \ db=[z,x]$ and $dc=[x,y]$.
By Example~\ref{eg-alg-fat-wedge}, as algebras $HU(L,d) \isom
U(\freeL_{ab} \langle [x],[y],[z] \rangle \amalg \freeL \langle [w]
\rangle)$ where $w = [a,x] + [b,y] + [c,z]$.
Therefore, as algebras
\[ H_*(\lY;R) \isom \left( H_*(\lX;R) / I \right) \amalg \freeT \langle [w]
\rangle.
\]
Thus $f$ is nice but not inert.

From Example~\ref{eg-wedge-of-spheres} $H_*(\lX;R) \isom UL_X$.
Let $[L^W_X] \subset L_X$ be the Lie ideal generated by $L^W_X$.
Then as algebras
\[ H_*(\lY;R) \isom U \left( (L_X / [L^W_X]) \amalg \freeL \langle [w]
\rangle \right).
\]
Therefore $f$ is a semi-inert attaching map.
\eolBox

\section{Correspondence between topology and algebra}
\label{section-correspondence} 

Let $R = \Q$ or $\Fp$ where $p>3$ or $R \subset \Q$ be a subring
containing $\frac{1}{6}$.
Let $W \xto{f} X \to Y$ be the cofibration in the previous section
where $W = \bigvee_{j \in J} S^{n_j}$, $f = \bigvee_{j \in J} \alpha_j$ and
$H_*(\lX;R)$ is torsion-free.
Let $\hat{\alpha}_j$ denote the adjoint of $\alpha_j$.
Let $\AH{X}$ and $\AH{Y}$ be the corresponding AH-models given in
Section~\ref{section-ah}. 
Recall that $\AH{Y} = \AH{X} \amalg \freeT \langle y_j \rangle_{j \in
J}$.

Recall from Section~\ref{section-basic-top} that $L_X$ is the image of the
Hurewicz map.
Assume that $H_*(\lX;R) \isom UL_X$ as algebras.
Let $d':R\{y_j\} \to Z\AH{X} \to H\AH{X} \isomto UL_X$ be the induced map.
Recall from Section~\ref{section-ah} that $d\p y_j =
h_X(\hat{\alpha}_j) \in L_X$. 
That is, $\AH{Y}$ is a dga extension of $(\AH{X},L_X)$ (see
Section~\ref{section-inert-extns}).
Therefore the notation $\eL = (L_X \amalg \freeT \langle y_j
\rangle_{j \in J}, d\p)$ defined in Section~\ref{section-cap} is
consistent with the notation defined in
Section~\ref{section-inert-extns}.

In fact if $V_1 = R\{y_j\}_{j \in J}$ then $d\p V_1 \subset L^W_X$ and
the Lie subalgebra of $L_X$ generated by $d\p V_1$ is $L^W_X$.
Recall that $[V] \subset L_X$ denotes the Lie ideal generated by $V
\subset L_X$.
Therefore $[L^W_X] = [d\p V_1]$.
Similarly if $R \subset \Q$ then \feachpnP, $[\bar{L}^W_X] =
[\bar{d}\p \bar{V}_1]$. 

Recall from Section~\ref{section-cap} that if $R = \Q$ or $\Fp$ then
$f$ is free if $[L^W_X] \subset L_X$ is a free Lie 
algebra, and if $R \subset \Q$ then $f$ is free if \feachpnP,
$[\bar{L}^W_X] \subset \bar{L}_X$ is a free Lie algebra. 
Also recall from Section~\ref{section-inert-extns} that if $R = \Q$ or
$\Fp$ then $\AH{Y}$ is a free dga extension of $\AH{X}$ if $[d\p V_1]
\subset L_X$ is a free Lie algebra and if $R \subset \Q$ then $\AH{Y}$
is a free dga extension of $\AH{X}$ if \feachpnP, $[\bar{d}\p
\bar{V}_1] \subset \bar{L}_X$ is a free Lie algebra. 

Since $[L^W_X] = [d\p V_1]$ and $[\bar{L}^W_X] = [\bar{d}\p
\bar{V}_1]$, the conditions for $f$ being free and $\AH{Y}$ being a
free dga extension of $\AH{X}$ coincide. 
Likewise the inert and semi-inert conditions (see Sections
\ref{section-cap} and \ref{section-inert-extns}) coincide.

\section{Loop space decomposition theory} \label{section-hsb}
\label{section-lsd} 

Let $X$ denote a simply-connected CW complex of finite type.
Let $R \subset \Q$ be a subring of the rationals containing $1/2$ and
$1/3$.
Let $P$ be the set of primes invertible in $R$.
Let $\F = \Q$ or $\Fp$, where $p \notin P$.

An important idea in homotopy theory is to show that a given loop
space is homotopy equivalent to a product of other
spaces~\cite{cmn:aarhus, husemoller:cmns, cohen:decompositionsInHandbook}.
In particular, one can try and show that a space is homotopy
equivalent to a product of simpler \emph{atomic spaces} (defined below).
These simpler spaces will have the property that their loop spaces are
not homotopy equivalent to any nontrivial product. 
For example, if one can show that $\lY \approx \prod_i Y_i$ then
$\pi_n(Y) \isom \pi_{n-1}(\lY) \isom \prod_i \pi_{n-1} (Y_i)$.

A space $Y$ is called \emph{atomic} if it is $r$-connected, 
$\pi_{r+1}(Y)$ is a cyclic abelian group, and any self-map $f:Y\rightarrow Y$
inducing an isomorphism on $\pi_{r+1}(Y)$ is a homotopy equivalence.
Let $\mathcal{S} = \{ S^{2m-1},\,\Omega S^{2m+1} \:|\: m \geq 1\}$.
The spaces in $\mathcal{S}$ are atomic.
Let $\PiS$ be the collection of spaces homotopy equivalent
to a weak product of spaces in $\mathcal{S}$.
The properties of the homotopy groups of spaces in $\PiS$ are
determined by the properties of the homotopy groups of spheres.
For example if $\lX \in \PiS$ then $X$ `satisfies the statement of the
Moore conjecture'~\cite{selick:mooreConj}.
We make a brief digression explaining this statement.

Let $X$ be a finite simply-connected CW-complex.
The homotopy groups of $X$, $\{\pi_n(X)\}_{n \geq 2}$,
are finitely-generated abelian groups. 
As such, each has a torsion-free subgroup $\Z^{r_n}$, and for each prime
$p$, a finite subgroup $G_{n,p}$ of the elements with order $p^t$ for
some $t$.
Since each $p$-torsion subgroup $G_{n,p}$ is finite, there is a number
$m$ such that $p^m \cdot G_{n,p} = 0$. 
After tensoring $\pi_n(X)$ with the rational numbers $\Q$, only the
torsion-free subgroup remains. 

Call an integer $M$ an \emph{exponent} for a group $G$ if $M \cdot x=0 \
\forall x \in G$. 
Call $p^M$ a \emph{homotopy exponent for $X$ at the prime $p$} if
$\forall n$, $p^M \cdot G_{n,p} = 0$.
Call $p^M$ an \emph{eventual $H$-space exponent for $X$ at the prime
$p$} if $p^M$ is an exponent for $\Omega^N X$ for sufficiently large
$N$. 
The existence of an eventual $H$-space exponent implies the existence of a
homotopy exponent.

We say that $X$ is \emph{rationally elliptic} if $\pi_*(X) \tensor
\Q$ is finite.
Otherwise $X$ is said to be \emph{rationally hyperbolic}.
The Rational Dichotomy Theorem of F\'{e}lix, Halperin and Thomas
\cite{fht:homotopyLieAlg, felix:dichotomie, fht:rht} 
states that if $X$ is rationally hyperbolic then $\pi_*(X)\tensor \Q$
grows exponentially. 
Similarly $H_*(\lX;\Q)$ grows polynomially for
rationally elliptic spaces and exponentially for rationally hyperbolic
spaces. 

In the late 1970's, J.C. Moore conjectured a deep connection between
the rational homotopy groups $\pi_*(X) \tensor \Q$ and the $p$-torsion
subgroups for each prime $p$~\cite{selick:mooreConj}.

\begin{conj}[The Moore conjecture] \label{conj-wmc}
Let $X$ be a simply-connected finite CW-complex. \\
(a) if $X$ is rationally elliptic then $X$ has an eventual H-space
exponent at every prime $p$, and \\
(b) if $X$ is rationally hyperbolic then $X$ does not have a homotopy
exponent at any prime $p$. 
\end{conj}

This conjecture has been verified to hold for all but finitely many
primes for elliptic spaces~\cite{mcgibbonWilkerson} and
for 2-cones~\cite{anick:cat2}, but there have been only sparse results for
spaces outside these two classes~\cite{neisendorferSelick:someEgs,
selick:mooreSerre}.

\begin{lemma} \label{lemma-pis}
If $\lX \in \PiS$ then for all primes $p$,
$X$ is rationally elliptic if and only $X$ has a homotopy exponent at
$p$ if and only if $X$ has an eventual H-space exponent at $p$.
\end{lemma}

\begin{proof}
By assumption $\lX \approx \prod_{i \in I} T_i$ for some $T_i \in
\mathcal{S}$. 
Then $X$ is rationally elliptic if and only if $I$ is finite.
By \cite{cmn:annals109, cmn:annals110} and \cite{gray:sphereOfOrigin} this coincides with the
existence of a homotopy exponent and an eventual H-space exponent.
\end{proof}

From this it follows immediately:

\begin{cor}
Let $R \subset \Q$ be a subring containing $\frac{1}{6}$. 
Let $W \to X \to Y$ be a cell attachment satisfying the assumptions of
Theorem~\ref{thm-hurewiczR}. 
Then $\forall p \notin P_Y$ (see Section~\ref{section-ip}),
$Y$ is rationally elliptic iff $Y$ has a homotopy exponent at $p$ iff
$Y$ has an eventual $H$-space exponent at $p$.
\end{cor}

Let $R \subset \Q$ be a subring with invertible primes $P$ and let $X$
be a simply-connected topological space.
Then there exists a topological analogue to the algebraic localization
of a $\Z$-module at $R$.
We will denote the \emph{localization} at $R$ by $X_R$ which we will also
call the localization away from $P$.

There are many localization constructions. 
Perhaps the most widely used is that of Bousfield and
Kan~\cite{bousfieldKan:book} (outlined in \cite{selick:book}). 
We present a simple construction for the localization of
CW-complexes~\cite{fht:rht}.
One can construct $S^n_R$ to be the infinite mapping
telescope (see Section~\ref{section-basic-top}) of the sequence of
maps
\[ S^n \xto{j_1} S^n \xto{j_2} S^n \xto{j_3} \ldots S^n \xto{j_n}
\ldots
\]
where $j_k$ is the degree $m_k$ map where $m_k$ is the product of the
first $k$ primes in $P$.
Note that $e^{n+1} \isom (S^n \times I) / (S^n \times \{0\})$.
We can localize the $(n+1)$-cell by letting $e^{n+1}_R  =(S^n_R \times
I) / (S^n_R \times \{0\})$.
Since we can localize cells and spheres we can localize any
CW-complex.
One can check that this construction is functorial~\cite{fht:rht}.
Furthermore $\pi_*(X_R) = \pi_*(X) \tensor R$ and $H_*(X_R) = H_*(X)
\tensor R$.

The following is an expanded version of the Hilton-Serre-Baues
Theorem \cite{baues:commutator, anick:slsd}. 
Recall from Section~\ref{section-basic-top} that if $H_*(\lX; R)$ is
torsion-free,  we have maps $h_X$ and $h_X \tensor \F$.
Since $h_X$, $h_X \tensor \F$ are maps of Lie algebras taking
Samelson products to commutator brackets~\cite{samelson:products},
they induce maps $\tilde{h}_X: U\pi_*(\lX) \tensor R \to H_*(\lX; R)$ and
$\tilde{h}_X \tensor \F: U\pi_*(\lX) \tensor \F \to H_*(\lX; \F)$.

\begin{thm}[The Hilton-Serre-Baues Theorem] \label{thm-hsb}
The following are equivalent.

(a) There exists an $R$-equivalence $\tilde{\lambda}: \prod_i T_i \to
\lX$ where $T_i \in \mathcal{S}$ and the factors correspond\,\footnote{ 
\label{footnote-correspond} 
For $\lambda_i \in \pi_*(X) \tensor R$, let $\hat{\lambda}$
be the adjoint and $\{h(\hat{\lambda}_i)\}$ a $R$-basis for $L_X$.
If $n_i$ is even then $T_i = S^{n_i-1}$ and $\tilde{\lambda}_i =
\hat{\lambda}_i: T_i \to \lX$.
If $n_i$ is odd then $T_i = \Omega S^{n_i}$ and $\tilde{\lambda}_i =
\Omega\lambda_i: T_i \to \lX$.
}
to an
$R$-basis of the image of $h_X$. 
That is, localized at $R$, $\lX \in \PiS$.

(b) $H_*(\lX; R)$ is $R$-free and the map $\tilde{h}_X$ is
surjective.
That is, $H_*(\lX; R)$ is generated as an algebra by
the image of $h_X$.

(c) $H_*(\lX; R)$ is torsion-free and equal to a universal enveloping
algebra generated by $\im(h_X)$.

(d) $H_*(\lX; R)$ is torsion-free and $\forall p \notin P$ the map
$\tilde{h}_X \tensor \Fp$ is surjective.
That is, $H_*(\lX; \Fp)$ is generated as an algebra by the image of
$h_X \tensor \Fp$. 

(e) $H_*(\lX; R)$ is torsion-free and for all $p \notin P$, $H_*(\lX; \Fp)$ is
equal to a universal enveloping algebra generated by $\im(h_X \tensor \Fp)$.
\end{thm}

\begin{proof}
(a) $\Leftrightarrow$ (b) is the usual statement of the
theorem (\cite[Lemma V.3.10]{baues:commutator}, see also \cite[Lemma 3.1]{anick:cat2}). \\
(a) $\implies$ (c) 
By assumption, we have that $H_*(\lX; R) \isom H_*(\prod T_i;R)$. 
Use the Eilenberg-Zilber theorem and the K\"{u}nneth theorem
to write this as $\tensor_i H_*(T_i;R)$.
Let $m$ be an odd number. Then by the Bott-Samelson Theorem
(Theorem~\ref{thm-bott-sam}) $H_*(\Omega S^m;R) \isom U \freeL \langle
x_{m-1} \rangle$.  
In addition $H_*(S^m;R) \isom U \freeL_{ab} \langle x_m \rangle$, where
$\freeL_{ab}$ denotes the abelian Lie algebra.
Therefore $H_*(\lX;R) \isom \tensor_i UL_i \isom U(\prod_i L_i)$, where $L_i$ is generated by a Hurewicz image.
\\
(c) $\implies$ (b) is trivial. \\
Also (c) $\implies$ (e) $\implies$ (d) is trivial. \\
(d) $\implies$ (b) $\tilde{h}_X \tensor \Q$ is surjective by the
Milnor-Moore theorem. Since $\tilde{h}_X \tensor \Fp$ is surjective
for all $p \notin P$, it follows that $\tilde{h}_X$ is surjective.
\end{proof}

\section{Implicit primes} \label{section-ip}

Let $R \subset \Q$ be a subring with invertible
primes $P \supset \{2,3\}$.

Consider the cofibration $W \xrightarrow{f= \bigvee \alpha_j} X \to Y$, 
where $W = \bigvee_{j \in J} S^{n_j}$ is a finite-type wedge of spheres,
$H_*(\lX;R)$ is torsion-free and as algebras $H_*(\lX;R) \isom UL_X$.

We may need to exclude those primes $p$ for which an attaching map $\alpha_j
\in \pi_*(X)$ includes a term with $p$-torsion. 
We will define this set of \emph{implicit primes} below.

Recall from Section~\ref{section-basic-top} that $h_X: \pi_*(\lX)
\tensor R \to L_X \subset H_*(\lX;R)$ is a Lie algebra map.
Assume that there exists a Lie algebra map $\sigma_X: L_X \to
\pi_*(\lX)\tensor R$ such that $h_X \circ \sigma_X = id_{L_X}$.

By the Milnor-Moore theorem~\cite{milnorMoore:hopfAlgebras}, $h_X,
\sigma_X$ are rational isomorphisms, so $\im(\sigma_X \circ h_X - id)$
is a torsion element of $\pi_*(\lX) \tensor R$.
Let $\gamma_j = \sigma_X h_X(\hat{\alpha_j}) - \hat{\alpha}_j$ where
$\hat{\alpha}_j: S^{n_j-1} \to \lX$ is the adjoint of $\alpha_j$.
Then $t_j \gamma_j = 0 \in \pi_*(\lX) \tensor R$ for some $t_j>0$.
Let $t_j$ be the smallest such integer.

Define $P_Y$, the set of \emph{implicit primes} of $Y$ as follows.
A prime $p$ is in $P_Y$ if and only if $p \in P$ or 
$p | t_j$ for some $j \in J$.
We can invert these primes in $R$ to get a new ring $R\p =
\Z[{P_Y}^{-1}]$ with invertible primes $P_Y$.
Localized away from $P_Y$, $\hat{\alpha}_j = \sigma_X h_X
(\hat{\alpha}_j) = \sigma_X(dy_j)$. 

The implicit primes have the following properties.

\begin{lemma} \label{lemma-ip}
(a) Let $\{x_i\}$ be a set of Lie algebra generators for $L_X$ and let
$\beta_i = \sigma_X x_i$.
If all of the attaching maps are $R$-linear combinations of
iterated Whitehead products of the maps $\beta_i$, then $P_Y = P$. \\
(b) If $P = \{2,3\}$ and $n = \dim (Y)$ then the implicit primes are
bounded by $\max (3,n/2)$. 
\end{lemma}

\begin{proof}
(a) Since $\sigma_X$ is a Lie algebra map, in this case each
$\hat{\alpha}_j \in \im (\sigma_X)$.
So $\hat{\alpha}_j = \sigma_X a_j$ for some $a_j \in L_X$.
Then $\sigma_X h_X \hat{\alpha}_j = \sigma_X h_X \sigma_X a = \sigma_X
a = \hat{\alpha}_j$. 
Therefore $\forall j$, $\gamma_j = 0$ and $P_Y = P$. \\
(b) We assumed that $H_*(\lX;R) \isom UL_X$ as algebras.
By the Hilton-Serre-Baues Theorem (Theorem~\ref{thm-hsb}) $\lX \in \PiS$. 
So $\lX \approx \prod S^{2m_i-1} \times \prod \Omega S^{2m\p_k+1}$ where
$m_i, \ m\p_k \geq 1$.
Take $\hat{\alpha}_j$ to be the map $\hat{\alpha}_j: S^{n_j} \to \prod
S^{2m_i-1} \times \prod \Omega S^{2m\p_k+1}$.
Thus let $\hat{\alpha}_j = (g_i, g\p_k)$ where $g_i: S^{n_j} \to S^{2m_i-1}$
and $g\p_k: S^{n_j} \to \Omega S^{2m\p_k+1}$.
So $g_i \in \pi_{n_j}{S^{2m_i-1}}$ and $g\p_k \in \pi_{n_j+1}{S^{2m\p_k+1}}$.

Now in $\pi_*(S^n)$ the first $p$-torsion element occurs in
$\pi_{n+2p-3}(S^n)$.
Since $\pi_l(S^1)=0$ for $l>1$ so we can assume $n \geq 3$.
So there exists $p$-torsion in $\pi_l(S^n)$ only if $l \geq n+2p-3
\Leftrightarrow p \leq (l-n+3)/2 \leq l/2$ (since $n \geq 3$).
Hence $\hat{\alpha}_j$ contains $p$-torsion only if $p \leq (n_j+1)/2$.

Since $n_j \leq n-1$ the lemma follows.
\end{proof}

\section{Ganea's fiber-cofiber construction} \label{section-ganea}

We review Ganea's construction~\cite{ganea:construction}.
Let 
\begin{equation} \label{eqn-gan-fib} 
F \xto{\iota} X \xto{\rho} Y 
\end{equation}
be a fibration.
Let $X\p$ be the cofiber of $F \xto{\iota} X$.
Since~\eqref{eqn-gan-fib} is a fibration $\rho \circ \iota = *$.
Therefore $\rho$ extends to a map $\rho\p: X\p \to Y$.
Let $F\p$ be the homotopy fiber of $\rho\p$.

\begin{thm}[Ganea's Theorem~\cite{ganea:construction}]
$F\p$ is weakly homotopy equivalent to $F * \lY$ (where $A*B$ is the
\emph{join} of $A$ and $B$~\cite{selick:book}).
\end{thm}

Rutter~\cite{rutter:ganea} strengthened this theorem to the following.

\begin{thm}[\cite{rutter:ganea}] \label{thm-rutter}
$F\p \approx F * \lY$.
\end{thm}

\begin{rem}
If $F \to X \to Y$ is a fibration of CW-complexes then by Whitehead's
Theorem (Theorem~\ref{thm-whitehead}) these two theorems are equivalent.
\end{rem}

Mather~\cite{mather:ganea} generalized Ganea's construction as follows.
Given the fibration~\eqref{eqn-gan-fib} above and a map $W\xto{f} X$
such that $\rho \circ f \simeq *$, let $X\p$ be the cofiber of $f$.
Since $\rho \circ f \simeq *$, $\rho$ can be extended to a map $\rho\p:
X\p \to Y$.
Let $F\p$ be the homotopy fiber of $\rho'$.
This construction is also used in~\cite{felixThomas:fcof}

Since $\rho \circ f \simeq *$ there is a lifting $g: W \to F$.
Let $K$ be the cofiber of $g$.
\begin{thm}[\cite{mather:ganea}]
There is a cofibration $K \to F\p \to F * \lY$.
\end{thm}

Since Ganea's construction begin and ends with a fibration we can
iterate the construction to get the following commutative diagram of
fibrations.
\[
\xymatrix{
X \ar[d]_{\rho} \ar[r] & X_1 \ar[dl]_{\rho_1} \ar[r] & X_2
\ar[dll]_{\rho_2} \ar[r] & \ldots \ar[r] & X_n \ar[dllll]^{\rho_n}
\ar[r] & \ldots \\
Y
}
\]
Let $F_n$ be the homotopy fiber of $\rho_n$.

For example starting with the path-space fibration $\lX \to PX \to X$
the iterated Ganea construction together with Theorem~\ref{thm-rutter}
gives that $F_n \approx (\lX)^{*n}$ (the $n$-th fold join of $\lX$).

The iterated Ganea construction is particularly useful for studying
Lusternik-Schnirel\-mann category \cite{ganea:construction,
james:lscatInHandbook, fht:rht, selick:book}.

\part[Homology of DGA's and Loop Spaces]{Homology of Differential
Graded Algebras and Loop Spaces}

\chapter[The Homology of Differential Graded Algebras]{The Homology of
DGA's} \label{chapter-HUL} 

In this chapter we prove our main algebraic results: Theorems
\ref{thm-a} and \ref{thm-b}.


Let $R = \Q$ or $\Fp$ where $p>3$ or $R$ is a subring of $\Q$
containing $\frac{1}{6}$. 
Recall that if $R \subset \Q$ then $P$ is the set of invertible primes
in $R$ and $\nP = \{p \in \Z \ | \ p \text{ is prime and } p \notin
P\} \cup \{0\}$. 
Let $(\check{A},\check{d})$ be a connected finite-type dga over $R$ which is
$R$-free.
Let $Z\check{A}$ denote the subalgebra of cycles of $\check{A}$.
Let $\bA=(\check{A} \amalg \freeT V_1, d)$ be a connected finite-type dga
over $R$ where $V_1$ is a free $R$-module and $d\!\!\mid_{\check{A}} = \check{d}$
and $dV_1 \subset Z\check{A}$. 

Assume that as algebras $H(\check{A},\check{d}) \isom UL_0$ as
algebras for some Lie algebra $L_0$ which is a free $R$-module.
There is an induced map $d\p: V_1 \xto{d} Z\check{A} \to H\check{A}
\isomto UL_0$.  
Assume that $L_0$ can be chosen such that $d\p V_1 \subset L_0$.
In other words assume that $\bA$ is a dga extension of $((\check{A},
\check{d}), L_0)$ (see Section~\ref{section-inert-extns}).

Let $\eL = (L_0 \amalg \freeL V_1, d\p)$.
Then $\eL$ is a bigraded dgL where the Lie algebra $L_0$ and the free
$R$-module $V_1$ are in degrees $0$ and $1$ respectively and the
differential $d$ has bidegree $(-1,-1)$. 
Subscripts of bigraded objects will denote degree, eg. $M_0$ is the
component of $M$ in degree $0$. 

The following lemma is a well-known fact, and the subsequent two
lemmas are parts of lemmas from~\cite{anick:cat2}.
We remind the reader that all of our $R$-modules have finite type.

\begin{lemma} \label{lemma-isom}
Let $R \subset \Q$.
A homomorphism $f:M \to N$ is an isomorphism if and only if \feachpnP, $f
\tensor \Fp$ is an isomorphism.
\end{lemma}

Let $L$ be a connected bigraded dgL.
The inclusion $L \incl UL$ induces a natural map 
\begin{equation} \label{eqn-psi} 
\psi: UHL \to HUL. 
\end{equation}

\begin{lemma}[{\cite[Lemma 4.1]{anick:cat2}}] \label{lemma-psiF}
Let $R$ be a field of characteristic $k$ where $k=0$ or $k>3$.
Then the map $\psi$ in \eqref{eqn-psi} is an isomorphism in degrees
$0$ and $1$.
\end{lemma}

\begin{lemma}[{\cite[Lemma 4.3]{anick:cat2}}] \label{lemma-psiR}
Let $R$ be a subring of $\Q$ containing $\frac{1}{6}$.
Suppose that $HUL$ is $R$-free in degree $0$ and $1$.
Then $HL$ is $R$-free in degrees $0$ and $1$ and the map $\psi$ in
\eqref{eqn-psi} is an isomorphism in degrees $0$ and $1$.
\end{lemma}

Let $R = \F$ where $\F = \Q$ or $\Fp$ with $p>3$.
The Hilbert series of an $\F$-module is given by the power series
$A(z) = \Sigma_{n=0}^{\infty} (\Rank_{\F} A_n) z^n$.
Assuming that $A_0 \neq 0$, the notation $(A(z))^{-1}$ denotes the
power series $1/(A(z))$. 

Recall from Section~\ref{section-inert-extns} that $\bA$ is a filtered dga
under the increasing filtration given by 
$F_{-1}\bA = 0$, $F_0\bA = \check{A}$, and for $i\geq 0$, $F_{i+1}\bA =
\sum_{j=0}^i F_{j}\bA \cdot V_1 \cdot F_{i-j}\bA$.
We showed that this gives a first quadrant spectral sequence of algebras:
\[  
(E^0(\bA),d^0) = \gr(\bA) \implies E^{\infty} = \gr(H\bA) 
\]
where 
$E^0_{p,q}(\bA) = [F_p(\bA)/F_{p-1}(\bA)]_{p+q}$.

It is easy to check (see Section~\ref{section-ss}) that $(E^1,d^1)
\isom U\eL$ and hence $E^2 \isom HU\eL$. 
The following theorem follows from the main result of Anick's
thesis \cite[Theorem 3.7]{anick:thesis}.
Anick's theorem holds under either of two hypotheses. 
We will use only one of these.

\begin{thm} \label{thm-multF}
Let $R = \F$.
If the two-sided ideal $(d\p V_1) \subset UL_0$ is a free $UL_0$-module
then the above spectral sequence collapses at the $E^2$ term.
That is, $\gr(H\bA) \isom HU\eL$ as algebras.
Furthermore the multiplication map
\[ \nu: \freeT(\psi(H\eL)_1) \tensor (HU\eL)_0 \to HU\eL
\]
is an isomorphism and $(HU\eL)_0 \isom UL_0/(d\p V_1)$.
In addition,
\begin{equation} \label{eqn-a-formula} 
H\bA(z)^{-1} = HU\eL(z)^{-1} \\
= (1+z)(HU\eL)_0(z)^{-1} - z(UL_0)(z)^{-1} - V_1(z).
\end{equation}
\end{thm}

\begin{proof}
\cite[Theorem 3.7]{anick:thesis} shows that the spectral sequence collapses as
claimed and that the multiplication map $\freeT
W \tensor (HU\eL)_0 \to HU\eL$ is an isomorphism where $W$
is a basis for $(HU\eL)_1$ as a right $(HU\eL)_0$-module.
By Lemma~\ref{lemma-psiF} and the Poincar\'{e}-Birkhoff-Witt Theorem the
homomorphism $\psi(H\eL)_1 \tensor (HU\eL)_0 \to (HU\eL)_1$ induced by
multiplication in $HU\eL$ is an isomorphism.
So we can let $W = \psi(H\eL)_1$.

The remainder of the theorem follows directly from
\cite[Theorem 3.7]{anick:thesis}.
\end{proof}

\begin{cor} \label{cor-rfree}
If $R \subset \Q$ and \feachpnP, $(\bar{d}\bar{V}_1) \subset U\bar{L}_0$
is a free $U\bar{L}_0$-module then $H\bA$ is $R$-free iff $HU\eL$ is
$R$-free iff $L_0/[d\p V_1]$ is $R$-free. 
\end{cor}

\begin{proof}
First $(HU\eL)_0 \isom UL_0/(d\p V_1) \isom U(L_0/[d\p V_1])$.
So $(HU\eL)_0$ is $R$-free if and only if $L_0/[d\p V_1]$ is $R$-free.

Let $A_p(z)$, $B_p(z)$ and $C_p(z)$ be the left, middle and right parts of
\eqref{eqn-a-formula} for $\F = \F_p$ with $\F_0 = \Q$.
Then $H\bA$ is $R$-free iff $\fpnP$ $A_p(z) = A_0(z)$ and 
$HU\eL$ is $R$-free iff $\fpnP$ $B_p(z) = B_0(z)$. 
Since $UL_0$ and $V_1$ are $R$-free, $(HU\eL)_0$ is $R$-free iff
$\fpnP$ $C_p(z) = C_0(z)$.
Since $\fpnP$ $A_p(z) = B_p(z) = C_p(z)$ this proves the corollary.
\end{proof}

We now prove a version of Theorem~\ref{thm-multF} for subrings of $\Q$.

\begin{thm} \label{thm-multR}
Let $R \subset \Q$.
If $L_0/[d\p V_1]$ is $R$-free and \feachpnP, $(\bar{d}\bar{V}_1) \subset
U\bar{L}_0$ is a free $U\bar{L}_0$-module, then $H\bA$ is $R$-free
and the multiplication map
\[ \nu: \freeT(\psi(H\eL)_1) \tensor (HU\eL)_0 \to HU\eL
\]
is an isomorphism.
Also $\gr(H\bA) \isom HU\eL$ as algebras and
$(HU\eL)_0 \isom UL_0/(d\p V_1)$. 
\end{thm}

\begin{proof}
Since $L_0/[d\p V_1]$ is $R$-free, by Corollary~\ref{cor-rfree} so are
$HU\eL$ and $H\bA$.
It follows from the Universal Coefficient Theorem that $\fpnP$, $H\bA \tensor
\Fp \isom H(\bA \tensor \Fp)$ and $HU\eL \tensor \Fp \isom HU(\eL
\tensor \Fp)$. 
Hence $\fpnP$, $(HU\eL)_0 \tensor \Fp \isom (HU(\eL\tensor \Fp))_0$. 
Using Lemmas \ref{lemma-psiR} and \ref{lemma-psiF} 
\[ \psi(H\eL)_1 \tensor \F_p \isom (H\eL)_1 \tensor \Fp \isom
H(\eL\tensor \Fp)_1 \isom \psi H(\eL\tensor\Fp)_1.
\]
Thus $\fpnP$ 
\[ \nu \tensor \Fp: \freeT(\psi(H\eL)_1 \tensor \Fp) \tensor (HU\eL)_0
\tensor \Fp \to H\bA \tensor \Fp
\]
corresponds under these isomorphisms to the multiplication map
\[ \freeT(\psi(H(\eL\tensor \Fp))_1) \tensor (HU(\eL \tensor \Fp))_0 \to
H(\bA \tensor \Fp). 
\]
But this is an isomorphism by Theorem~\ref{thm-multF}.
Therefore $\nu$ is an isomorphism by Lemma~\ref{lemma-isom}.

The last two isomorphisms also follow from Theorem~\ref{thm-multF}.
\end{proof}

The next lemma will prove that if the Lie ideal $[d\p V_1] \subset
L_0$ is a free Lie algebra then the hypothesis in Anick's Theorem
(Theorem~\ref{thm-multF}) holds. That is, $(d\p V_1)$ is a free $UL_0$-module.

\begin{lemma} \label{lemma-I}
Let $A=UL$ over a field \F. 
Assume $I$ is a two-sided ideal of $A$ generated by a Lie ideal $J$ 
of $L$ which is a free Lie algebra, $\freeL W$.
Then the multiplication maps $A \tensor W \to I$ and $W \tensor A \to
I$ are isomorphisms of left and right $A$-modules respectively.
\end{lemma}

\begin{proof}
From the short exact sequence of Lie algebras 
\[ 0 \rightarrow J \rightarrow L \rightarrow L/J \rightarrow 0 \]
we get the short exact of sequence of Hopf algebras \cite[Theorem
10.5.3]{selick:book}  
\[ \F \to U(J) \rightarrow U(L) \rightarrow U(L/J) \to \F \]
and so $UL \isom UJ \tensor U(L/J)$ as $\F$-modules. 
Since $J$ is a free Lie algebra $\mathbb{L}W$, $UJ \cong \freeT W$. It
is also a basic fact that $U(L/J) \isom UL/I$.
Hence we have that 
\begin{equation} \label{eqn-Avs1} A \isom TW \tensor A/I 
\end{equation}
as $\F$-modules.
Furthermore
\begin{equation} \label{eqn-Avs2} A \isom I \oplus A/I 
\end{equation}
as $\F$-modules.

Let $M(z)$ denote the Hilbert series for the $\F$-module $M$, and to
simplify the notation let $B = A/I$.
Then from equations \eqref{eqn-Avs1} and \eqref{eqn-Avs2} we have the
following (using $(TW)(z) = 1/(1-W(z))$).
\[ B(z) = A(z) (1 - W(z)), \quad
I(z) = A(z) - B(z).
\]
Combining these we have $I(z) = A(z) W(z)$. That is, $I \isom A \tensor W$
as $\F$-modules.

Let $\mu: A \tensor W \to I$ be the multiplication map. To show that
it is an isomorphism it remains to show that it is either injective
or surjective.

We claim that $\mu$ is surjective.
Since $I$ is the ideal in $A$ generated by $W$, any $x \in I$ can be
written as 
\begin{equation} \label{eqn-awb}
x= \Sigma_i a_i w_i b_{i_1} \cdots b_{i_{m_i}} \text{, where }
a_i \in A, \ w_i \in W \text{ and } b_{i_k} \in L.
\end{equation}
Each such expression gives a sequence of numbers $\{m_i\}$.
Let $M(x) = \min \max_i(m_i)$, where the minimum is taken over all
possible ways of writing $x$ as in~\eqref{eqn-awb}.
We claim that $M(x) = 0$.

Assume that $M(x) = t>0$.
Then $x = x\p + \Sigma_i a_i w_i b_{i_1} \cdots b_{i_t}$, where $M(x\p)
<t$. 
Now $w_i b_{i_1} = [w_i, b_{i_1}] \pm b_{i_1} w_i$. 
Furthermore since $J$ is a Lie ideal $[w_i, b_{i_1}] \in J \isom 
\freeL W$, so  
\[ [w_i, b_{i_1}] = \Sigma_j c_j [[w_{j_1}, \ldots , w_{j_{n_j}}] 
= \Sigma_k d_k w_{k_1} \cdots w_{k_{N_k}} 
= \Sigma_l a_l w_l,
\]
where $a_l \in A$ and $w_l \in W$.
So $x = x\p + \Sigma_i \Sigma_{l_i} a_i a_{l_i} w_{l_i} b_{i_2} \cdots
b_{i_t}$.
But this is of the form in~\eqref{eqn-awb} and shows that $M(x) \leq
t-1$ which is a contradiction.

Therefore for $x \in I$ $M(x)=0$ and we can write $x = \Sigma_i a_i
w_i$ where $a_i \in A$ and $w_i \in W$.
Then $x \in \im(\mu)$ and hence $\mu$ is an isomorphism.

Since $A$ is associative $\mu$ is a map of left $A$-modules.

The second isomorphism follows similarly.
\end{proof}

We are now almost ready to prove Theorems \ref{thm-a} and \ref{thm-b}.
Recall that $\bA = (\check{A} \amalg \freeT V_1, d)$ where $d\check{A} \subset
\check{A}$ and $dV_1 \subset Z\check{A}$.
Also $H(\check{A},\check{d}) \isom UL_0$ as algebras and if
$d':V_1 \to UL_0$ is the induced map then $d'(V_1) \subset L_0$.
Let $\eL = (L_0 \amalg \freeL V_1,d\p)$ with $d\p L_0 = 0$.
Recall that there is a map $\underline{\psi}: UH\eL \to HU\eL$.
$\bA$ is a filtered dga under the increasing filtration given by
$F_{-1}\bA = 0$, $F_0\bA = \check{A}$, and for $i\geq 0$, $F_{i+1} =
\sum_{k=0}^i F_k\bA \cdot V_1 \cdot F_{i-k}\bA$.
We prove one last lemma.

\begin{lemma} \label{lemma-quotientf}
There exists a quotient map
\begin{equation} \label{eqn-quotientf} 
f: F_1 H\bA \to (HU\eL)_1.
\end{equation}
Given $R\{\bar{w}_i\} \subset (H\eL)_1$ there exist cycles
$\{w_i\} \subset F_1 \bA$ such that $f([w_i])=\bar{w}_i$.
\end{lemma}

\begin{proof}
By Theorem~\ref{thm-multF} $(\gr(H\bA))_1 \isom (HU\eL)_1$. 
So there is a quotient map 
\[ f: F_1 H\bA \to (\gr(H\bA))_1 \isom (HU\eL)_1.
\]
By Lemma~\ref{lemma-psiF} $(H\eL)_1 \isom (\underline{\psi}H\eL)_1
\subset (HU\eL)_1$.
So for each $\bar{w}_i$ one can choose a representative cycle $w_i \in
Z F_1 \bA$ such that $f([w_i]) = \bar{w}_i$.
\end{proof}

We are now ready to prove Theorem~\ref{thm-a}. 
We prove a slightly more detailed form of the theorem which we now state.

\begin{thm}[Theorem~\ref{thm-a}] 
Let $R = \F$ where $\F= \Q$ or $\Fp$ with $p>3$.
Let $(\check{A}, \check{d})$ be a connected finite-type dga and let
$V_1$ be a connected finite-type $\F$-module with a map $d:V_1 \to
\check{A}$. 
Let $\bA = (\check{A} \amalg \freeT V_1, d)$.
Assume that there exists a Lie algebra $L_0$ such that
$H(\check{A},\check{d}) \isom UL_0$ as algebras and $d\p V_1 \subset
L_0$ where $d\p$ is the induced map.
Also assume that $[d\p V_1] \subset L_0$ is a free Lie algebra.
That is, $\bA$ is a \emph{dga extension} of
$((\check{A},\check{d}),L_0)$ which is \emph{free}. 
Let $\eL = (L_0 \amalg \freeL V_1, d\p)$. \\ 
(a) Then as algebras
\[ \gr(H\bA) \isom U ((H\eL)_0 \sdp \freeL (H\eL)_1)
\]
with $(H\eL)_0 \isom L_0/[d\p V_1]$ as Lie algebras.  \\
(b) Furthermore if $\bA$ is \emph{semi-inert} (that is, there is a
free $R$-module $K$ such that $(H\eL)_0 \sdp \freeL (H\eL)_1 \isom
(H\eL)_0 \amalg \freeL K$) then as algebras
\[ H\bA \isom U ( (H\eL)_0 \amalg  \freeL K\p )
\]
for some $K\p \subset F_1 H\bA$ such that $f: K\p \isomto K$, where
$f$ is the quotient map in Lemma~\ref{lemma-quotientf}.
\end{thm}

\begin{proof}
Let $R=\F$ where $\F=\Q$ or $\Fp$ with $p>3$.
Assume $[d\p V_1] \subset L_0$ is a free Lie algebra.
  
\noindent 
(a) By Lemma~\ref{lemma-I} $(d\p V_1) \subset UL_0$ is a free $UL_0$-module.
So we can apply Theorem~\ref{thm-multF} to show that
$\gr(H\bA) \isom HU\eL$ as algebras and that the
multiplication map
\[ \nu: \freeT(\underline{\psi}(H\eL)_1) \tensor (HU\eL)_0 \to HU\eL
\]
is an isomorphism.

By Lemma~\ref{lemma-psiF} $(HU\eL)_0 \isom U(H\eL)_0$ and
$\underline{\psi}(H\eL)_1 \isom (H\eL)_1$.
By the definition of homology $(H\eL)_0 \isom L_0/[d\p V_1]$.

$(H\eL)_0$ acts on $(H\eL)_1$ via the adjoint action.
Let $L\p = (H\eL)_0 \sdp \freeL (\underline{\psi}(H\eL)_1)$.
There is an induced Lie algebra map $u:L\p \to HU\eL$ which gives an
induced algebra map $\tilde{u}: UL\p \to HU\eL$.
To simplify the notation we will refer to $(H\eL)_0$ and
$\underline{\psi} (H\eL)_1$ by $L\p_0$ and $L\p_1$ respectively.

Recall that as $R$-modules, $L\p \isom L\p_0 \times \freeL L\p_1$.
The Poincar\'{e}-Birkhoff-Witt Theorem shows that the multiplication
map
\[ \phi: \freeT L\p_1 \tensor (HU\eL)_0 \isom U\freeL L\p_1 \tensor U L\p_0
\to UL\p
\]
is an isomorphism.
Since $\tilde{u}$ is an algebra map $\nu = \tilde{u} \circ \phi$.
Therefore $\tilde{u}: UL\p \to HU\eL$ is an isomorphism.
Hence  $\gr(H\bA) \isom  UL\p$ as algebras.
This finishes the first part of the Theorem.

\noindent 
(b) Recall that $L\p_0$ acts on $L\p_1  = (\underline{\psi}H\eL)_1
\isom (H\eL)_1$ via the adjoint action.
Assume that $\bA$ is semi-inert. That is, $\exists \{\bar{w}_i\} \subset L\p_1$
such that $L\p \isom L\p_0 \amalg \freeL K$, where $K = R\{\bar{w}_i\}$. 
Recall from (a) that $HU\eL \isom \gr(H\bA)$ and that the inclusions
$L\p_0 \subset HU\eL$ and $\bar{w}_i \in HU\eL$ induce a Lie algebra map 
\[ u: L\p \to \gr(H\bA).
\]

By Lemma~\ref{lemma-quotientf} $\exists \ w_i \in F_1 \bA$ such that
$f([w_i]) = \bar{w}_i$ where $f$ is the map in~\eqref{eqn-quotientf}.
Let $K\p = R\{[w_i]\} \subset F_1 H\bA$, and let $L\pp = L\p_0 \amalg
\freeL K\p$. 
Then $f: K\p \xrightarrow{\isom} K$ and $f$ induces an isomorphism
$L\pp \xrightarrow{\isom} L\p$.

By part (a) $L\p_0 \subset (\gr H\bA)_0$. 
Since $F_{-1} H\bA = 0$, $(\gr H\bA)_0 \subset H\bA$, so $L\p_0
\subset H\bA$.
Since $L\p_0$, $K\p \subset H\bA$ there are induced maps 
\[
\xymatrix{
L\pp \ar[r]^-{\eta} \ar[d] & H\bA \\
UL\pp \ar[ur]_{\theta}
},
\]
where $\theta$ is an algebra map.

Grade $L\pp$ by letting $L\p_0$ be in degree $0$ and $K\p$ be in
degree $1$.
This also filters $L\pp$.
Then $\eta$ is a map of filtered objects.

From this we get the following commutative diagram
\[
\xymatrix{
\gr(L\pp) \ar[d] \ar[r]^-{\gr(\eta)} & \gr(H\bA) \\
\gr(UL\pp) \ar[ur]^{\gr(\theta)} \ar[d]_{\isom} \\
U\gr(L\pp) \ar[uur]_{\rho}
}.
\]

Now $\gr(L\pp) \isom L\pp \isom L\p$ and $\gr(\eta)$ corresponds to $u$ under
this isomorphism. 
So $\rho$ corresponds to $\tilde{u}$ which is an isomorphism.
Thus $\gr(\theta)$ is an isomorphism, and hence $\theta$ is an
isomorphism. 
Therefore $H\bA \isom UL\pp$ which finishes the proof.
\end{proof}


We will prove a slightly more detailed form of Theorem~\ref{thm-b}
which we now state.

\begin{thm}[Theorem~\ref{thm-b}] \label{thm-b'}
Let $R \subset \Q$ be a subring containing $\frac{1}{6}$.
Let $(\check{A}, \check{d})$ be a connected finite-type dga and let
$V_1$ be a connected finite-type free $R$-module with a map $d:V_1 \to
\check{A}$. 
Let $\bA = (\check{A} \amalg \freeT V_1, d)$.
Assume that $H(\check{A},\check{d})$ is $R$-free and that there exists
a Lie algebra $L_0$ such that $H(\check{A},\check{d}) \isom UL_0$ as
algebras and $d\p V_1 \subset L_0$ where $d\p$ is the induced map.
Also assume that $L_0/[d\p V_1]$ is a free $R$-module and that
\fanypnP, $[d\p V_1] \subset L_0$ is a free Lie algebra. 
That is, $\bA$ is a \emph{dga extension} of
$((\check{A},\check{d}),L_0)$ which is \emph{free}. 
Let $\eL = (L_0 \amalg \freeL V_1, d\p)$. \\ 
(a) Then $H\bA$ and $\gr(H\bA)$ are $R$-free and as algebras
\[ \gr(H\bA) \isom U ( (H\eL)_0 \sdp \freeL (H\eL)_1 )
\]
with $(H\eL)_0 \isom L_0/[d\p V_1]$ as Lie algebras.
Additionally $(H\eL)_0 \sdp \freeL (H\eL)_1 \isom \underline{\psi}
H\eL$ as Lie algebras.\\
(b) Furthermore if $\bA$ is \emph{semi-inert} (that is, there is a
free $R$-module $K$ such that $(H\eL)_0 \sdp \freeL (H\eL)_1 \isom
(H\eL)_0 \amalg \freeL K$) then as algebras
\[ H\bA \isom U ( (H\eL)_0 \amalg  \freeL K\p )
\]
for some $K\p \subset F_1 H\bA$ such that $f: K\p \isomto K$, where
$f$ is the quotient map in Lemma~\ref{lemma-quotientf}.
\end{thm}

\begin{proof}
We follow the argument in the proof of Theorem~\ref{thm-a} to prove
Theorem~\ref{thm-b}, but add an additional argument to show that
$\underline{\psi} H\eL \isom (H\eL)_0 \sdp \freeL (H\eL)_1$ as Lie algebras.

Let $R$ be a subring of $\Q$ containing $\frac{1}{6}$.
Assume $L_0/[d\p V_1]$ is $R$-free and that \feachpnP,
$[\bar{d}\bar{V}_1] \subset \bar{L}_0$ is a free Lie algebra.

\noindent
(a) By Corollary~\ref{cor-rfree}, $H\bA$ and $HU\eL$ are $R$-free and
by Lemma~\ref{lemma-I}, \feachpnP, $(\bar{d}\bar{V}_1) \subset U\bar{L}_0$ is a
free $U\bar{L}_0$-module.
So we can apply Theorem~\ref{thm-multR} to show that $\gr(H\bA)
\isom HU\eL$ as algebras and that the
multiplication map
\[ \nu: \freeT(\psi(H\eL)_1) \tensor (HU\eL)_0 \to HU\eL
\]
is an isomorphism.

By Lemma~\ref{lemma-psiR} $(HU\eL)_0 \isom U(H\eL)_0$ and by the
definition of homology $(H\eL)_0 \isom L_0/[d\p V_1]$.

Let $N = \underline{\psi}(H\eL)$ and
define
\[
L\p = N_0 \sdp \freeL N_1.
\]
Note that $L\p_0 = N_0$ and $L\p_1 = N_1$.
There is a Lie algebra map $u: L\p \to N$ and
an induced algebra map $\tilde{u}: UL\p \xrightarrow{Uu}
UN \to HU\eL$. 

Again the Poincar\'{e}-Birkhoff-Witt Theorem shows that the
multiplication map
\[ \phi: \freeT N_1 \tensor (HU\eL)_0 \isom U\freeL N_1 \tensor U N_0 \to
UL\p
\]
is an isomorphism.
Since $\tilde{u}$ is an algebra map $\nu = \tilde{u} \circ \phi$.
Thus $\tilde{u}$ is an isomorphism.
Therefore $HU\eL \isom UL\p$ as algebras, as in the proof of
Theorem~\ref{thm-a}. 

Unlike the proof of Theorem~\ref{thm-a}, we will show that $u: L\p \to
N$ is an isomorphism. 
Let $\iota:N \incl HU\eL$ be the inclusion.
Since the composition $L\p \xrightarrow{u} N
\xrightarrow{\iota} HU\eL \isom UL\p$ is 
the inclusion $L\p \incl UL\p$, $u$ is injective.
So as $R$-modules $N \isom L\p \oplus N/L\p$.
Since $L\p$ and $N$ are $R$-free, so is $N/L\p$.

Recall that the composition $\iota \circ u$ induces the isomorphism
$\tilde{u}: UL\p \xrightarrow{Uu} UN \to HU\eL$.
Tensor these maps with $\Q$ to get the commutative diagram
\begin{equation} \label{eqn-cdq}
\xymatrix{
UL\p\tensor\Q \ar[r]^{Uu\tensor\Q} \ar[dr]_{\isom} & UN \tensor \Q
\ar[d] \\
& HU\eL \tensor \Q }.
\end{equation}
It is a classical result that the natural map
\begin{equation} \label{eqn-psiq} 
\underline{\psi}_{\Q}: UH(\eL \tensor \Q) \xrightarrow{\isom} HU(\eL
\tensor \Q)
\end{equation} 
is an isomorphism (see \cite[Theorem 21.7(i)]{fht:rht} for example). 
Notice that
\[ N \tensor \Q = (\underline{\psi}H\eL) \tensor \Q
\isom \underline{\psi}_{\Q} H(\eL \tensor \Q) \isom H(\eL \tensor \Q)
\]
and $HU\eL \tensor \Q \isom HU(\eL \tensor \Q)$.
Under these isomorphisms the vertical map in~\eqref{eqn-cdq}
corresponds to the isomorphism in~\eqref{eqn-psiq}.

Therefore $Uu \tensor \Q$ is an isomorphism and hence $u\tensor\Q$ is
surjective.
As a result $\coker u = N/L\p$ is a torsion $R$-module.
But we have already shown that $N/L\p$ is $R$-free.
Thus $N/L\p = 0$ and $N \isom
L\p$. 
Hence $HU\eL \isom UN$.

If $\bA$ is semi-inert then the proof of (b) is the
same as for Theorem~\ref{thm-a}(b).
\end{proof}

As claimed in the introduction, we prove that the following corollary
follows from Theorems \ref{thm-a} and \ref{thm-b}.

\begin{cor}[Corollary~\ref{cor-semiInertEqnAlg}]
Let $\bA$ be a semi-inert dga extension satisfying the hypotheses of
either Theorems \ref{thm-a} or \ref{thm-b}.
Then
\begin{equation*} 
K\p (z) = V_1(z) + z[UL_0(z)^{-1} - U(H\eL)_0(z)^{-1}].
\end{equation*}
\end{cor}

\begin{proof}
By Theorem~\ref{thm-a} or Theorem~\ref{thm-b},
\[ H\bA \isom U((H\eL)_0 \amalg \freeL K\p) \isom U(H\eL)_0 \amalg
\freeT K\p.
\]
Since $(A \amalg B)(z)^{-1} = A(z)^{-1} + B(z)^{-1} - 1$~\cite[Lemma
5.1.10]{lemaire:monograph} it follows that 
\[ H\bA(z)^{-1} = U(H\eL)_0(z)^{-1} - K'(z)
\]
Therefore, using Anick's formula~\ref{eqn-a-formula} we have that
\begin{eqnarray*}
K'(z) & = & U(H\eL)_0(z)^{-1} - H\bA(z)^{-1} \\
& = & U(H\eL)_0(z)^{-1} - (1+z)(HU\eL)_0(z)^{-1} + z UL_0(z)^{-1} +
V_1(z) \\
& = & V_1(z) + z[U(H\eL)_0(z)^{-1} - UL_0(z)^{-1}]
\end{eqnarray*}
\end{proof}

\ignore
{

\section{The algebraic structure}

In this section we give an example which illustrates why part (i) of
Theorems \ref{thm-a} and \ref{thm-b} gives the algebra structure only
for the associated graded object $H\bA$ and not $H\bA$ itself. 

We will consider the case when $R = \Q$ and $\bA = U\bL$.
In this case, the Milnor-Moore Theorem~\cite{milnorMoore:hopfAlgebras}
tells us that $HU\bL \isom U\psi H\bL \isom UH\bL$.
From Theorems \ref{thm-a} and \ref{thm-b} we have that $\gr(H\bL)
\isom (H\eL)_0 \sdp \freeL (H\eL)_1$ as Lie algebras, where $(H\eL)_0
\isom L_0/[dV_1]$ [***]. 
Equivalently, by Lemma~\ref{lemma-sdp} we have the split short exact
sequence of Lie algebras
\[ 0 \to \freeL (H\eL)_1 \to \gr(H\bL) \to L_0/[dV_1] \to 0.
\]

We will show that in general it is not true that as Lie algebras $H\bL
\isom \freeL K \sdp L_0/[dV_1]$ for some $\Q$-module $K$ (this was
reported to be the case in \cite[Theorem 3]{felixThomas:fcof}). 
In the following example we will see that if $M$ is an $R$-module
complement of $L_0/[dV_1]$ then $M$ cannot be a Lie ideal.

\begin{eg}
$(\check{L}, \check{d}) = (\freeL \langle x_1,x_2,y_1,y_2,u_1,u_2,v_1,v_2 \rangle,
\check{d})$ where for $i=1,2$, $|x_i| = |y_i| = 2$, $\check{d} x_i = \check{d} y_i = 0$,
$\check{d} u_i = [[x_1,x_2],y_i]$ and $\check{d} v_i = [[y_1,y_2],x_i]$.
$\bL = (\check{L} \amalg \freeL \langle a,b \rangle, d)$ where $da = x_1$
and $db = [y_1,y_2]$. [xxx $[d\pV_1]$ is not a free Lie algebra]
\end{eg}

Let $\{ F_i \bL \}$ be the usual filtration.
Let $L_0 = H(\check{L}, \check{d})$. 
Abusing notation we will refer to the homology classes represented by
$x_i$ and $y_i$ by $x_i$ and $y_i$.
Let $\eL = (L_0 \amalg \freeL \langle a,b
\rangle, d\p)$ where $d\p|_{\check{L}} = 0$ and $d\p a = x_1$ and $d\p b =
[y_1, y_2]$.

One can check that for $i=1,2$, $[a,x_2,y_i]$ and  $[b,x_i]$ are cycles
in $\eL$.
Let $\bar{\alpha}_i$ and $\bar{\beta}_i$ be the corresponding homology
classes in $(H\eL)_1$.

Now 
\begin{eqnarray*} 
d\p [[a,b],x_2] & = & [[x_1,b],x_2] - [[a,[y_1,y_2]],x_2] \\
& = & [[x_1,b],x_2] - [[a,x_2],[y_1,y_2]] \\
& = & - [[b,x_1],x_2] - [[[a,x_2],y_1],y_2] + [[[a,x_2],y_2],y_1].
\end{eqnarray*}
Therefore 
\begin{equation} \label{eqn-barGamma}
-[\bar{\beta}_1,x_2] - [\bar{\alpha}_1,y_2] + [\bar{\alpha}_2,y_1] =
 0 \text{ in } (H\eL)_1.
\end{equation}

We will show that the corresponding sum in $H\bL$ is not zero. 
For $i=1,2$, $\bar{\alpha}_i$ is represented by the cycle $\alpha_i =
[[a,x_2],y_i] - u_i$ in $\bL$ and $\bar{\beta}_i$ is represented by the
cycle $\beta_i = [b,x_i] - v_i$.

Let $\gamma = -[\beta_1,x_2] - [\alpha_1,y_2] + [\alpha_2,y_1] \in
Z F_1 \bL$ which is the sum in $F_1 \bL$ corresponding to
\eqref{eqn-barGamma}.
Let $\epsilon = [v_1,x_2] - [v_2,x_1] + [u_1,y_2] - [u_2,y_1] \in Z
F_0 \bL$.
The one can check $\gamma$ is not a boundary but that 
\[ d([[a,b],x_2] - [a,v_2]) = \gamma - \epsilon.
\]
As a result $[\gamma] = [\epsilon] \in F_0 H\bL = L_0/[dV_1]$.
Therefore we cannot write $H\bL \isom M \times L_0/[dV_1]$ as
$R$-modules, where $M$ is a Lie ideal.

This difficulty arises from the fact that $\bL$ is not a bigraded dgL.

}

\chapter{Application to Cell Attachments} \label{chapter-attach}

In this chapter we prove two of our main topological results: Theorems
\ref{thm-c} and \ref{thm-d}.

Let $f: W \stackrel{f}{\rightarrow} X$ be a map of finite-type
simply-connected CW-complexes where $W = \bigvee_{j \in J} S^{m_j}$
and $f = \bigvee_{j \in J} \alpha_j$.
Let $\hat{\alpha}_j: S^{m_j-1} \to \lX$ denote the adjoint of $\alpha_j$.
Let $Y$ be the adjunction space 
\[ Y = W \cup_f \left( \bigvee_{j \in J} e^{m_j+1} \right).
\]

Let $R = \Q$ or $\Fp$ where $p >3$ or $R$ is a subring of $\Q$
containing $\frac{1}{6}$.
Recall that $P$ is the set of primes invertible in $R$, $\nP =
\{p \in \Z \ | \ p \text{ is prime and }  p \notin P\} \cup \{0\}$ and
$\F_0 = \Q$. 

Recall from Section~\ref{section-ah} that we can choose Adams-Hilton
models $\AH{X}$ and $\AH{Y}$ for $X$ and $Y$ such that
\begin{equation} \label{eqn-lly} 
\AH{Y} = \AH{X} \amalg \freeT \langle y_j \rangle_{j \in J}. 
\end{equation}
These come with algebra isomorphisms $H_*(\lX;R) \isom H\AH{X}$
and $H_*(\lY;R) \isom H\AH{Y}$. 

Filter $\AH{Y}$ by taking $F_{-1}\AH{Y} = 0$, $F_0\AH{Y} = \AH{X}$ and
for $i\geq 0$, $F_i\AH{Y} = \sum_{k=0}^i F_{k}\AH{Y} \cdot R\{y_j\}_{j
\in J} \cdot F_{i-k}\AH{Y}$.
This filtration makes $\AH{Y}$ a filtered dga.

Recall from Section~\ref{section-ss} that this filtration induces a
first quadrant multiplicative spectral sequence converging from
$\gr(\AH{Y})$ to $\gr(H\AH{Y})$.

Recall that $L_X = \im(h_X: \pi_*(\lX)\tensor R \to H_*(\lX;R))$.
Assume that $H_*(\lX;R) \isom UL_X$ as algebras and
that it is $R$-free. 
Then $(E^1,d^1) \isom (UL_X \amalg U \freeL \langle y_j
\rangle_{j\in J}, d\p)$ 
where $d\p$ is determined by the induced map $d':R\{y_j\} \to Z\AH{X}
\isomto H_*(\lX;R) \isomto UL_X$.
Recall from Section~\ref{section-ah} that $d\p y_j =
h_X(\hat{\alpha}_j) \in L_X$. 
Therefore $\AH{Y}$ is a dga extension (see
Section~\ref{section-inert-extns}) of $(\AH{X},L_X)$. 

Furthermore, $L^W_X$ is the Lie subalgebra of $L_X$ generated by
$R\{d\p y_j\}_{j \in J}$.
Therefore $[L^W_X] = [R\{d\p y_j\}]$.
Let $\eL = (L_X \amalg \freeL \langle y_j \rangle, d\p)$.
Then $(E^1, d^1) \isom U\eL$.

We can prove Theorem~\ref{thm-c} by applying Theorem \ref{thm-a} to
$\AH{Y}$ and $\eL$. 
For convenience we restate Theorem~\ref{thm-c} here.

\begin{thm}[Theorem \ref{thm-c}]
Let $R = \F$ where $\F = \Q$ or $\Fp$ with $p > 3$.
Let $X$ be a finite-type simply-connected CW-complex such that
$H_*(\lX;R)$ is torsion-free and as algebras $H_*(\lX;R) \isom UL_X$
where $L_X$ is the Lie algebra of Hurewicz images. 
Let $W = \bigvee_{j \in J} S^{n_j}$ be a finite-type wedge of spheres and
let $f:W \to X$.
Let $Y = X \cup_f \left(\bigvee_{j \in J} e^{n_j+1}\right)$.
Assume that $[L^W_X] \subset L_X$ is a free Lie algebra.
That is, $f$ is \emph{free}. \\
(a) Then as algebras 
\[ \gr(H_*(\lY; \F)) \isom U( {L}^X_Y \sdp \freeL (H\eL)_1)
\]
with ${L}^X_Y \isom {L}_X/[{L}^W_X]$ as Lie algebras. \\ 
(b) Furthermore if $f$ is \emph{semi-inert} then as algebras 
\[ H_*(\lY;\F) \isom U ( L^X_Y \amalg \freeL K\p )
\]
for some $K\p \subset F_1 H_*(\lY; \F)$. 
\end{thm}

\begin{proof}

\noindent
(a) Let $\AH{Y}$ be the Adams-Hilton model given above.
Therefore $H_*(\lY;\F) \isom \AH{Y}$ as algebras.
$\AH{Y}$ is a dga extension of $(\AH{X},L_X)$.
Since $[L^W_X] \subset L_X$ is a free Lie algebra, $\AH{Y}$ is a free
dga extension of $(\AH{X},L_X)$.
Thus by Theorem~\ref{thm-a} we have the algebra isomorphism
\[
\gr(H\AH{Y}) \isom U( (H\eL)_0 \sdp \freeL (H\eL)_1)
\]
with $(H\eL)_0 \isom L_X / [L^W_X]$.
Therefore 
\begin{equation} \label{eqn-f0hul} 
F_0 H\AH{Y} \isom (\gr(H\AH{Y}))_0 \isom U(H\eL)_0 \isom U(L_X/[L^W_X]). 
\end{equation}
It remains to show that $(H\eL)_0 \isom L^X_Y$.

The inclusion $i: \AH{X} \isomto F_0 \AH{Y}$ induces a map $H(i):
H\AH{X} \to F_0 H\AH{Y}$.
Now under the isomorphism \eqref{eqn-f0hul} and $UL_X \isomto
H\AH{X}$ the  map $H(i)$ corresponds to a map $UL_X \to
U(L_X/[L^W_X])$ where $U(L_X/[L^W_X]) \subset UL_Y$.
It is easy to check that this is the canonical map.
In other words $L^X_Y \isom L_X/[L^W_X]$.
Therefore $(H\eL)_0 \isom L^X_Y$.


\noindent
(b) Assume that $f$ is semi-inert.
That is, there exists an $\F$-module $K$ such that $L^X_Y \sdp \freeL
(H\eL)_1 \isom L^X_Y \amalg \freeL K$. 
Then by Theorem~\ref{thm-a}(b), $H_*(\lY;\F) \isom U (L^X_Y \amalg
\freeL K\p)$ for some $K\p \subset F_1 HU\bL = F_1 H_*(\lY;\F)$.
\end{proof}

We recall the statement of Theorem~\ref{thm-d}.

\begin{thm}[Theorem~\ref{thm-d}]
Let $R \subset \Q$ be a subring containing $\frac{1}{6}$.
Let $X$ be a finite-type simply-connected CW-complex such that
$H_*(\lX;R)$ is torsion-free and as algebras $H_*(\lX;R) \isom UL_X$
where $L_X$ is the Lie algebra of Hurewicz images. 
Let $W = \bigvee_{j \in J} S^{n_j}$ be a finite-type wedge of spheres and
let $f:W \to X$.
Let $Y = X \cup_f \left(\bigvee_{j \in J} e^{n_j+1}\right)$.
Assume that $L_X/[L^W_X]$ is $R$-free and that \feachpnP,  
$[\bar{L}^W_X] \subset \bar{L}_X$ is a free Lie algebra. 
That is, $f$ is \emph{free}. \\
(a) Then $H_*(\lY;R)$ and $\gr(H_*(\lY;R))$ are torsion-free and as algebras 
\[ \gr(H_*(\lY; R)) \isom U ( L^X_Y \sdp \freeL (H\eL)_1 )
\] 
with $L^X_Y \isom L_X/[L^W_X]$ as Lie algebras. \\ 
(b) If in addition $f$ is \emph{semi-inert} then 
as algebras 
\[ H_*(\lY; R) \isom U ( L^X_Y \amalg \freeL K\p )
\] 
for some $K\p \subset F_1 H_*(\lY; R)$. 
\end{thm}

\begin{proof}
The proof Theorem~\ref{thm-d} is exactly the same as the proof of
Theorem \ref{thm-c}, except that it uses Theorem~\ref{thm-b} instead
of Theorem~\ref{thm-a}.
\end{proof}

In the introduction we claimed that the following corollary follows
from Theorems \ref{thm-c} and \ref{thm-d}.

\begin{cor}[Corollary \ref{cor-semiInertEqnTop}]
Let $f$ be a semi-inert attaching map satisfying the hypotheses of
either Theorems \ref{thm-c} or \ref{thm-d}.
Then
\begin{equation}
K\p (z) = \tilde{H}_*(W)(z) + z[UL_X(z)^{-1} - U(L^X_Y)(z)^{-1}].
\end{equation}
\end{cor}

\begin{proof}
By Theorem~\ref{thm-c} or Theorem~\ref{thm-d},
\[ H_*(\lY;R) \isom U(L^X_Y \amalg \freeL K\p) \isom UL^X_Y \amalg
\freeT K\p.
\]
Since $(A \amalg B)(z)^{-1} = A(z)^{-1} + B(z)^{-1} - 1$~\cite[Lemma
5.1.10]{lemaire:monograph} it follows that 
\[ H_*(\lY;R)(z)^{-1} = UL^X_Y(z)^{-1} - K'(z)
\]

Let $V_1 = R\{y_j\}_{j \in J}$, where $y_j$ is given
in~\eqref{eqn-lly}.
Then $V_1(z) = \tilde{H}_*(W)(z)$.
Recall that $(H\eL)_0 \isom L^X_Y$.
Therefore, using Anick's formula~\eqref{eqn-a-formula} (with $L_0 =
L_X$) we have that 
\begin{eqnarray*}
K'(z) & = & UL^X_Y(z)^{-1} - H_*(\lY;R)(z)^{-1} \\
& = & UL^X_Y(z)^{-1} - (1+z)U(H\eL)_0(z)^{-1} + z UL_X(z)^{-1} +
V_1(z) \\
& = & \tilde{H}_*(W)(z) + z[UL^X_Y(z)^{-1} - UL_X(z)^{-1}]
\end{eqnarray*}
\end{proof}

\chapter{Constructing Spherical Hurewicz Maps} \label{chapter-hurewicz}

In this chapter we prove our other two main topological results:
Theorems \ref{thm-hurewiczF} and \ref{thm-hurewiczR}. 

Let $R = \Fp$ with $p >3$ or $R \subset \Q$ be a subring containing
$\frac{1}{6}$. 

Consider the map $W \xrightarrow{f} X$ where $W = \bigvee_{j \in J}
S^{m_j}$ is a finite-type wedge of spheres, $f = \bigvee_{j \in J} \alpha_j$,
$H_*(\lX;R)$ is torsion-free and as algebras $H_*(\lX;R) \isom UL_X$.
Let $Y$ be the adjunction space
\[ Y = X \cup_f \left( \bigvee_{j \in J} e^{m_j+1} \right).
\]
Assume that $f$ is free. That is, $[L^W_X]$ is a free Lie algebra.
Recall that $h_X: \pi_*(\lX) \tensor R \to L_X \subset H_*(\lX;R)$ is
a Lie algebra map.
Assume that there exists a Lie algebra map $\sigma_X: L_X \to
\pi_*(\lX)\tensor R$ such that $h_X \circ \sigma_X = id_{L_X}$.

If $R \subset \Q$ recall from Section~\ref{section-ip} that there is a
set of implicit primes $P_Y$ containing the invertible primes in $R$.
By replacing $R$ with $R\p=\Z[{P_Y}^{-1}]$ if necessary, we may assume
that there are no non-invertible implicit primes. 
This implies that $\forall j \in J$, $\sigma_X h_X
\hat{\alpha}_j = \hat{\alpha}_j$.

If $R = \Fp$ assume that $\forall j \in J$, $\sigma_X h_X
\hat{\alpha}_j = \hat{\alpha}_j$.

By Theorems \ref{thm-c} and \ref{thm-d}, $H_*(\lY;R)$ is torsion-free
and $\gr(H_*(\lY;R)) \isom U (L^X_Y \sdp \freeL (H\eL)_1)$ as algebras.
From this we want to show that $H_*(\lY;R) \isom UL_Y$ as algebras.

This situation closely resembles that of torsion-free spherical two-cones,
and we will generalize Anick's proof for that situation~\cite{anick:cat2}.

Recall that $\eL = L_X \amalg \freeL \langle y_j \rangle_{j \in J}$,
$\underline{\psi}: H\eL \to HU\eL$, $(\underline{\psi} H\eL)_1 \isom
(H\eL)_1$, and $HU\eL \isom \gr(H\AH{Y})$.
Let $K\p = (H\eL)_1$ and $\{\bar{w}_i\}$ be a $R$-module basis for
$K\p$.
Then $\bar{w}_i \in (H\eL)_1 = (H(L_X \amalg
\freeL \langle y_j \rangle_{j \in J}, d\p))_1$.
So each $\bar{w}_i$ is represented by a sum of brackets each with one
$y_j$ and other elements in $L_X$.
Using Jacobi identities we can write 
\begin{equation} \label{eqn-ckvk} 
\bar{w}_i = [\Sigma_{k=1}^s c_k v_k] \text{ where } c_k \in R
\text{ and } v_k = \left[ \left[ y_{j_k}, [x_{k_1}], \ldots,
[x_{k_{n_k}}] \right] \right. 
\end{equation}
where $x_{k_l} \in \AH{X}$, $[x_{k_l}] \in L_X$ and
the bracket is defined inductively by
$[[a, \ldots b, c] = [[[a, \ldots b], c]$.

Define $\gamma_i \in F_1 \AH{Y}$ by
\begin{equation} \label{eqn-ckuk} 
\gamma_i = \Sigma_{k=1}^s c_k u_k \text{ where } 
u_k = [[ y_{j_k}, x_{k_1}, \ldots, x_{k_{n_k}} ].
\end{equation}

We will use $\gamma_i$ and Adams-Hilton models to construct a map
$g_i: S^{m+1} \rightarrow Y$, whose Hurewicz image `modulo lower
filtration' is $\bar{w}_i$.

The following geometric construction is a slight generalization 
of~\cite[Proposition 5.4]{anick:cat2} whose proof is essentially the same.
It is the central construction of this chapter and will give us the
desired map $S^{m+1} \rightarrow Y$.

\begin{prop} \label{prop-sphericalHurewiczMap}
Let $W \rightarrow X \rightarrow Y$ be as above.
Let $\gamma = \Sigma_{k=1}^s c_k u_k$ as in~\eqref{eqn-ckuk} with
$[x_{k_i}] \in L_X$. 
Then there exists a map 
\[ g: S^{m+1}_s \rightarrow Y \]
which has an AH model $\AH{g}$ satisfying $\AH{g}(b_k) = c_k u_k$ for 
$1\leq k \leq s$, and $\AH{g}(b_0) \in \AH{X}$.
\end{prop}

\begin{proof}
By Lemma~\ref{lemma-ah} there exist maps 
$g_k: (D^{m+1}_0,S^m_0) \rightarrow (Y,X)$ for $1 \leq k \leq s$ with
models $\AH{g_k}(ii) = c_k u_k$.
In addition 
\[g_k|_{S^m_0} = \epsilon_k c_k [[\alpha_{j_k}, \beta_{k_1}, \ldots
\beta_{k_{n_k}}] \] 
where $\epsilon_k = \pm 1$ and $\beta_i: S^{n_i+1} \rightarrow X$ 
is the adjoint of $\sigma([x_i]): S^{n_i} \rightarrow \lX$
for $1\leq k \leq s$, 
and $\psi_{g_k}(i) \in CU_*(\Omega X)$ for $1\leq k\leq s$.
Let \[ g' = g_1 \vee \ldots \vee g_s : \left(\bigvee_{k=1}^{s} D^{m+1}_0, 
\bigvee_{k=1}^{s} S^m_0\right) \rightarrow (Y,X) \]
for which one can choose $\AH{g'} (b_k) = c_k u_k$ for $1\leq k \leq s$. 

Restricting $g'$ we get a map $g_0 : \left( \bigvee_{k=1}^{s} S^m_0 \right) \rightarrow X$,
and $(\AH{g'},\psi_{g'})$ extends a valid AH model $(\AH{g_0},\psi_{g_0})$
for $g_0$.

We will show that $g_0$ can also be extended to a map
\[ g'': \left( \bigvee_{k=1}^{s} S^m_0 \right) \cup_{\Sigma_{k=1}^{s}
\iota_k} e^{m+1}  \rightarrow X,\]
where $\iota_k$ is the inclusion of the $k$th sphere into the wedge.
It will follow at once that there is an AH model for $g''$ which extends
$(\AH{g_0},\psi_{g_0})$ and has $\AH{g''}(b_0) \in \AH{X}$, 
$b_0$ denoting the $m$-dimensional generator of 
$\AH{\bigvee_{k=1}^{s} S^m_0 \cup_{\Sigma_{k=1}^{s} \iota_k} e^{m+1}}$. 

To prove the existence of $g''$, it suffices to show that 
$g_{0\#} (\Sigma_{k=1}^{s}\iota_k)$ vanishes in $\pi_m(X)$.
\[ g_{0\#}(\Sigma_{k=1}^{s}\iota_k) = \Sigma_{k=1}^s \epsilon_k c_k 
[[\alpha_{j_k}, \beta_{k_1}, \ldots \beta_{k_{n_k}}] \] 
Recall from \eqref{eqn-ckvk} that $\bar{w}_i = [ \Sigma_{k-1}^s c_k
v_k]$.
Since $\Sigma_{k=1}^s c_k v_k$ is a cycle in $\eL$,
\begin{eqnarray*} 
0 & = & d \Sigma_{k=1}^s c_k v_k \\
& = & \Sigma_{k=1}^s c_k \left[ \left[ d\p y_{j_k}, [x_{k_1}], \ldots
[x_{k_{n_k}}] \right] \right. \\
& = & \Sigma_{k=1}^s c_k \left[ \left[ h_X(\hat{\alpha}_{j_k}),
[x_{k_1}], \ldots [x_{k_{n_k}}] \right] \right. . 
\end{eqnarray*}
Since $\sigma_X h_X \hat{\alpha}_j = \hat{\alpha}_j$ applying the map
$\sigma_X$ gives 
\[0 = \Sigma_{k=1}^s c_k \left[ \left[ \hat{\alpha}_{j_k},
\sigma_X([x_{k_1}]), \ldots, \sigma_X([x_{k_{n_k}}]) \right] \right. 
\] 
and adjointing gives that $g_{0\#}(\Sigma_{k=1}^s \iota_k) = 0$.

Now $g'$ and $g''$ are compatible extensions of $g_0$, so together they
define a map $g: S^{m+1}_s \rightarrow Y$.
Furthermore the corresponding AH models are compatible so they give
a valid AH model for $g$ with the desired properties.
\end{proof}

We review and continue our construction.
Starting with $\bar{w}_i \in (H\eL)_1$ we choose $\gamma_i$ as above
and use Proposition~\ref{prop-sphericalHurewiczMap} to construct $g_i: 
S^{m+1}_s \to Y$.  
We then have the following map
\begin{equation} \label{eqn-phi}
\phi_i: S_0^{m+1} \stackrel{\Psi}{\to} S_s^{m+1} \stackrel{g_i}{\to} Y 
\end{equation}
where $\Psi: S^{m+1}_0 \rightarrow S^{m+1}_s$ is the homeomorphism
from~\eqref{eqn-Psi}.
Taking $\AH{\phi_i} = \AH{g_i} \circ \AH{\Psi}$ we have
$\AH{\phi_i} (i) = \Sigma_{k=1}^s c_k u_k - \AH{g_i}b_0$.
Letting $\lambda_i = \AH{g_i}(b_0)$ we have $\AH{\phi_i}(i) = \gamma_i -
\lambda_i$ and hence $h_Y(\phi_i) = [\gamma_i - \lambda_i] \in L_Y$.

Furthermore $\gamma_i \in F_1 \AH{Y}$ and $\lambda_i \in \AH{X} =
F_0 \AH{Y}$.
Since $L_Y$ inherits a filtration from $H\AH{Y}$ which is in turn
induced by the filtration on $\AH{Y}$, $[\gamma_i - \lambda_i] \in
F_1 L_Y$.
Recall the quotient map~\eqref{eqn-quotientf} $f: F_1(H\AH{Y}) \to
\gr(H\AH{Y})_1 \isomto (HU\eL)_1$.
By construction $f([\gamma_i - \lambda_i]) = \bar{w}_i$.

Recall that $K\p = R\{\bar{w}_i\}$.
Therefore we have an injection $K\p \incl (\gr(L_Y))_1$.
Furthermore by construction this is a map of $L^X_Y$-modules.
Thus we have proved the following.

\begin{prop} \label{prop-u1}
Let $W \to X \to Y$ and $K\p$ be as above.
Then there exists an injection of $L^X_Y$-modules
\[ u_1: K\p \incl (\gr(L_Y))_1.
\]
\end{prop}

\ignore
{
\begin{proof}
Also recall from Section~\ref{section-ss} that the isomorphism
$\theta$ is the isomorphism between the $E^{\infty}$ and $E^2$ terms
of the spectral sequence $\gr(\AH{Y}) \implies \gr(H\AH{Y})$.
By the definitions of $v_k$ and $u_k$ in \eqref{eqn-ckvk} and
\eqref{eqn-ckuk} $[\Sigma c_k v_k]$ and $f([\Sigma c_k u_k -
\lambda])$ are both represented by the same cycle $\gr(w_i) \in
\gr(\AH{Y})$. 
Therefore we have that
\[ f([\gamma_i - \lambda_i]) = f([\Sigma c_k u_k - \lambda_i]) =
[\Sigma c_k v_k] = \bar{w}_i.
\]
Since $[\gamma_i - \lambda_i] = h_Y(\phi_i) \in L_Y$ this proves the
lemma.
\end{proof}
}

We are now ready to prove Theorems \ref{thm-hurewiczF} and
\ref{thm-hurewiczR}. For convenience we restate the theorems.

\begin{thm}[Theorem~\ref{thm-hurewiczF}]
Let $R = \F$ where $\F = \Q$ or $\Fp$ with $p >3$.
Let $\bigvee_{j \in J} S^{n_j} \xto{\bigvee \alpha_j} X$ be a cell
attachment satisfying the hypotheses of Theorem~\ref{thm-c}.
Let $Y = X \cup_{\bigvee \alpha_j} \left(\bigvee e^{n_j+1}\right)$
and let $\hat{\alpha}_j$ denote the adjoint of $\alpha_j$. 
In addition assume that there exists a map $\sigma_X$ right inverse to
$h_X$ and that $\forall j \in J$, $\sigma_X h_X \hat{\alpha}_j =
\hat{\alpha}_j$. 
Then the canonical algebra map
\[ UL_Y \to H_*(\lY;\F)
\] 
is a surjection.
That is, $H_*(\lY;\F)$ is generated as an algebra by Hurewicz images.
\end{thm}

\begin{proof}
Let $R = \Fp$ where $p > 3$.
%
Let $g: \gr(H_*(\lY;R)) \isomto UL\p$ be the algebra isomorphism given by
Theorem~\ref{thm-c}(i) where $L\p = L^X_Y \sdp \freeL K\p$ with $K\p =
(H\eL)_1$. 
Note that $L\p_0 = L^X_Y$ and that $L\p_1 = K\p$.
We will show that the canonical map $UL_Y \to H_*(\lY;R)$ is surjective.

We have an injection of Lie algebras 
\begin{equation} \label{eqn-l0'}
u_0: L^X_Y \incl F_0 L_Y \isomto (\gr(L_Y))_0.
\end{equation}

Since $\forall j \ h_X \sigma_X \hat{\alpha}_j = \hat{\alpha}_j$, by
Proposition~\ref{prop-u1} and get an injection of $L^X_Y$-modules 
$u_1: K\p \incl (\gr(L_Y))_1$.
So for $x \in L^X_Y$ and $y \in K\p$, $u_1([y,x]) = [u_1(y), u_0(x)]$.
Thus by Lemma~\ref{lemma-sdp-extension} $u_{0}$ and $u_{1}$ can be
extended to a Lie algebra map $u: L\p \rightarrow \gr(L_Y)$.

The inclusion $L_Y \hookrightarrow H_*(\lY;R)$ induces a map between 
the corresponding graded objects, $\chi: \gr(L_Y) \rightarrow
\gr(H_*(\lY;R))$.

We claim that for $j=0$ and $1$, $g\circ \chi \circ u_{j}$ is the ordinary 
inclusion of $L\p_{j}$ in $UL\p$.
By~\eqref{eqn-l0'} $u_0$ is an injection.
In addition, restricted to degree $0$, $g$ is the identity and $\chi$
is just the ordinary inclusion.
For $j=1$, $g \chi u_1 \bar{w}_i = g \gr([\gamma_i - \lambda_i]) =
f([w_i]) = \bar{w}_i$. 
It follows that $g \circ \chi \circ u$ is the standard inclusion $L\p
\incl UL\p$. 
Since $g \circ \chi \circ u$ is an injection, so is $u$.

By Lemma~\ref{lemma-Ugr=grU} the canonical map $U\gr(L_Y) \isomto
\gr(UL_Y)$ is an algebra isomorphism. 
Now $u$ and $\chi$ induce the maps $Uu$ and  $\bar{\chi}$
in the following diagram.
\[ \xymatrix{
UL\p \ar[r]^-{Uu} \ar[dr]^{\isom}_{g^{-1}} & U\gr(L_Y) \ar[r]^-{\cong} 
\ar[d]^{\bar{\chi}} & \gr(UL_Y) \ar@{-->}[dl]^{\tilde{\chi}} \\
& \gr(H_*(\lY;R)) }
\]
Since we showed that $g \chi u$ is the ordinary inclusion $L\p \incl
UL\p$ the diagram commutes.
Since $g^{-1}$ is surjective, the induced map $\tilde{\chi}$ is surjective.
Since $\tilde{\chi}$ is the associated graded map to the canonical map
$UL_Y \to H_*(\lY;R)$ and the filtrations are bicomplete, the associated
ungraded map is also surjective. 
So the canonical map $UL_Y \to H_*(\lY;R)$ is surjective which
finishes the proof. 
\end{proof}

The proof of Theorem~\ref{thm-hurewiczF} is still valid in the case
where $R \subset \Q$ containing $\frac{1}{6}$. 
However in this case we can tensor with $\Q$ and make use of results
from rational homotopy theory. 
As a result Theorem~\ref{thm-hurewiczR} is a stronger version of
Theorem~\ref{thm-hurewiczF}.

\begin{thm}[Theorem~\ref{thm-hurewiczR}]
Let $R \subset \Q$ be a subring containing $\frac{1}{6}$.
Let $\bigvee_{j \in J} S^{n_j} \xto{\bigvee \alpha_j} X$ be a cell
attachment satisfying the hypotheses of Theorem~\ref{thm-d}.
Let $Y = X \cup_{\bigvee \alpha_j} \left(\vee e^{n_j+1}\right)$.
Furthermore assume that there exists a map $\sigma_X$ right inverse to
$h_X$.
Let $P_Y$ be the set of \emph{implicit primes} and let $R\p =
\Z[{P_Y}^{-1}]$. 
Then \\
(i)  $H_*(\lY;R\p)$ is torsion-free and as algebras 
\[ H_*(\lY;R\p) \isom UL_Y
\] 
where $\gr(L_Y) \isom L^X_Y \sdp \freeL (H\eL)_1$ as Lie algebras, and \\ 
(ii) localized away from $P_Y$, $\lY \in \PiS$. \\
(iii) If in addition $f$ is semi-inert then localized away from $P_Y$,
$L_Y \isom L^X_Y \amalg \freeL \hat{K}$ as Lie algebras for some
$\hat{K} \subset F_1 L_Y$, and 
there exists a map $\sigma_Y$ right inverse to $h_Y$. 
\end{thm}

\begin{rem}
Using Remark~\ref{rem-hurewiczF} we could have tried proving
Theorem~\ref{thm-hurewiczR}(i) using Theorem~\ref{thm-hurewiczF} and
the Hilton-Serre-Baues Theorem (Theorem~\ref{thm-hsb}), however this
is not the approach we use here.
\end{rem}

\begin{proof}[Proof of Theorem~\ref{thm-hurewiczR}]
Let $R \subset \Q$ containing $\frac{1}{6}$.

\noindent
(i) Since there are no non-invertible implicit primes we have that
$\sigma_X h_X \hat{\alpha}_j = \hat{\alpha}_j$.
As a result we can copy the proof of Theorem~\ref{thm-hurewiczF}
verbatim, except we use Theorem~\ref{thm-d}(i) instead of
Theorem~\ref{thm-c}(i). 

Recall that $\gr(H_*(\lY;R)) \isom UL\p$ and that we constructed a Lie
algebra map $u: L\p \to \gr(L_Y)$ and showed that it is an injection.
We claim that for $R \subset \Q$, $u$ is an isomorphism. 
$H_*(\lY;R)$ and $\gr(H_*(\lY;R))$ have the same Hilbert series.
Also since $H_*(\lY; R)$ is torsion-free, it has the same Hilbert
series as $H_*(\lY; \Q)$.
Let $S$ be the image of $h_Y \tensor \Q$. 
Then $S$, $L_Y$ and $\gr(L_Y)$ have the same Hilbert series.
By the Milnor-Moore Theorem (Theorem~\ref{thm-milnorMoore}), $H_*(\lY;
\Q) \isom US$. 
So by the Poincar\'{e}-Birkhoff-Witt Theorem, $S$ has the same Hilbert
series as $L\p$, and hence $\gr(L_Y) \isom L\p$ as $R$-modules. 
Since $u: L\p \to \gr(L_Y)$ is an injection it follows that it is an
isomorphism.  

By Lemma~\ref{lemma-Ugr=grU} the canonical map $U\gr(L_Y) \isomto
\gr(UL_Y)$ is an algebra isomorphism. 
Now $u$ and $\chi$ induce the maps $Uu$ and  $\bar{\chi}$
in the following diagram.
\[ \xymatrix{
UL\p \ar[r]^-{Uu} \ar[dr]^{\isom}_{g^{-1}} & U\gr(L_Y) \ar[r]^-{\cong} 
\ar[d]^{\bar{\chi}} & \gr(UL_Y) \ar@{-->}[dl]^{\tilde{\chi}} \\
& \gr(H_*(\lY;R)) }
\]
Since we showed that $g \chi u$ is the ordinary inclusion $L\p \incl
UL\p$ the diagram commutes.
Since $u$ is an isomorphism, so is $Uu$.
Hence the induced map $\tilde{\chi}$ is an isomorphism.
Since $\tilde{\chi}$ is the associated graded map to the canonical map
$UL_Y \to H_*(\lY;R)$ and the filtrations are bicomplete, the associated
ungraded map is also an isomorphism. 
So $H_*(\lY;R) \isom UL_Y$ as algebras.

The map $u$ gives the desired Lie algebra isomorphism $\gr(L_Y) \isom
L\p = L^X_Y \sdp \freeL (H\eL)_1$, which finishes the proof of (i).

\noindent
(ii) is equivalent to (i) by the Hilton-Serre-Baues Theorem (Theorem~\ref{thm-hsb}).
It remains to prove (iii).

\noindent
(iii) Before we construct the desired map $\sigma_Y$ we strengthen the
result in (i) in the semi-inert case. 

Assume that $f$ is semi-inert.
Recall the situation from Theorem~\ref{thm-d}(ii).
We have that $\gr(H_*(\lY;R)) \isom UL\p$ where $L\p \isom L^X_Y
\amalg \freeL K\p$ where $K\p \isom R\{\bar{w}_i\} \subset
(\gr(H_*(\lY;R)))_1$.
For each $\bar{w}_i$ let $[w_i]$ be an inverse image under the
quotient map $H_*(\lY;R) \to \gr(H_*(\lY;R))$.
Let $K\pp = \{[w_i]\}$.
By Theorem~\ref{thm-d}(ii) as algebras $H_*(\lY;R) \isom UL\pp$ where
$L\pp = L^X_Y \amalg \freeL K\pp$ (see Theorem~\ref{thm-b'}(b)).
Since $K\pp \xrightarrow{\isom} K\p$, there is an induced Lie algebra
isomorphism $L\pp \xrightarrow{\isom} L\p$.
So $\gr(H_*(\lY;R)) \isom H_*(\lY;R)$ as algebras. 

Let 
\[ \hat{L} = L^X_Y \amalg \freeL \hat{K} \text{ where } \hat{K}=R\{[\gamma_i
- \lambda_i]\} \subset F_1 L_Y.
\]
Recall that $[w_i] - [\gamma_i - \lambda_i] \in F_0 (H_*(\lY;R))$ and that
$f([\gamma_i - \lambda_i]) = \bar{w}_i$ where $f$ is the quotient map
from~\eqref{eqn-quotientf}.
So $f: \hat{K} \isomto K\p$, which induces a Lie algebra isomorphism
$\hat{L} \isomto L\p$.
This in turn induces the algebra isomorphism $U\hat{L} \isomto UL\p
\isom H_*(\lY;R)$. 

Since $U\hat{L} \isom H_*(\lY;R)$ there is an injection $L_Y \incl
U\hat{L}$.
Also, since $\hat{K} \subset L_Y$ there is a canonical Lie algebra map
$u: \hat{L} \to L_Y$.
These fit into the following commutative diagram.
\[
\xymatrix{
\hat{L} \ar[rr] \ar[dr]_u & & U\hat{L} \\
& L_Y \ar[ur]
}
\]
It follows that $u$ is an injection.

We claim that $u$ is an isomorphism. 
Since $H_*(\lY; R)$ is torsion-free, it has the same Hilbert
series as $H_*(\lY; \Q)$.
Let $S$ be the image of $h_Y \tensor \Q$. 
Then $S$ and $L_Y$ have the same Hilbert series.
By the Milnor-Moore Theorem (Theorem~\ref{thm-milnorMoore}), $H_*(\lY;
\Q) \isom US$. 
So by the Poincar\'{e}-Birkhoff-Witt Theorem (Theorem~\ref{thm-pbw}),
$S$ has the same Hilbert series as $\hat{L}$, and hence $L_Y \isom
\hat{L}$ as $R$-modules.  
Since $u: \hat{L} \to L_Y$ is an injection it follows that it is an
isomorphism.  

Therefore $L_Y \isom L^X_Y \amalg \freeL \hat{K}$ as Lie algebras,
with $\hat{K} \subset F_1 L_Y$ and
$H_*(\lY;R) \isom UL_Y$ as algebras. 

\begin{rem}
If we generalize a conjecture of Anick~\cite[Conjecture 4.9]{anick:cat2}
from spherical two-cones to our adjunction space $Y$ then the Lie
algebra isomorphism $L_Y \isom L\p$ proves the conjecture in the
semi-inert case. 
\end{rem}

We are now in a position to construct a map $\sigma_Y$ right inverse
to $h_Y$.
Let $i$ denote the inclusion $X \to Y$.
Consider the composite map 
\[ F: [L^W_X] \incl L_X \xrightarrow{\sigma_X} \pi_*(\lX) \tensor R
\xrightarrow{(\Omega i)_{\#}} \pi_*(\lY) \tensor R.
\]
We claim that $F=0$. Since $F$ is a Lie algebra map it is sufficient
to show that $F(L^W_X) = 0$. 
\[ F(L^W_X) = (\Omega i)_{\#} \sigma_X (R\{h_X(\hat{\alpha}_j)\}).
\]
Since there are no implicit primes $\sigma_X h_X \hat{\alpha}_j =
\hat{\alpha}_j$. 
By the construction of $Y$, $\Omega i \circ \hat{\alpha}_j \simeq 0$.
So $F=0$ as claimed.

Therefore there is an induced map $G: L^X_Y \isom L_X/[L^W_X] \to
\pi_*(\lY) \tensor R$. 
That $h_Y \circ G$ is the inclusion map can be seen from the following
commutative diagram. 
\[
\xymatrix{
L_X \ar@(ur,dr)@{.>}[r]^-{\sigma_X} \ar[d] & \pi_*(\lX) \tensor R \ar[dd]^{(\Omega
i)_{\#}} \ar[l]^-{h_X} \\ 
L^X_Y \ar[dr]^{G} \ar[d] \\
L_Y & \pi_*(\lY) \tensor R \ar[l]_-{h_Y}
}
\]

Now construct $\sigma_Y: L_Y \to \pi_*(\lY) \tensor R$ as follows.
We have shown that $L_Y \isom L^X_Y \amalg \freeL K$ for some $K
\subset F_1 L_Y$.  
Since $F_1 L_Y \isom \hat{K} = R\{[\gamma_i - \lambda_i]\}$ it follows
that $K = R\{[\gamma_l + \lambda_l]\}_{l \in L}$ for some $L$.
Recall that $\exists \ \phi_l \in \pi_*(\lY) \tensor R$ such that
$h_Y(\phi_l) = [\gamma_l + \lambda_l]$. 
Let $\sigma_Y \!\! \mid_{L^X_Y} = G$ and let $\sigma_Y([\gamma_l +
\lambda_l]) = \phi_l$.
Now extend $\sigma_Y$ canonically to a Lie algebra map on $L_Y$.

We finally claim that $h_Y \sigma_Y = id_{L_Y}$.
Since $h_Y \sigma_Y$ is a Lie algebra map it suffices to check that it
is the identity for the generators. 
\[ h_Y \sigma_Y L^X_Y = h_Y G L^X_Y = L^X_Y, \quad
h_Y \sigma_Y [\gamma_l + \lambda_l] = h_Y \phi_l = [\gamma_l +
\lambda_l]
\]
Therefore $\sigma_Y$ is the desired Lie algebra map right inverse to $h_Y$.
\end{proof}

\part{Applications}
\chapter{Free Subalgebras} \label{chapter-freeLa}

In this chapter we prove Theorem~\ref{thm-schreier2} and give a second
proof of a special case of this theorem.

\section{A generalized Schreier property}

In chapter~\ref{chapter-HUL} we calculated the homology of $U\eL =
U(L_0 \amalg \freeL V_1, d)$, a two-level
universal enveloping algebra over a field $\F$ in which the Lie ideal
$[dV_1] \subset L_0$ was a free Lie algebra. 
The result was a universal enveloping algebra on a Lie algebra 
\[ L = L_0 \sdp \freeL L_1.
\]
We remark that an equivalent description of $L$ is the fact that it
has homogeneous generators and relations in degrees $0$ and $1$.

In this chapter we give a simple criterion which we prove guarantees
that a Lie subalgebra of $L$ is free. 

If $R \subset \Q$ is a subring containing $\frac{1}{6}$,
Theorem~\ref{thm-hurewiczR} shows that under certain conditions,
$H_*(\lY;R\p)$ is torsion-free, as algebras  $H_*(\lY;R\p) \isom UL_Y$
and as Lie algebras $\gr(L_Y) \isom L^X_Y \sdp \freeL V$ for some free
$R\p$-module $V$. 
Let $p$ be a noninvertible prime in $R\p$ and recall the notation
$\bar{M} = M \tensor \Fp$ for any $R\p$-module $M$. 
It follows that $\gr(\bar{L}_Y) \isom \bar{L}^X_Y \sdp \freeL
\bar{V}$.

We will give a simple criterion guaranteeing that a Lie subalgebra $J
\subset \bar{L}_Y$ is a free Lie algebra.

It is a well-known fact that any (graded) Lie subalgebra of a (graded)
Lie algebra is a free Lie algebra
\cite{mikhalev:subalgebras, shirshov:subalgebras, mikhalevZolotykh}
(this is referred to as the \emph{Schreier property}).
In this chapter we generalize this result to the following.

\begin{thm}[Theorem~\ref{thm-schreier2}]
Over a field $\F$, let $L$ be a finite-type graded Lie algebra with filtration
$\{F_k L\}$ such that $\gr(L) \isom L_0 \sdp \freeL V_1$ as Lie
algebras, where $L_0 = F_0 L$ and $V_1  = F_1 L / F_0 L$.
Let $J \subset L$ be a Lie subalgebra such that $J \cap F_0 L = 0$.
Then $J$ is a free Lie algebra.
\end{thm}

Before proving this theorem,  we prove the following lemma.

\begin{lemma} \label{lemma-schreier}
Let $J$ be a finite-type filtered Lie algebra such that $\gr(J)$ is a
free Lie algebra.
Then $J$ is a free Lie algebra.
\end{lemma}

\begin{proof}
By assumption there is an $\F$-module $\bar{W}$ such that $\gr(J)
\isom \freeL \bar{W}$.

Let $\{\bar{w}_i\}_{i \in I} \subset \gr(J)$ be an $\F$-module basis
for $\bar{W}$.
Let $m_i = \deg(\bar{w}_i)$.
That is, $\bar{w}_i \in F_{m_i}J / F_{m_i-1} L$.
For each $\bar{w}_i$ choose a representative $w_i \in F_{m_i} J$.
Let $W = \F\{w_i\}_{i \in I} \subset J$.

Then there is a canonical map $\phi: \freeL W \to J$.
Grade $\freeL W$ by letting $w_i \in W$ be in degree $m_i$.
Then $\phi$ is a map of filtered objects and there is an induced map
$\theta: \freeL W \to \gr(J)$.
However the composite map
\[ \freeL W \xto{\theta} \gr(J) \isomto \freeL \bar{W}
\]
is just the canonical isomorphism $\freeL W \isomto \freeL \bar{W}$.
So $\theta$ is an isomorphism.

\ignore
{
Let $\tilde{W}$ be a the graded $\F$-module on a basis
$\{\tilde{w}_i\}_{i \in I}$ where $|\tilde{w}_i| = |\bar{w}_i| =
|w_i|$.
Let $f(\tilde{w}_1, \ldots \tilde{w}_n) \in (\freeL \tilde{W})_m$.

Then $f(w_1, \ldots w_n) \in F_m J$.
If $f(w_1, \ldots w_n) = 0$ then $f(w_1, \ldots w_n) \in F_{m-1}J$ and
hence $0 = f(\bar{w}_1, \ldots \bar{w}_n) \in \gr(J)$.
But this contradicts the fact that $\gr(J) \isom \freeL \bar{W}$.
Therefore $f(w_1, \ldots w_n) \neq 0$.

Thus the canonical map $\phi: \freeL W \to J$ is injective. 
As $\F$-modules $\freeL W \isom \freeL \bar{W} \isom \gr(J) \isom J$.
}

Therefore $\phi$ is an isomorphism and $J$ is a free Lie algebra.
\end{proof}

\begin{proof}[Proof of Theorem~\ref{thm-schreier2}]
The filtration on $L$ filters $J$ by letting 
\[ F_k J = J \cap F_k L.
\]
From this definition it follows that the inclusion $J \incl L$ induces
an inclusion $\gr(J) \incl \gr(L)$.
So $\gr(J) \incl \gr(L) \isom L_0 \sdp \freeL V_1$.
Since $J \cap F_0 L = 0$ it follows that $(\gr J)_0 = 0$ and $\gr(J)
\incl (\gr J)_{\geq 1} \isom \freeL V_1$.  
By the Schreier property $\gr(J)$ is a free Lie algebra.

Thus by Lemma~\ref{lemma-schreier}, $J$ is a free Lie algebra.
\end{proof}

The following corollary is a special case of this theorem.

\begin{cor} \label{cor-schreier}
Over a field $\F$, if $J \subset L_0 \sdp \freeL (L_1)$ is a Lie
subalgebra such that $J \cap L_0 = 0$ then $J$ is a free Lie algebra.\
\end{cor}

Note that since $J$ is not necessarily homogeneous with respect to
degree $J \cap L_0 = 0$ does not imply that $J \subset \freeL L_1$.

\section{A second proof}

%

In this section we give an independent proof
Corollary~\ref{cor-schreier} which does not use the Schreier property.
As such, it proves the Schreier property as a special case.

\begin{thm} \label{thm-freeJinL}
Let $L = L_0 \sdp \freeL L_1$ be a finite-type bigraded Lie algebra
over a field $\F$ with generators and relations in degrees $0$ and $1$.
Let $J$ be a graded subalgebra of $L$ not necessarily graded with
respect to degree, such that $J \cap L_0 = 0$.
Then $J$ is a free Lie algebra.
\end{thm}

Let $J$ be a Lie subalgebra of $L$ generated by elements which are
homogeneous with respect to dimension, but not necessarily with
respect to degree. So $J$ is graded by dimension and filtered by
degree. Let $J_m$ be the subspace of $J$ in degrees $\leq m$ and
$J_{m,n}$ the component of $J_m$ in dimension $n$.
Let $J_{m,<n} = \bigoplus_{k<n} J_{m,k}$.

Let $\{x_i\}$ be a set of generators homogeneous with respect to
dimension for a filtered Lie algebra $J$. 
Let $J^{j}$ denote the Lie subalgebra of $J$ generated by $\{x_i\}_{i
\neq j}$.
The set $\{x_i\}$ is called \emph{reduced} if for each generator 
$x_i \in J_{m,n} \backslash J_{m-1,n}$, 
\begin{equation} \label{eqn-reduced}
\nexists \: y \in J^i_{*,n} \: \mid \: x_i - y \in J_{m-1,n}.
\end{equation}

\begin{rem}
Note that reduced $\implies$ \emph{minimal}.
That is, for each $x_i$,
$\nexists \: y \in J^i$ such that $x_i -y = 0$.
We also remark that the standard definition of a reduced
set~\cite{mikhalevZolotykh} 
(in the context where $L$ is a free Lie algebra) is the condition that
none of the $x\p_j$ lie in the Lie subalgebra of $L$ generated by
$\{x\p_i\}_{i\neq j}$.
We will see in the proof of Lemma~\ref{lemma-Nielsen} that this condition
follows from our condition~\eqref{eqn-reduced}.
\end{rem}

\begin{lemma} \label{lemma-reduced}
A connected, graded, filtered, finite-type Lie algebra $J$ has a
reduced set of generators.
\end{lemma}

\begin{proof}
We deduce by induction on dimension.
Assume we have a reduced set of generators homogeneous with respect to
dimension for the
Lie subalgebra generated by $J_{*,<n}$.
Since $J$ has finite-type, $J_{*,n}$ is a finite vector space. 
Therefore $J_{*,n} \isom J_{M,n}$ for some $M$.

We now deduce by induction on degree.
Assume we have a reduced set of generators homogeneous 
with respect to dimension for the
Lie subalgebra generated by $J_{*,<n} \oplus J_{m-1,n}$.
Choose a basis $\{y_j\}$ for the finite vector space 
$J_{m,n} \ominus J_{m-1,n}$.

We finally deduce by induction on $j$.
Assume we have a reduced set of generators $\{x_i\}$ homogeneous with
respect to dimension, for the
Lie subalgebra generated by $V := J_{*,<n} \oplus J_{m-1,n} \bigoplus_{j<k}
\F\{y_j\}$, which we denote $\langle x_i \rangle$. 
More generally we will let $\langle S \rangle$ denote the Lie subalgebra
of $J$ generated by the set $S$.

Consider $y_k$.
If $\exists y \in \langle x_i \rangle_{*,n}$ such that $y_k -y \in
J_{m-1,n}$ then since $J_{m-1,n} \subset \langle x_i \rangle$,
$\exists x \in \langle x_i \rangle$ such that $y_k -y = x$.
Therefore $y_k = x + y \in \langle x_i \rangle$, and $\F\{y_k\} \subset
\langle x_i \rangle$.
Thus $\{x_i\}$ is a reduced set of generators for the Lie subalgebra
generated by $V \oplus \F\{y_k\}$.

If $\nexists$ such a $y$ we claim that $\{x_i\}_{i\in I}
\cup \{y_k\}$ is a reduced set of generators for the Lie subalgebra
generated by $V \oplus \F\{y_k\}$.
Clearly it is a set of generators; it remains to show that it is
reduced.
We need to check that each generator satisfies
condition~\eqref{eqn-reduced}. 

By assumption, $y_k$ satisfies the condition.
By the induction hypothesis and since $y_k \in J_{*,n}$, all $x_i \in
J_{*,<n}$ satisfy the condition.
It remains to prove that each $x_i \in J_{m,n}$ satisfies
condition~\eqref{eqn-reduced}.

\nopagebreak[1]
Assume $\exists j$, and $\exists y \in \langle \{x_i\}_{i \neq j} \cup \{y_k\}
\rangle_{*,n}$ such that $x_j \in J_{l,n} \backslash J_{l-1,n}$ (where $l \leq
m$) and $x_j -y \in J_{l-1,n}$. 
For degree reasons, $y \in J_{l,n} \subset J_{m,n}$ and since $y_k \in J_{m,n}$
and $J$ is connected,
$y = x + a y_k$, where $x \in \langle \{x_i\}_{i\neq j}
\rangle_{*,n}$ and $a \in \F$.
Therefore $x_j -x - ay_k \in J_{l-1,n}$.
By assumption $a \neq 0$.
Thus $y_k - a^{-1}(x_j - x) \in J_{l-1,n} \subset J_{m-1,n}$.
But $a^{-1}(x_j-x) \in \langle x_i \rangle_{*,n}$, which contradicts
the assumption that $y_k$ satisfies condition~\eqref{eqn-reduced}.
Hence $V \oplus \F\{y_k\}$ has a reduced set of generators.

So by induction, the Lie subalgebra generated by $J_{*,<n} \oplus J_{m,n}$
has a reduced set of generators.
Then by induction the Lie subalgebra generated by $J_{*,\leq n}$ has a
reduced set of generators.
Finally by induction $J$ has a reduced set of generators.
\end{proof}

Let $J$ be a Lie subalgebra of the semi-direct product $L = L_0 \sdp
\freeL(L_1)$ such that $J \cap L_0 = 0$.
$L$ is bigraded. 
Let $L_j$ denote the component of $L$ in degree $j$.
By the previous lemma, $J$ has a reduced set of generators
$\{x_i\}$ homogeneous with respect to dimension, where $x_i =
\Sigma_{j=0}^{n_i} x_i^{(j)}$, with $x_i^{(j)} \in L_j$, 
$x_i^{(n_i)} \neq 0$ and $n_i \geq 1$.
Let $x\p_i = x^{(n_i)}_i$.
Note that the $x_i^{(j)}$ are homogeneous with respect to degree.

\begin{lemma} \label{lemma-Nielsen}
There does not exist a nonzero $f(y_1, \ldots y_n) \in \freeL \langle
y_1, \ldots y_n \rangle$ with $|y_i| = |x_i|$ such that $f(x\p_1,
\ldots x\p_n) = 0$.
\end{lemma}

\begin{proof}
Assume such an $f$ exists.
Let $(a_1, \ldots, \hat{a}_j, \ldots a_n)$ denote the list $(a_1,
\ldots, a_n)$ with $a_j$ omitted. 
Since $n_i \geq 1$, $x\p_i \in L_{\geq 1} \isom \freeL(L_1)$.
Since free graded Lie algebras satisfy the Nielsen 
property \cite[Theorem 14.1]{mikhalevZolotykh} and each $x\p_i$ is
homogeneous with respect to degree,
$\exists \ 0 \neq g(y_1, \ldots \hat{y}_j, \ldots y_n) \in 
\freeL \langle y_1, \ldots \hat{y}_j, \ldots y_n \rangle$ with $|y_i| =
|x_i|$ such that $x\p_j = g(x\p_1, \ldots \hat{x\p}_j, \ldots x\p_n)$
for some $j$.
Without loss of generality, we can assume that $g(x\p_1, \ldots
\hat{x\p}_j, \ldots x\p_n)$ is a sum of nonzero monomials in dimension
$|x_j|$ and degree $n_j (= \deg(x\p_j))$.

By linearity, 
\begin{align*}
g(x_1, \ldots \hat{x}_j, \ldots x_n) &= \alpha + g(x\p_1, \ldots
\hat{x\p}_j, \ldots x\p_n), && \alpha \in J_{n_j-1} \\
&= \alpha + x\p_j, && x\p_j = x_j - \Sigma_{k=0}^{n_j-1} x_j^{(k)} \\
&= \beta + x_j, && \beta \in J_{n_j-1}.
\end{align*}
Therefore $x_j - g(x_1, \ldots \hat{x}_j, \ldots x_n) \in J_{n_j-1}$.
But this contradicts the assumption that $\{x_i\}$ is reduced.
\end{proof}

\begin{proof}[Proof of Theorem \ref{thm-freeJinL}]
Let $J \subset L_0 \sdp \freeL L_1$ be a Lie subalgebra such that $J
\cap L_0 = 0$.
By Lemma~\ref{lemma-reduced} $J$ has a reduced set of generators $\{x_i\}$.

Assume $\exists \ 0 \neq f(y_1, \ldots y_n) \in \freeL \langle y_1,
\ldots y_n \rangle$ with $|y_i| = |x_i|$ such that $f(x_1, \ldots x_n)
= 0$.
Without loss of generality assume that $f(x_1, \ldots x_n)$ is a sum
of nonzero monomials $\{f_j(x_1, \ldots x_n)\}_{j \in J}$ in a fixed
dimension.

Since $f_j \in \freeL \langle y_1, \ldots y_n \rangle$, by the
previous lemma $f_j(x\p_1, \ldots x\p_n) \neq 0$.
Each of these is in a fixed degree $M_j$.

By linearity and for degree reasons, $f_j(x_1, \ldots x_n) = \alpha_j
+ f_j(x\p_1, \ldots x\p_n)$, with $\alpha_j \in J_{M_j-1}$.
Let $\{f_k\}_{k \in K}$ be the subset of the monomials $\{f_j\}_{j \in
J}$ such that $f_k(x\p_1, \ldots x\p_n)$ is in the highest degree, say $M$.
Then 
\[
\begin{split}
0 = f(x_1, \ldots x_n) & = \Sigma_{j \in J} ( \alpha_j + f_j(x\p_1,
\ldots x\p_n) ) \\
& = \alpha + \Sigma_{k \in K} f_k(x\p_1, \ldots x\p_n), \quad \alpha
\in J_{M-1}.
\end{split}
\]
Therefore $\Sigma_{k \in K} f_k(x\p_1, \ldots x\p_n) = 0$.
But this contradicts the previous lemma.
\end{proof}

\ignore{


In Section~\ref{section-hurewicz} we constructed Hurewicz images from a
bigraded Lie algebra $L$ over a field $\F$ with generators and
relations in degrees $0$ and $1$. 
However, the Lie algebra of Hurewicz images, $\bar{L}_Y$ is
(potentially) more complicated than $L$.
It is filtered though not necessarily graded by degree and lies in $UL$
and not necessarily in $L$.
However the component of $\bar{L}_Y$ in degree $0$ equals $L_0$.

In the remainder of this section we will prove the following.

\begin{thm} \label{thm-freeJinLY}
Let $\bar{L}_Y$ be the Lie algebra of Hurewicz images over \F of a space $Y$
of the type described in Section~\ref{section-hurewicz}.
Let $J$ be a Lie subalgebra of $\bar{L}_Y$ such that $J \cap (\bar{L}_Y)_0
= 0$, then $J$ is a free Lie algebra.
\end{thm}

Let $\{z_j\}$ be a basis for $L_1$. We will assume without loss of
generality that these elements, and others below, are homogeneous in
dimension.
By Corollary~\ref{cor-L+}, $L_+ \isom \freeL \langle z_j
\rangle$.

For each $z_j$, using Proposition~\ref{prop-lsul} choose $z_{j,h} \in
\bar{L}_Y$ such that $z_j = z_{j,h} + z_{j,0}$, $z_{j,0} \in UL\p_0$.
To simplify the notation write $y_j = z_{j,h}$.
Then $y_j = z_j - z_{j,0}$.
Let $\langle y_j \rangle \subset \bar{L}_Y$ be the Lie subalgebra
generated by $\{y_j\}$.

\begin{lemma} \label{lemma-yfree}
$\langle y_j \rangle \isom \freeL \langle y_j \rangle$ as Lie
algebras.
Furthermore, $\forall 0 \neq x \in \langle y_j \rangle$, $x = \Sigma_{j=0}^m
x^{(j)}$, $x^{(j)} \in (UL)_j$, $m \geq 1$, and $0 \neq x^{(m)} \in L_m$.
\end{lemma}

\begin{proof}
Let $0 \neq f(w_1,\ldots,w_n) \in \freeL \langle w_1, \ldots, w_n
\rangle$, with $|w_i| = |y_i| = |z_i|$.
Since $L_+ \isom \freeL \langle z_j \rangle$, $f(z_1, \ldots, z_n) \neq
0$.
Let $M$ be the degree of $f(z_1,\ldots z_n)$. 
Then by linearity, $f(y_1, \ldots, y_n) = \alpha + f(z_1,\ldots,z_n)$,
where $\alpha \in (UL)_{<M}$.
Thus for degree reasons, $f(y_1,\ldots y_n) \neq 0$.

If $0 \neq x \in \langle y_j \rangle \isom \freeL \langle y_j
\rangle$, then $x = f(y_1, \ldots y_n) = \alpha + f(z_1, \ldots,
z_n)$, where $\alpha \in (UL)_{<m}$ as above.
Then $x^{(m)} = f(z_1,\ldots,z_n) \in L_m$ is nonzero.
\end{proof}

\begin{lemma} \label{lemma-ly}
$\bar{L}_Y \isom (\bar{L}_Y)_0 \oplus \langle y_j \rangle$ as \F\!\!-modules.
\end{lemma}

\begin{proof}
As \F\!\!-modules, $\bar{L}_Y \isom \gr(\bar{L}_Y) \isom L \isom L_0
\oplus \freeL (L_1) \isom L_0 \oplus \freeL \langle z_j \rangle \isom
(\bar{L}_Y)_0 \oplus \freeL \langle y_j \rangle$.
\end{proof}

\begin{prop}
If $x \in \bar{L}_Y$ but $x \notin (\bar{L}_Y)_0$ then 
$x = \Sigma_{j=0}^m x^{(j)}$, $x^{(j)} \in (UL)_j$, $m \geq 1$, and 
$0 \neq x^{(m)} \in L_m$.  
\end{prop}

\begin{proof}
It follows from Lemma~\ref{lemma-ly} that any $x \in \bar{L}_Y$ can be
written as $x=x_a + x_b$, 
where $x_a \in (\bar{L}_Y)_0 \isom L_0$, and $x_b \in \langle y_j \rangle$.
Since $x \notin (\bar{L}_Y)_0$, $x_b \neq 0$. 
By Lemma~\ref{lemma-yfree}, $x_b = \Sigma_{j=0}^m 
x^{(j)}$, $x^{(j)} \in (UL)_j$, $m \geq 1$, and $0 \neq x^{(m)} \in
L_m$.
The statement of the Proposition follows.
\end{proof}

\begin{proof}[Proof of Theorem~\ref{thm-freeJinLY}]
By Lemma~\ref{lemma-reduced}, $J \subset \bar{L}_Y$ has a reduced set
of generators $\{x_i\}$.
By the previous proposition, $x_i = \Sigma_{j=0}^{m_i} x_i^{(j)}$, $x_i^{(j)}
\in (UL)_j$, $m \geq 1$, and  $0 \neq x_i^{(m_i)} \in L_m$.  
Let $x\p_i = x_i^{(m_i)}$.

From this, the proof of Lemma~\ref{lemma-Nielsen} holds as does the proof of
Theorem~\ref{thm-freeJinL}. 
Therefore $J$ is a free Lie algebra.
\end{proof}

}

\chapter{An Algebraic Ganea Construction} \label{chapter-ganea}

In this section we apply our results to give an algebraic version of
Ganea's Fiber-Cofiber construction \cite{ganea:construction,
rutter:ganea} (see Section~\ref{section-ganea}).
All of the isomorphisms in this chapter are algebra or Lie algebra
isomorphisms.

\section{The rational case} \label{section-ganeaQ}

In the case when $R = \Q$ there is an exact correspondence between Lie
algebras and rational spaces (see~\cite{fht:rht} for a reference).
Using this correspondence, Ganea's construction gives a result about
Lie algebras which we explain below. We will see that this result is a
special case of our results.

Let $(L_0,0)$ be a dgL with a Lie ideal $J \subset L_0$ such that $J$
is a free Lie algebra.
Then as explained in~\cite{fht:rht} the short exact sequence of Lie algebras
\[ 0 \to J \to L_0 \to L_0/J \to 0
\] 
models a fibration
\[ F \xto{f} X \to B.
\]
In Ganea's construction, we take the cofiber $X\p = X \cup_f CF$,
where $CF$ is the cone on $F$,
and the fiber $F\p$ of the fibration $X\p \to B$.

Choose a set of generators $S$ for the Lie ideal $J$.
Let $V_0$ be the free $R$-module with basis $S$.
Let $V_1$ be an $R$-module such that there is a bijection $d: V_1 \to
V_0$ where $d$ reduces dimension by $1$. 
Let $L = (L_0 \amalg \freeL V_1, d)$ be the dgL where $dL_0 = 0$ and
$dV_1$ is give above.
Then $L$ models $X\p$~\cite{fht:rht}.
Let $J\p = \ker(L \onto L_0/J)$.
Then $J\p$ models $F\p$~\cite{fht:rht}.

The long exact homology sequence splits into short exact sequences
which can be collected to give the short exact sequence
of Lie algebras  
\[ 0 \to HJ\p \to HL \to L_0/J \to 0.
\]
By Ganea's Theorem \cite{ganea:construction} $F\p \approx F * \Omega
B \approx \Sigma(F \wedge \Omega B)$ which rationally is a wedge of spheres.
So since $J\p$ models $F\p$, $HJ\p$ should be a free Lie algebra.

Using our results and without using Ganea's Theorem,  we will see what
this free Lie algebra is. 
In our terminology $UL$ is a dga extension of $((UL_0,0),L_0)$.
By Lemma~\ref{lemma-dA} and Theorem~\ref{thm-b}, as algebras $HUL
\isom U \psi HL \isom UHL$ (since the map $\psi: UHL \to HUL$ is an
isomorphism when $R = \Q$~\cite[Theorem 21.7(i)]{fht:rht}), where $HL
\isom L_0/J \sdp \freeL (HL)_1$. 
By Lemma~\ref{lemma-sdp} $HL$ fits into the split short exact sequence of
Lie algebras
\[
0 \to \freeL (HL)_1 \to HL \to L_0/J \to 0.
\]
So we see that $HJ\p \isom \freeL (HL)_1$.

\begin{rem}
In~\cite{felixThomas:fcof} F\'{e}lix and Thomas also do this calculation using
Ganea's Theorem.
They show that $HJ\p \isom \freeL K$ where $K$ is the kernel of the
canonical homomorphism $U(L/J) \tensor V_0 \to J/J^2$.
\end{rem}

\section{The construction}

Let $R = \Fp$ with $p>3$ or $R \subset \Q$ containing $\frac{1}{6}$.
For the cases when $R \neq \Q$ we do not have an exact correspondence
between Lie models and spaces. 
However our results allow us to give an algebraic analogue to Ganea's
construction. 

Let $(L_0,0)$ be a dgL with a Lie ideal $J \subset L_0$ such that $J$
is a free Lie algebra.
Choose a set of generators $S$ for the Lie ideal $J$.
Let $V_1$ be an $R$-module such that there is a bijection $d: V_1 \to
R\{S\}$ where $d$ reduces dimension by $1$. 
Let $L = (L_0 \amalg \freeL V_1, d)$ be the dgL where $dL_0 = 0$ and
$dV_1$ is given above.
Then $[dV_1] = J \subset L_0$.
In our terminology $UL$ is a dga extension of $((UL_0,0),L_0)$ which
is free.
By Lemma~\ref{lemma-dA} and Theorems \ref{thm-a} and \ref{thm-b}, as algebras
\[ HUL \isom U((HL)_0 \sdp \freeL (HL)_1)
\]
where $(HL)_0 \isom L_0/J$.
By Lemma~\ref{lemma-sdp}, $HUL \isom UL\p$ as algebras where 
\[ 0 \to \freeL (HL)_1 \to L\p \to L_0/J \to 0
\]
is a split short exact sequence of Lie algebras.
If $R \subset \Q$ then $L\p = \psi HL$ where $\psi$ is the natural map
$\psi: UHL \to HUL$.

As seen in Section~\ref{section-ganeaQ} this generalizes the algebraic
analogue of Ganea's construction obtained using rational homotopy
theory.

\section{Iterating the construction}

Let $R = \Fp$ where $p>3$ or $R \subset \Q$ containing $\frac{1}{6}$.
As with Ganea's construction (see Section~\ref{section-ganea}) this
algebraic construction can be iterated.
This can be used to obtain successively closer approximations to a
given Lie algebra.

Given the split short exact sequence of Lie algebras  
\[ 0 \to \freeL(HL^{(n)})_1 \to {L^{(n)}}\p \to L_0/J \to 0,
\]
choose a set $S^{(n)}$ of generators for the Lie ideal $\freeL
(HL^{(n)})_1$.
Let $V^{(n+1)}$ be an $R$-module such that there is a bijection $d:
V^{(n+1)} \to RS^{(n)}$ where $d$ reduces dimension by $1$. 
Let
\[ (L^{(n+1)}, d^{(n+1)}) = ({L^{(n)}}\p \amalg \freeL V^{(n+1)},
d^{(n+1)}),
\]
where $d^{(n+1)}L^{(n)} = 0$, and $d^{(n+1)} V^{(n+1)} \subset
{L^{(n)}}\p$ is given above. 
Then $[d^{(n+1)} V^{(n+1)}] = \freeL (HL^{(n)})_1 \subset {L^{(n)}}\p$.
By Lemma~\ref{lemma-dA} and Theorems \ref{thm-a} and \ref{thm-b}, as
algebras $HUL^{(n+1)} \isom U{L^{(n+1)}}\p$ where 
\[ 0 \to \freeL (HL^{(n+1)})_1 \to {L^{(n+1)}}\p \to L_0/J \to 0
\]
is a split short exact sequence of Lie algebras.
By construction the lowest dimension in which ${L^{(n)}}\p$ and
$L_0/J$ differ increases to $\infty$ as $n \to \infty$.

\ignore{

Alternatively let $\check{S}^{(n)} \subset ZL^{(n)}$ be a set of
cycles representing $S^{(n)}$.
Let $V^{(n+1)}$ be an $R$-module such that there is a bijection $d:
V^{(n+1)} \to R\check{S}^{(n)}$ which lowers dimension by one.
Let
\[ (L^{(n+1)}, d^{(n+1)}) = ({L^{(n)}} \amalg \freeL V^{(n+1)},
d^{(n+1)}),
\]
where $d^{(n+1)} |_{L^{(n)}} = d^{(n)} |_{L^{(n)}}$ and
$d^{(n+1)} V^{(n+1)} \subset ZL^{(n)}$ is given above.
Then there is an induced dgL $(\hat{L}^{(n+1)}, \hat{d}^{(n+1)}) =
({L^{(n)}}\p \amalg \freeL V^{(n+1)}, \hat{d}^{(n+1)})$ where
$\hat{d}^{(n+1)} {L^{(n)}}\p = 0$ and $\hat{d}^{(n+1)} V^{(n+1)} \subset
(HL^{(n)})_1$ where the differential is the composite map
$\hat{d}^{(n+1)}: V^{(n+1)} \xto{d^{(n+1)}} ZL^{(n)} \onto HL^{(n)}$.
So $\hat{d}^{(n+1)} V^{(n+1)} = S^{(n)}$ and $[\hat{d}^{(n+1)}
V^{(n+1)}] = \freeL (HL^{(n)})_1$. 
Then $(L^{(n+1)}, d^{(n+1)})$ is a dgL extension of $L^{(n)}$ which is
free.
By Theorems \ref{thm-a} and \ref{thm-b}
\[ [xxx] \gr(HU(L^{(n+1)}, d^{(n+1)})) \isom U{L^{(n+1)}}\p
\] 
where ${L^{(n+1)}}\p$ fits into the split short exact sequence of Lie 
algebras
\[
\freeL (H\hat{L}^{(n+1)})_1 \to {L^{(n+1)}}\p \to L_0/J.
\]

Let $(L,d) = \varinjlim (L^{(n)},d^{(n)})$.
Then $\gr(HU(L,d)) \isom U(L_0/J)$.
Since this is in degree $0$ we have that $HU(L,d) \isom U(L_0/J)$ as
algebras. 
If $L_0 \isom \freeL (V^{(0)})$ then  $(L,d) \isom \freeL
(\bigoplus_{k\geq 0} V^{(k)})$. 
Thus $U(L,d)$ is a free model for $U(L_0/J)$.
Since any Lie algebra can be written as $L_0/J$ for some free Lie
algebra $L_0$, this construction produces a free model in which
the homology of each stage in the construction is understood.

}

\section{A more general construction}

Ganea's construction can be generalized to taking more
general cofibers than $X \cup_f CF$ \cite{mather:ganea,
felixThomas:fcof} (see Section~\ref{section-ganea}).
We can do the same in our algebraic construction.
Let $L = (L_0 \amalg \freeL V_1, d)$ where $dL_0 = 0$ and $dV_1
\subset J \subset L_0$. 
Then $[dV_1] \subset J \subset L_0$.
Since $[dV_1]$ is a subalgebra of the free Lie algebra $J$ it is also
a free Lie algebra.
So $UL$ is a free dga extension of $((UL_0,0),L_0)$.
Therefore
\[ HUL \isom U( (HL)_0 \sdp \freeL (HL)_1), \text{ where } (HL)_0
\isom L/[dV_1].
\]

\chapter {Topological Examples and Open Questions} \label{chapter-eg}

In this chapter we use our results to analyze various topological
examples. 
We conclude with some open questions.

\section{Topological examples}

All of the isomorphisms in this section are isomorphisms of algebras
or Lie algebras. 

The following \emph{spherical three-cone} $Y$, illustrates our results.

\begin{eg} \label{eg-semiInert3coneConstruct}
Let $R = \Z[\frac{1}{6}]$. 
Let  $A = S^3 \vee S^3$ and let $Z = A \cup_{\alpha_1 \vee \alpha_2}
(e^8 \vee e^8)$ where the attaching maps are given by the iterated
Whitehead products $\alpha_1 = [[\iota_a,\iota_b],\iota_a]$ and
$\alpha_2 = [[\iota_a,\iota_b],\iota_b]$ of the inclusions of the
$3$-spheres.

It is known (see Example~\ref{eg-wedge-of-spheres}) that $L_A \isom
 \freeL \langle x,y \rangle$ where  $|x|=|y|=2$.
It is also known (see \cite[Example 4.1]{halperinLemaire:attach} or
 Example~\ref{eg-semi-inert-top}) that 
\[
H_*(\Omega Z;R) \isom U(L\p_0 \amalg \freeL \langle \bar{w} \rangle) \text{,
where } L\p_0 = \freeL \langle x,y \rangle/J,
\] 
with $J$ the Lie ideal of all brackets of length $\geq 3$ and $|\bar{w}| = 9$.
By Lemma~\ref{lemma-ip}, $P_Z = \{2,3\}$.
Using \cite{anick:cat2} or Theorem~\ref{thm-hurewiczR}, 
\[
H_*(\Omega Z;R) \isom UL_Z \text{ where } L_Z \isom L^A_Z \amalg
\freeL \langle w \rangle
\]
with $L^A_Z \isom L\p_0$ and $w = h_Z(\hat{\omega})$ where
$\hat{\omega}$ is the adjoint of some map $\omega: S^{10} \to Z$.
By Theorem~\ref{thm-hurewiczR} there is a map $\sigma_Z$ right inverse
to $h_Z$.

For $i=1,2$ let $Z_i$ be two copies of $Z$.
Let $X = Z_1 \vee Z_2$, $W = S^{28} \vee S^{28}$ and let $f = \beta_1
\vee \beta_2$ where $\beta_1 = [[\omega_1, \omega_2], \omega_1]$ and 
$\beta_2 = [[\omega_1, \omega_2], \omega_2]$. 
Let 
\[ Y = X \cup_f (e^{29} \vee e^{29}).
\]

Now, 
\begin{equation} \label{eqn-LXinSIeg} 
L_X \isom L_{Z_1} \amalg L_{Z_2} \isom L^{A_1}_{Z_1} \amalg
L^{A_2}_{Z_2} \amalg \freeL \langle w_1, w_2 \rangle.
\end{equation}
Therefore by Theorem~\ref{thm-schreier2}, $f$ is a free attaching map.
Thus $Y$ satisfies the hypotheses of Theorem~\ref{thm-d}.

By Corollary~\ref{cor-semiInertEqnTop}, if $f$ is semi-inert then 
\[
K\p (z) = \tilde{H}_*(W)(z) + z[UL_X(z)^{-1} - U(L^X_Y)(z)^{-1}].
\]
$L_X$ is given by~\eqref{eqn-LXinSIeg} and $L^X_Y \isom L_X /
[h_X(\hat{\beta_1}), h_X(\hat{\beta_2})]$. 
Since the Hurewicz images are $[[w_1,w_2],w_1]$ and $[[w_1,w_2],w_1]$
which are contained in $\freeL \langle w_1, w_2 \rangle$ we have
\[ 
L^X_Y \isom L^{A_1}_{Z_1} \amalg L^{A_2}_{Z_2} \amalg \tilde{L}
\text{ where } \tilde{L} = \freeL \langle w_1, w_2 \rangle / [
R \{ [[w_1,w_2],w_1],[[w_1,w_2],w_2] \} ].
\]

Since $(A \amalg B)(z)^{-1} = A(z)^{-1} + B(z)^{-1} -
1$~\cite[Lemma 5.1.10]{lemaire:monograph} it follows that
\[
K\p(z) = \tilde{H}_*(W)(z) + z[U\freeL \langle w_1,w_2 \rangle -
U\tilde{L}(z)^{-1}].
\]
Now $U\freeL \langle w_1,w_2 \rangle \isom \freeT \langle w_1,w_2
\rangle$ and as $R$-modules $\tilde{L} \isom R\{
w_1,w_2,[w_1,w_1],[w_1,w_2],[w_2,w_2]\}$ and $U\tilde{L} \isom \mathbb{S}
\tilde{L}.$
Therefore
\begin{eqnarray*}
K\p(z) & = & 2z^{28} + z \left[ 1-2z^9 -
\frac{(1-z^{18})^3}{(1+z^9)^2} \right] \\ 
& = & z^{37}.
\end{eqnarray*}

Recall that $\eL = (L_X \amalg \freeL \langle e,g \rangle, d\p)$ where
$d\p e = [[w_1,w_2],w_1]$ and $d\p g = [[w_1,w_2],w_2]$.
Also recall (from Theorem~\ref{thm-c}) that $(H\eL)_0 \isom L^X_Y$.
Let $\bar{u} = [e,w_2] + [g,w_1]$ (with $|\bar{u}| = 37$).
Then it is easy to check that
\[ (H\eL)_0 \sdp \freeL (H\eL)_1 \isom (H\eL)_0 \amalg \freeL \langle
\bar{u} \rangle
\]
and $f$ is indeed semi-inert.
Since $X = Z_1 \vee Z_2$ we have the following commutative diagram.
\[
\xymatrix{
\pi_*(\lX) \tensor R \ar[r]^-{\isom} \ar[d]^{h_X} &
\pi_*(\lZ_1) \tensor R \oplus \pi_*(\lZ_2) \tensor R \ar[d]^{h_{Z_1}
\oplus h_{Z_2}} \\
H_*(\lX;R) \ar[r]^-{\isom} & H_*(\lZ_1;R) \oplus H_*(\lZ_2;R)
}
\]
Let $\sigma_X$ be the map corresponding to $\sigma_{Z_1} \oplus
\sigma_{Z_2}$ under these isomorphisms.
Since $\sigma_{Z_1}$ and $\sigma_{Z_2}$ are right inverses to
$h_{Z_1}$ and $h_{Z_2}$ it follows that $\sigma_X = \sigma_{Z_1} \oplus
\sigma_{Z_2}$ is right inverse to $h_X$.
By Lemma~\ref{lemma-ip}, $P_Y = \{2,3\}$.
As a result by Theorem~\ref{thm-hurewiczR},
\[ H_*(\lY;R) \isom UL_Y \text{ where } L_Y \isom L^X_Y \amalg \freeL
\langle u \rangle
\]
with $u = h_Y(\hat{\mu})$ for some map $\mu: S^{38} \to Y$.
Furthermore $\lY \in \PiS$ and there exists a map $\sigma_Y$ right
inverse to $h_Y$.  

Note that
\[
L_Y \isom L_1^1 \amalg L_1^2 \amalg L_2 \amalg \freeL \langle u
\rangle
\]
where $L_1^1 \isom L_1^2 \isom \freeL \langle x,y \rangle / J_1$ and
$L_2 \isom \freeL \langle w_1,w_2 \rangle / J_2$
with $J_1$ and $J_2$ the Lie ideals of all brackets of length $\geq
3$. 
\eolBox
\end{eg}

\begin{eg} \label{eg-semiInertNCones}
An infinite family of finite CW-complexes constructed out of
semi-inert attaching maps
\end{eg}

The construction in the previous example can be extended inductively. 
By induction, we will construct spaces $X_n$ and maps $\omega_n:
S^{\lambda_n} \to X_n$ for $n \geq 1$ such that $X_n$ is an $n$-cone
constructed out of a sequence of semi-inert attaching maps.
Given $\omega_n$, let $w_n = h_{X_n}([\omega_n])$ and given $w_i^a$ and
$w_i^b$, 
let $L_i = \freeL \langle w_i^a, w_i^b \rangle/ J_i$ where $J_i$ is
the Lie ideal generated by $[[w_i^a,w_i^b],w_i^a]$ and
$[[w_i^a,w_i^b],w_i^b]$. 
That is $L_i$ is the quotient of $\freeL \langle w_i^a, w_i^b \rangle$
where all brackets of length three are equal to zero.

Let $R = \Z[\frac{1}{6}]$.
Begin with $X_1 = S^3$ and $\lambda_1 = 3$.  
Let $\omega_1: S^{\lambda_1} \to X_1$ be the identity map.

Given $X_n$, let $X_n^a$ and $X_n^b$ be two copies of $X_n$.
For $n\geq 1$, let 
\[
X_{n+1} = X_n^a \vee X_n^b \cup_{f_{n+1}} \left(
e^{\kappa_{n+1}} \vee e^{\kappa_{n+1}} \right),
\]
where $\kappa_{n+1} = 3 \lambda_n - 1$ and $f_{n+1} = [[\omega_n^a,
\omega_n^b], \omega_n^a] \vee [[\omega_n^a, \omega_n^b], \omega_n^b]$.

By the same argument as in the previous example, 
$f_{n+1}$ is a semi-inert cell attachment and
there exists a map 
\[
\omega_{n+1}: S^{\lambda_{n+1}} \to X_{n+1}
\] 
where $\lambda_{n+1} = 4 \lambda_n - 2$, such that 
\[
H_*(\lX_{n+1}; R) \isom UL_{X_{n+1}} \text{ where } L_{X_{n+1}} =
\Bigl( \coprod_{ \substack{1 \leq i \leq n \\ 1 \leq j \leq 2^{n-i}} }
L^j_i \Bigr) 
\amalg \freeL \langle w_{n+1} \rangle
\]
with $L^j_i$ a copy of $L_i$ and $w_{n+1} = h_{X_{n+1}} ([\omega_{n+1}])$.
\eolBox

\begin{eg} \label{eg-2cones}
Spherical $2$-Cones
\end{eg}

A \emph{spherical $2$-cone} is an \emph{adjunction space} (see
Section~\ref{section-basic-top}) of the form
\[ \Bigl( \bigvee_{j \in J_1} S^{m_j} \Bigr) \cup_f \Bigl( \bigvee_{j \in J_2}
e^{n_j+1} \Bigr).
\]
Equivalently it is the adjunction space of a cell attachment $W_2
\xto{f_2} X_1$, where $W_2= \bigvee_{j \in J_2} S^{n_j}$ and $X_1$ are
finite-type wedges of spheres and $f_2 = \bigvee_{j \in J_2} \alpha_j^{(2)}$.
This is the situation studied by Anick in~\cite{anick:cat2}.
Let $R = \F$ where $\F = \Q$ or $\Fp$ with $p>3$,
or let $R \subset \Q$ be a subring containing $\frac{1}{6}$.
By the Hilton-Milnor Theorem $H_*(\lX_1;R) \isom UL_{X_1}$ and
$L_{X_1} \isom \freeL V$ for some free $R$-module $V$.
$X_1$ has Adams-Hilton model $U(\freeL V, 0)$.
Since $L_{X_1} \isom \freeL V$ we can define a map $\sigma_{X_1}:
L_{X_1} \to \pi_*(\lX_1) \tensor R$ right inverse to $h_{X_1}$ by
choosing inverse images of $V$ and extending to a Lie algebra map.
Let $\eL_2 = (L_{X_1} \amalg \freeL \langle y_j^{(2)} \rangle_{j \in
J_2}, d\p)$ where $d\p y_j^{(2)} = h_{X_1}(\hat{\alpha}_j^{(2)})$.

If $R = \F$, $[L^{W_2}_{X_1}] \subset L_{X_1}$ is a Lie subalgebra of
a free Lie algebra and so is automatically a free Lie algebra. 
Since the Adams-Hilton model for $X_1$ has zero differential, by
Lemma~\ref{lemma-dA} $H_*(\lX_2;\F) \isom \gr(H_*(\lX_2;\F))$.
Applying Theorem~\ref{thm-c}(i) we get the following isomorphism of
algebras
\[ H_*(\lX_2;\F) \isom U \left( {L}^{X_1}_{X_2} \sdp \freeL (H\eL_2)_1
\right).
\]
Furthermore by Theorem~\ref{thm-hurewiczF} if $\sigma_{X_1} h_{X_1}
\hat{\alpha}_j^{(2)} = \hat{\alpha}_j^{(2)}$ then the canonical map
\[ UL_{X_2} \to H_*(\lX_2;\F)
\] 
is a surjection.

If $R \subset \Q$ and $L_{X_1}/[L^{W_2}_{X_1}]$ is $R$-free then we
can apply Theorem~\ref{thm-d}(i) to get that $H_*(\lX_2;R)$ is
torsion-free and 
\[ H_*(\lX_2;R) \isom U \left(L^{X_1}_{X_2} \sdp \freeL (H\eL_2)_1 \right).
\]
Furthermore since $\exists \ \sigma_{X_1}$ right inverse to $h_{X_1}$,
then by Theorem~\ref{thm-hurewiczR}, 
\[ H_*(\lX_2;\Z[{P_{X_2}}^{-1}]) \isom UL_{X_2}
\]
where $\gr(L_{X_2}) \isom L^{X_1}_{X_2} \sdp \freeL (H\eL_2)_1$ and localized
away from $P_{X_2}$, $\lX_2 \in \PiS$.
This result is given in \cite{anick:cat2}.

In either case if $f_2$ is semi-inert then by part (ii) of Theorems
\ref{thm-c} and \ref{thm-d}
\[ H_*(\lX_2;R) \isom U \left(L^X_Y \amalg \freeL K^{(2)} \right)
\]
for some $K^{(2)} \subset F_1 H_*(\lX_2;R)$.
If $R \subset \Q$ then by Theorem~\ref{thm-hurewiczR}(iii)
$L_{X_2} \isom L^X_Y \amalg \freeL (H\eL_2)_1$
and there exists a map $\sigma_{X_2}$ right inverse to $h_{X_2}$.

This last Lie algebra isomorphism proves a conjecture of
Anick~\cite[Conjecture 4.9]{anick:cat2} in the semi-inert case.
\eolBox

\begin{eg}
Spherical $3$-Cones
\end{eg}

Consider the adjunction space $X_3$ of a cell attachment $W_3
\xrightarrow{f_3} X_2$ where $X_2$ is the spherical $2$-cone given above, 
$W_3= \bigvee_{j \in J_3} S^{n_j}$ is a finite-type wedge of
spheres, and $f_3 = \bigvee_{j \in J_3} \alpha_j^{(3)}$.
Let $R$ be a subring of $\Q$ containing $\frac{1}{6}$.
Assume that $H_*(\lX_2;R) \isom UL_{X_2}$ and that there exists a map
$\sigma_{X_2}$ right inverse to $h_{X_2}$.
Let $\eL_3 = (L_{X_2} \amalg \freeL \langle y_j^{(3)} \rangle_{j \in
J_3}, d\p)$ where $d\p y_j^{(3)} = h_{X_2}(\hat{\alpha}_j^{(3)})$.


Assume $L_{X_2}/[L^{W_3}_{X_2}]$ is $R$-free and
\fanypnP, $[\bar{L}^{W_3}_{X_2}] \cap \bar{L}^{X_1}_{X_2} = 0$.
By Example~\ref{eg-2cones}, $\gr(L_{X_2}) \isom L^{X_1}_{X_2} \sdp
\freeL V^{(2)}$ for some free $R$-module $V^{(2)}$.
It follows that \fanypnP,  $\gr(\bar{L}_{X_2}) \isom \bar{L}^{X_1}_{X_2}
\sdp \freeL \bar{V}^{(2)}$.
Then by Theorem~\ref{thm-schreier2}, \fanypnP, $[\bar{L}^{W_3}_{X_2}]$
is a free Lie algebra.
That is, $f_3$ is a free cell-attachment.
Thus by Theorem~\ref{thm-d}(i) we get that $H_*(\lX_3;R)$ is
torsion-free and 
\[ \gr(H_*(\lX_3;R)) \isom U \bigl( L^{X_2}_{X_3} \sdp \freeL
V^{(3)}_1 \bigr).
\]
Furthermore since $\exists \ \sigma_{X_2}$ right inverse to $h_{X_2}$ then by
Theorem~\ref{thm-hurewiczR} 
\[ H_*(\lX_3;\Z[{P_{X_3}}^{-1}]) \isom UL_{X_3}
\]
where $\gr(L_{X_3}) \isom L^X_Y \sdp \freeL (H\eL_3)_1$ and localized
away from $P_{X_3}$, $\lX_3 \in \PiS$.

If in addition $f_3$ is semi-inert then by part (ii) of
Theorem \ref{thm-d},
\[ H_*(\lX_3;R) \isom U \left(L^{X_2}_{X_3} \amalg \freeL K^{(3)} \right),
\]
for some $K^{(3)} \subset F_1 H_*(\lX_3;R)$, and by
Theorem~\ref{thm-hurewiczR}(iii), localized away from $P_{X_3}$ there
exists a map $\sigma_{X_3}$ right inverse to $h_{X_3}$.
\eolBox

\begin{eg} 
Spherical $N$-cones 
\end{eg}

Let $W_k \xrightarrow{f_k} X_{k-1} \to X_k$ for $1 \leq k \leq N$ be a
sequence of cell attachments (with adjunction space $X_k$) with $X_0 =
*$, $W_k$ is a finite-type wedge of spheres and $X_N \approx X$.
Assume $R \subset \Q$ containing $\frac{1}{6}$.
$\forall k$ assume that $L_{X_{k-1}}/[L^{W_k}_{X_{k-1}}]$ is $R$-free
and that \fanypnP, $[\bar{L}^{W_k}_{X_{k-1}}] \cap
\bar{L}^{X_{k-2}}_{X_{k-1}} = 0$. 
In addition $\forall k$ assume that if $f_k$ is free then there 
exists a map $\sigma_{X_{k-1}}$ right inverse to $h_{X_{k-1}}$.
By Theorem~\ref{thm-hurewiczR} a sufficient condition for this is that
$f_{k-1}$ is semi-inert.

Then by induction on Theorem \ref{thm-hurewiczR} and
Theorem~\ref{thm-schreier2}, $H_*(\lX;\Z[{P_X}^{-1}])$ is
torsion-free, $H_*(\lX;\Z[{P_X}^{-1}]) \isom UL_{X}$ as algebras and
localized away from $P_{X}$, $\lX \in \PiS$. 
\eolBox

\begin{eg} 
$W \to X \to Y$ where $X$ is \emph{coformal}.
\end{eg}

Let $W$ be a simply-connected finite-type wedge of spheres.
Let $X$ be a simply-connected topological space such that $H_*(\lX;R)$
is torsion-free and has finite type.
Let $R \subset \Q$ containing $\frac{1}{6}$ or $R = \Fp$ where $p>3$.
Call $X$ \emph{coformal} if $X$ has a \emph{model} (see
Section~\ref{section-ah}) $U\LL{X} =
U(L,0)$~\cite{neisendorferMiller:formal}.
For example, when $R=\Q$, Neisendorfer and
Miller~\cite{neisendorferMiller:formal} show that any compact
$n$-connected $m$-dimensional manifold with $m \leq 3n+1$, $n\geq 1$
and $\rank (PH(M;\Q)) \geq 2$ is coformal, where $P(\cdot)$ denotes
the subspace of primitive elements. 

Let $Y$ be the adjunction space of a cell attachment $W \to X$.
For a model for $Y$ one can take $U\LL{Y} = U(L \amalg \freeL V_1,
d)$, where $dL=0$ and $dV_1 \subset L$. 
Assume that $[dV_1] \subset L$ is a free Lie algebra.
Since the differential in $\LL{X}$ is zero, by
Lemma~\ref{lemma-dA}, $H_*(\lY;R) \isom \gr(H_*(\lY;R))$.

By Theorems \ref{thm-c}(i) and \ref{thm-d}(i) if $L/[dV_1]$ is torsion-free
then
\[
H_*(\lY;R) \isom U((H\LL{Y})_0 \sdp \freeL (H\LL{Y})_1) \text{ where }
(H\LL{Y})_0 \isom L/[dV_1].
\]
\eolBox

\begin{eg} $X$ is a finite-type CW-complex with only odd dimensional
cells.
\end{eg}

Let $X^{(n)}$ denote the $n$-skeleton of $X$.
From the CW-structure of $X$, there is a sequence of cell attachments for
$k\geq 1$. 
\[ W_{2k} \xto{f_{2k+1}} X^{(2k-1)} \to X^{(2k+1)}
\]
where $W_{2k}$ is a finite wedge of $(2k)$-dimensional spheres,
$X^{(2k+1)}$ is the adjunction space of the cell attachment and
$X^{(1)} = *$.  
Assume $R \subset \Q$ containing $\frac{1}{6}$.
Let $P_{X^{(n)}}$ be the set of implicit primes of $X^{(n)}$ (see
Section~\ref{section-ip}). 
We will show by induction that
\[
H_*(\lX^{(2k+1)}; \Z[{P_{X^{(2k+1)}}}^{-1}]) \isom UL_{X^{(2k+1)}}
\]
where $L_{X^{(2k+1)}} \isom \freeL V^{(2k+1)}$ with 
$V^{(2k+1)}$ concentrated in even dimensions.
In addition localized away from $P_{X^{(2k+1)}}$, $\lX^{(2k+1)} \in \PiS$
and there exists a map $\sigma_{X^{(2k+1)}}$ right inverse to
$h_{X^{(2k+1)}}$.

For $k=0$ these conditions are trivial. 
Assume they hold for $k-1$.
Let $\eL = (L_{X^{(2k-1)}} \amalg \freeL K, d\p)$ where $K$ is a free
$R$-module in dimension $2k$ corresponding to the spheres in $W_{2k}$.
For degree reasons $L^{W_{2k}}_{X^{(2k-1)}} = d\p(K) = 0$.
So $f_{2k+1}$ is automatically free.
Furthermore $\eL$ has zero differential so $H\eL = \freeL V^{(2k-1)}
\amalg \freeL K$.
Thus $f_{2k+1}$ is semi-inert. 
By Theorem~\ref{thm-hurewiczR},
$H_*(\lX^{(2k+1)}; \Z[{P_{X^{(2k+1)}}}^{-1}]) \isom UL_{X^{(2k+1)}}$, 
where 
\[
L_{X^{(2k+1)}} \isom L^{X^{(2k-1)}}_{X^{(2k+1)}} \amalg \freeL K
\isom L_{X^{(2k-1)}} \amalg \freeL K \isom \freeL (V^{(2k-1)} \oplus
K).
\]
Also by Theorem~\ref{thm-hurewiczR}, localized away from
$P_{X^{(2k+1)}}$, $\lX^{(2k+1)} \in \PiS$, and there exists a map
$\sigma_{X^{(2k+1)}}$ right inverse to  $h_{X^{(2k+1)}}$.

Therefore by induction $H_*(\lX;\Z[{P_X}^{-1}]) \isom UL_X$, where $L_X
\isom \freeL V$ and localized away from $P_X$, $\lX \in \PiS$. 
\eolBox

\section{Open questions}

The results of this thesis naturally lead to the following open
questions.
\begin{itemize}
\item 
Is there a dga extension $\bA$ which is free such that
$\gr(H\bA) \not\cong H\bA$ as algebras? 
\item 
If $H(\check{A},\check{d})$ is $R$-free and as algebras
$H(\check{A},\check{d}) \isom UL_0$ for some Lie algebra $L_0$ then is
there a dga extension $(\check{A} \amalg \freeT V, d)$ for which
one cannot choose a Lie algebra $L_0$ $d\p V \subset L_0$?
\item 
If so, do any of these come from cell attachments, or have
$(\check{A},\check{d}) = U(\check{L},\check{d})$?
\item 
Let $R = \Q$ or $\Fp$ where $p>3$ and let $(\hat{A},\hat{d})$ be
a free dga extension of  $((\check{A},\check{d}),L_0)$. 
By Theorem~\ref{thm-a}, as algebras $H(\hat{A},\hat{d}) \isom UL\p$
where $L\p = (H\eL)_0 \sdp \freeL (H\eL)_1$.  
Let $\bA = \hat{A} \amalg \freeT V$ be a dga extension of $\hat{A}$
with induced map $d\p: V \to L\p$ such that $[d\p V] \cap (H\eL)_0 = 0$.
By Theorem~\ref{thm-schreier2}, $\bA$ is a free dgL extension.
Do the relations which are obstructions to the semi-inertness of
$\bA$ necessarily come from the anti-commutativity and Jacobi relations
(as in Examples \ref{eg-semi-inert-extn} and \ref{eg-alg-fat-wedge})?
\item 
Is there a finite CW-complex $X$ such that $H_*(\lX;R)$ is
$R$-free and $H_*(\lX;R) \isom UL_X$ as algebras but there does not
exist a map $\sigma_X$ right inverse to $h_X$ even after localizing
away from finitely many primes?
\item
One of the forms of the semi-inert condition, namely that
$\gr_1(H\bA)$ is a free $\gr_0(H\bA)$-bimodule, can be applied to more
general differential graded algebras than those studied in this
thesis.
J.-M. Lemaire has suggested that the condition that the two-sided ideal
\begin{equation} \label{eqn-gen-free-condn}
(dV_1) \subset \check{A} \text{ is a free } \check{A} \text{-module}
\end{equation}
be taken as the free condition for more general dga's.
(For the dga's in this thesis, \eqref{eqn-gen-free-condn} follows from
the free condition by Lemma~\ref{lemma-I}.
Can Theorems \ref{thm-a} and \ref{thm-b} be generalized by using these
conditions? 
\end{itemize}

\addcontentsline{toc}{part}{Bibliography}
\bibliographystyle{alpha}
\bibliography{my}

\begin{thebibliography}{CMN79b}

\bibitem[AH56]{adamsHilton}
J.~F. Adams and P.~J. Hilton.
\newblock On the chain algebra of a loop space.
\newblock {\em Comment. Math. Helv.}, 30:305--330, 1956.

\bibitem[Ani82a]{anick:thesis}
David~J. Anick.
\newblock A counterexample to a conjecture of {S}erre.
\newblock {\em Ann. of Math. (2)}, 115(1):1--33, 1982.

\bibitem[Ani82b]{anick:stronglyFree}
David~J. Anick.
\newblock Noncommutative graded algebras and their {H}ilbert series.
\newblock {\em J. Algebra}, 78(1):120--140, 1982.

\bibitem[Ani86]{anick:torsion}
David~J. Anick.
\newblock A loop space whose homology has torsion of all orders.
\newblock {\em Pacific J. Math.}, 123(2):257--262, 1986.

\bibitem[Ani89]{anick:cat2}
David~J. Anick.
\newblock Homotopy exponents for spaces of category two.
\newblock In {\em Algebraic topology (Arcata, CA, 1986)}, volume 1370 of {\em
  Lecture Notes in Math.}, pages 24--52. Springer, Berlin, 1989.

\bibitem[Ani92]{anick:slsd}
David~J. Anick.
\newblock Single loop space decompositions.
\newblock {\em Trans. Amer. Math. Soc.}, 334(2):929--940, 1992.

\bibitem[Avr86]{avramov:torsion}
Luchezar~L. Avramov.
\newblock Torsion in loop space homology.
\newblock {\em Topology}, 25(2):155--157, 1986.

\bibitem[Bau81]{baues:commutator}
Hans~Joachim Baues.
\newblock {\em Commutator calculus and groups of homotopy classes}, volume~50
  of {\em London Mathematical Society Lecture Note Series}.
\newblock Cambridge University Press, Cambridge, 1981.

\bibitem[BFT92]{bft:gammaFunctor}
Hans~Joachim Baues, Yves F{\'e}lix, and Jean-Claude Thomas.
\newblock Presentations of algebras and the {W}hitehead {$\Gamma$}-functor for
  chain algebras.
\newblock {\em J. Algebra}, 148(1):123--134, 1992.

\bibitem[BK72]{bousfieldKan:book}
A.~K. Bousfield and D.~M. Kan.
\newblock {\em Homotopy limits, completions and localizations}.
\newblock Springer-Verlag, Berlin, 1972.
\newblock Lecture Notes in Mathematics, Vol. 304.

\bibitem[B{\o}g85]{bogvad:envelopingAlg}
Rikard B{\o}gvad.
\newblock The enveloping algebra of a graded {L}ie algebra of global dimension
  two contains a free subalgebra on two generators.
\newblock {\em J. Pure Appl. Algebra}, 38(2-3):213--216, 1985.

\bibitem[BS53]{bottSamelson:thm}
Raoul Bott and Hans Samelson.
\newblock On the {P}ontryagin product in spaces of paths.
\newblock {\em Comment. Math. Helv.}, 27:320--337 (1954), 1953.

\bibitem[CMN79a]{cmn:aarhus}
F.~R. Cohen, J.~C. Moore, and J.~A. Neisendorfer.
\newblock Decompositions of loop spaces and applications to exponents.
\newblock In {\em Algebraic topology, Aarhus 1978 (Proc. Sympos., Univ. Aarhus,
  Aarhus, 1978)}, volume 763 of {\em Lecture Notes in Math.}, pages 1--12.
  Springer, Berlin, 1979.

\bibitem[CMN79b]{cmn:annals110}
F.~R. Cohen, J.~C. Moore, and J.~A. Neisendorfer.
\newblock The double suspension and exponents of the homotopy groups of
  spheres.
\newblock {\em Ann. of Math. (2)}, 110(3):549--565, 1979.

\bibitem[CMN79c]{cmn:annals109}
F.~R. Cohen, J.~C. Moore, and J.~A. Neisendorfer.
\newblock Torsion in homotopy groups.
\newblock {\em Ann. of Math. (2)}, 109(1):121--168, 1979.

\bibitem[Coh95]{cohen:decompositionsInHandbook}
Frederick~R. Cohen.
\newblock Fibration and product decompositions in nonstable homotopy theory.
\newblock In {\em Handbook of algebraic topology}, pages 1175--1208.
  North-Holland, Amsterdam, 1995.

\bibitem[F{\'e}l89]{felix:dichotomie}
Yves F{\'e}lix.
\newblock La dichotomie elliptique-hyperbolique en homotopie rationnelle.
\newblock {\em Ast\'erisque}, (176):187, 1989.

\bibitem[F{\'e}l92]{felix:H*LoopX}
Yves F{\'e}lix.
\newblock Structure de l'alg\`ebre de {H}opf {$H\sb *(\Omega X;k)$}.
\newblock {\em Bull. Soc. Math. Belg. S\'er. A}, 44(1):35--50, 1992.

\bibitem[FHT82]{fht:homotopyLieAlg}
Yves F{\'e}lix, Stephen Halperin, and Jean-Claude Thomas.
\newblock The homotopy {L}ie algebra for finite complexes.
\newblock {\em Inst. Hautes \'Etudes Sci. Publ. Math.}, (56):179--202 (1983),
  1982.

\bibitem[FHT84]{fht:cat2}
Yves F{\'e}lix, Stephen Halperin, and Jean-Claude Thomas.
\newblock Sur l'homotopie des espaces de cat\'egorie {$2$}.
\newblock {\em Math. Scand.}, 55(2):216--228, 1984.

\bibitem[FHT90]{fht:lshCat12}
Yves F{\'e}lix, Stephen Halperin, and Jean-Claude Thomas.
\newblock Loop space homology of spaces of {LS} category one and two.
\newblock {\em Math. Ann.}, 287(3):377--386, 1990.

\bibitem[FHT01]{fht:rht}
Yves F{\'e}lix, Stephen Halperin, and Jean-Claude Thomas.
\newblock {\em Rational homotopy theory}, volume 205 of {\em Graduate Texts in
  Mathematics}.
\newblock Springer-Verlag, New York, 2001.

\bibitem[FL91]{felixLemaire:pontryagin}
Yves F{\'e}lix and Jean-Michael Lemaire.
\newblock On the {P}ontryagin algebra of the loops on a space with a cell
  attached.
\newblock {\em Internat. J. Math.}, 2(4):429--438, 1991.

\bibitem[FL92]{felixLemaire:2level}
Yves F{\'e}lix and Jean-Michael Lemaire.
\newblock On the homology of two-level differential algebras.
\newblock {\em Bull. Soc. Math. Belg. S\'er. B}, 44(1):29--33, 1992.

\bibitem[FT86]{felixThomas:cat2}
Yves F{\'e}lix and Jean-Claude Thomas.
\newblock Sur la structure des espaces de {LS} cat\'egorie deux.
\newblock {\em Illinois J. Math.}, 30(4):574--593, 1986.

\bibitem[FT88]{felixThomas:fcof}
Yves F{\'e}lix and Jean-Claude Thomas.
\newblock The fibre-cofibre construction and its applications.
\newblock {\em J. Pure Appl. Algebra}, 53(1-2):59--69, 1988.

\bibitem[FT89]{felixThomas:attach}
Yves F{\'e}lix and Jean-Claude Thomas.
\newblock Effet d'un attachement cellulaire dans l'homologie de l'espace des
  lacets.
\newblock {\em Ann. Inst. Fourier (Grenoble)}, 39(1):207--224, 1989.

\bibitem[FT93]{felixThomas:lshSmallCat}
Yves F{\'e}lix and Jean-Claude Thomas.
\newblock Loop space homology of spaces of small category.
\newblock {\em Trans. Amer. Math. Soc.}, 338(2):711--721, 1993.

\bibitem[Gan65]{ganea:construction}
Tudor Ganea.
\newblock A generalization of the homology and homotopy suspension.
\newblock {\em Comment. Math. Helv.}, 39:295--322, 1965.

\bibitem[Gra69]{gray:sphereOfOrigin}
Brayton Gray.
\newblock On the sphere of origin of infinite families in the homotopy groups
  of spheres.
\newblock {\em Topology}, 8:219--232, 1969.

\bibitem[Hal92]{halperin:uea}
Stephen Halperin.
\newblock Universal enveloping algebras and loop space homology.
\newblock {\em J. Pure Appl. Algebra}, 83(3):237--282, 1992.

\bibitem[HL87]{halperinLemaire:inert}
Stephen Halperin and Jean-Michel Lemaire.
\newblock Suites inertes dans les alg\`ebres de {L}ie gradu\'ees (``{A}utopsie
  d'un meurtre. {II}'').
\newblock {\em Math. Scand.}, 61(1):39--67, 1987.

\bibitem[HL95]{halperinLemaire:attach}
Stephen Halperin and Jean-Michel Lemaire.
\newblock The fibre of a cell attachment.
\newblock {\em Proc. Edinburgh Math. Soc. (2)}, 38(2):295--311, 1995.

\bibitem[HL96]{hessLemaire:nice}
Kathryn Hess and Jean-Michel Lemaire.
\newblock Nice and lazy cell attachments.
\newblock {\em J. Pure Appl. Algebra}, 112(1):29--39, 1996.

\bibitem[HS97]{hiltonStammbach:book}
P.~J. Hilton and U.~Stammbach.
\newblock {\em A course in homological algebra}, volume~4 of {\em Graduate
  Texts in Mathematics}.
\newblock Springer-Verlag, New York, second edition, 1997.

\bibitem[Hus80]{husemoller:cmns}
Dale Husemoller.
\newblock Splittings of loop spaces, torsion in homotopy, and double suspension
  (after {F}. {C}ohen, {J}. {C}. {M}oore, {J}. {N}eisendorfer and {P}. {S}elick
  [{P}aul {S}elick]).
\newblock In {\em Topology Symposium, Siegen 1979 (Proc. Sympos., Univ. Siegen,
  Siegen, 1979)}, volume 788 of {\em Lecture Notes in Math.}, pages 456--470.
  Springer, Berlin, 1980.

\bibitem[Jac79]{jacobson:lieAlgebras}
Nathan Jacobson.
\newblock {\em Lie algebras}.
\newblock Dover Publications Inc., New York, 1979.
\newblock Republication of the 1962 original.

\bibitem[Jam95]{james:lscatInHandbook}
Ioan~M. James.
\newblock Lusternik-{S}chnirelmann category.
\newblock In {\em Handbook of algebraic topology}, pages 1293--1310.
  North-Holland, Amsterdam, 1995.

\bibitem[Lem74]{lemaire:monograph}
Jean-Michel Lemaire.
\newblock {\em Alg\`ebres connexes et homologie des espaces de lacets}.
\newblock Springer-Verlag, Berlin, 1974.
\newblock Lecture Notes in Mathematics, Vol. 422.

\bibitem[Lem78]{lemaire:autopsie}
Jean-Michel Lemaire.
\newblock ``{A}utopsie d'un meurtre'' dans l'homologie d'une alg\`ebre de
  cha\^\i nes.
\newblock {\em Ann. Sci. \'Ecole Norm. Sup. (4)}, 11(1):93--100, 1978.

\bibitem[Mat75]{mather:ganea}
Michael Mather.
\newblock A generalisation of {G}anea's theorem on the mapping cone of the
  inclusion of a fibre.
\newblock {\em J. London Math. Soc. (2)}, 11(1):121--122, 1975.

\bibitem[May99]{may:bookCCAlgTop}
J.~Peter May.
\newblock {\em A concise course in algebraic topology}.
\newblock Chicago Lectures in Mathematics. University of Chicago Press,
  Chicago, IL, 1999.

\bibitem[McC01]{mccleary:usersGuide}
John McCleary.
\newblock {\em A user's guide to spectral sequences}, volume~58 of {\em
  Cambridge Studies in Advanced Mathematics}.
\newblock Cambridge University Press, Cambridge, second edition, 2001.

\bibitem[Mik85]{mikhalev:subalgebras}
Alexander~A. Mikhal{\"e}v.
\newblock Subalgebras of free colored {L}ie superalgebras.
\newblock {\em Mat. Zametki}, 37(5):653--661, 779, 1985.

\bibitem[MM65]{milnorMoore:hopfAlgebras}
John~W. Milnor and John~C. Moore.
\newblock On the structure of {H}opf algebras.
\newblock {\em Ann. of Math. (2)}, 81:211--264, 1965.

\bibitem[MW86]{mcgibbonWilkerson}
Charles~A. McGibbon and Clarence~W. Wilkerson.
\newblock Loop spaces of finite complexes at large primes.
\newblock {\em Proc. Amer. Math. Soc.}, 96(4):698--702, 1986.

\bibitem[MZ95]{mikhalevZolotykh}
Alexander~A. Mikhalev and Andrej~A. Zolotykh.
\newblock {\em Combinatorial aspects of {L}ie superalgebras}.
\newblock CRC Press, Boca Raton, FL, 1995.
\newblock With 1 IBM-PC floppy disk (3.5 inch; HD).

\bibitem[NM78]{neisendorferMiller:formal}
Joseph Neisendorfer and Timothy Miller.
\newblock Formal and coformal spaces.
\newblock {\em Illinois J. Math.}, 22(4):565--580, 1978.

\bibitem[NS82]{neisendorferSelick:someEgs}
Joseph Neisendorfer and Paul Selick.
\newblock Some examples of spaces with or without exponents.
\newblock In {\em Current trends in algebraic topology, Part 1 (London, Ont.,
  1981)}, volume~2 of {\em CMS Conf. Proc.}, pages 343--357. Amer. Math. Soc.,
  Providence, R.I., 1982.

\bibitem[Pop99]{popescu:UHLandHUL}
C{\u{a}}lin Popescu.
\newblock On {UHL} and {HUL}.
\newblock {\em Bull. Belg. Math. Soc. Simon Stevin}, 6(2):219--235, 1999.

\bibitem[Pop00]{popescu:2cones}
C{\u{a}}lin Popescu.
\newblock Characteristic zero loop space homology of two-cones.
\newblock {\em Bull. London Math. Soc.}, 32(5):600--608, 2000.

\bibitem[Qui69]{quillen:rht}
Daniel Quillen.
\newblock Rational homotopy theory.
\newblock {\em Ann. of Math. (2)}, 90:205--295, 1969.

\bibitem[Rut71]{rutter:ganea}
John~W. Rutter.
\newblock On a theorem of {T}. {G}anea.
\newblock {\em J. London Math. Soc. (2)}, 3:190--192, 1971.

\bibitem[Sam53]{samelson:products}
Hans Samelson.
\newblock A connection between the {W}hitehead and the {P}ontryagin product.
\newblock {\em Amer. J. Math.}, 75:744--752, 1953.

\bibitem[Sco]{scott:tfmmPreprint}
Jonathan~A. Scott.
\newblock {A Torsion-Free Milnor-Moore Theorem}.
\newblock arXiv:math.AT/ 0103223.

\bibitem[Sel83]{selick:mooreSerre}
Paul Selick.
\newblock On conjectures of {M}oore and {S}erre in the case of torsion-free
  suspensions.
\newblock {\em Math. Proc. Cambridge Philos. Soc.}, 94(1):53--60, 1983.

\bibitem[Sel88]{selick:mooreConj}
Paul Selick.
\newblock Moore conjectures.
\newblock In {\em Algebraic topology---rational homotopy (Louvain-la-Neuve,
  1986)}, volume 1318 of {\em Lecture Notes in Math.}, pages 219--227.
  Springer, Berlin, 1988.

\bibitem[Sel97]{selick:book}
Paul Selick.
\newblock {\em Introduction to homotopy theory}, volume~9 of {\em Fields
  Institute Monographs}.
\newblock American Mathematical Society, Providence, RI, 1997.

\bibitem[{\v{S}}ir53]{shirshov:subalgebras}
A.~I. {\v{S}}ir{\v{s}}ov.
\newblock Subalgebras of free {L}ie algebras.
\newblock {\em Mat. Sbornik N.S.}, 33(75):441--452, 1953.

\bibitem[Spa93]{spanier:book}
Edwin~H. Spanier.
\newblock {\em Algebraic topology}.
\newblock Springer-Verlag, New York, 1993.
\newblock Corrected reprint of the 1966 original.

\bibitem[Whi39]{whitehead:simplicialSpaces}
J.~H.~C. Whitehead.
\newblock Simplicial spaces, nuclei and {$m$}-groups.
\newblock {\em Proc. London Math. Soc. (2)}, 45:243--327, 1939.

\bibitem[Whi41]{whitehead:addingRelations}
J.~H.~C. Whitehead.
\newblock On adding relations to homotopy groups.
\newblock {\em Ann. of Math. (2)}, 42:409--428, 1941.

\bibitem[Whi78]{whitehead:book}
George~W. Whitehead.
\newblock {\em Elements of homotopy theory}, volume~61 of {\em Graduate Texts
  in Mathematics}.
\newblock Springer-Verlag, New York, 1978.

\end{thebibliography}

\end{document}